\documentclass[12pt]{amsart}
\usepackage[usenames,dvipsnames]{xcolor}
\usepackage{amsmath,amssymb,amsthm,graphicx,mathrsfs,url,latexsym,enumerate}
\usepackage{a4wide}
\usepackage{ifthen}
\usepackage{verbatim}
\usepackage{fancyhdr}
\usepackage{url}
\usepackage{bbm}
\usepackage[sans]{dsfont}
\usepackage{transparent}
\usepackage{amsbsy}
\usepackage[colorlinks=true,linkcolor=Red,citecolor=Green]{hyperref}
\usepackage{wrapfig}
\usepackage{tikz}

\definecolor{bleu}{RGB}{27,88,145}
\definecolor{mauve}{RGB}{138,20,79}

\renewcommand{\Im}{\operatorname{Im}}
\renewcommand{\Re}{\operatorname{Re}}

\newcommand{\Id}{\operatorname{Id}}

\newcommand{\E}{\mathbb E}

\newcommand{\supp}{\operatorname{supp}}

\newcommand{\diag}{\operatorname{diag}}
\renewcommand{\span}{\operatorname{span}}
\newcommand{\C}{\mathbb C}

\newcommand{\N}{\mathbb N}
\newcommand{\R}{\mathbb R}

\newcommand{\Hess}{\operatorname{Hess}}

\newcommand{\Ran}{\operatorname{Ran}}
\newcommand{\Ker}{\operatorname{Ker}}
\newcommand{\bsigma}{\boldsymbol{\sigma}}

\newcommand{\Tr}{{\operatorname{Tr}}}
\renewcommand{\div}{\operatorname{div}}
\def\<{\langle}
\def\>{\rangle}
\newcommand{\bp}{{\it Proof. }}
\newcommand{\ep}{\hfill $\square$}
\newcommand{\be}{\begin{equation}}
\newcommand{\ee}{\end{equation}}
\newcommand{\bes}{\begin{equation*}}
\newcommand{\ees}{\end{equation*}}
\renewcommand{\d}{\operatorname{d}}

\numberwithin{equation}{section}
\numberwithin{figure}{section}

\def\m{\mathbf{m}}
\def\s{\mathbf{s}}
\def\j{\mathbf{j}}

\def\um{\underline{\m}}

\everymath{\displaystyle}
\DeclareMathOperator{\p}{\partial}
\DeclareMathOperator{\eps}{\varepsilon}
\DeclareMathOperator{\phii}{\varphi}

\newcommand{\dv}{\ dv}
\newcommand{\vvert}[1]{\left\lVert#1\right\rVert}
\newcommand{\lb}{\operatorname{[\negthinspace [}}
\newcommand{\rb}{\operatorname{]\negthinspace ]}}

\newtheorem{theorem}{Theorem}
\newtheorem{defin}{Definition}[section]
\newtheorem{corollary}[defin]{Corollary}
\newtheorem{lemma}[defin]{Lemma}

\newtheorem{proposition}[defin]{Proposition}
\newtheorem{remark}[defin]{Remark}

\newtheorem{assumption}{Assumption}

\def\bbb{{\mathcal B}}\def\ccc{{\mathcal C}}\def\ddd{{\mathcal D}}
 
 \def\lll{{\mathcal L}}
\def\mmm{{\mathcal M}} \def\ppp{{\mathcal P}}

\def\uuu{{\mathcal U}}\def\vvv{{\mathcal V}}\def\www{{\mathcal W}}

\def\Cr{{\mathscr C}}\def\Dr{{\mathscr D}}
\def\Er{{\mathscr E}}

\def\Mr{{\mathscr M}}\def\Pr{{\mathscr P}}

\begin{document}

\title{Hypocoercivity and metastability of degenerate KFP equations at low temperature}
\maketitle

\begin{center}\textsc{Loïs Delande\footnote{CERMICS, École des Ponts Champs-sur-Marne, France}\footnote{MATHERIALS, Inria Paris, Paris, France }}\end{center}

\begin{abstract}
    We consider Kramers-Fokker-Planck operators with general degenerate coefficients. We prove semiclassical hypocoercivity estimates for a large class of such operators. Then, we manage to prove Eyring-Kramers formulas for the bottom of the spectrum of some particular degenerate operators in the semiclassical regime, and quantify the spectral gap separating these eigenvalues from the rest of the spectrum. The main ingredient is the construction of sharp Gaussian quasimodes through an adaptation of the WKB method.
\end{abstract}

\section{Introduction}

\subsection{Motivations}

When considering a cloud of particles in dimension $d$, one is often led to study kinetic equations, or stochastic processes. An example of such process is given by the Langevin dynamics:
\be\label{eq:SDE0}
\left\{
\begin{aligned}
&dx_t = v_tdt,\\
&dv_t = -\partial_xV(x_t)dt - 2v_tdt + \sqrt{2h}dB_t,
\end{aligned}
\right.
\ee
where $(x_t,v_t)$ denotes the position and velocity of the particles at time $t$, $V:\R^d\to\R$ is a smooth potential corresponding to an energetic field constraining the studied particles, $h>0$ is a semiclassical parameter (typically proportional to the temperature of the system), and $B_t$ is a Brownian motion in $\R^d$ representing random forces.

To obtain results on the long-time behavior of the solution to a stochastic differential equation (SDE), one can for example study its generator. For a general SDE
\be\label{eq:SDEgeneral}
dX_t = b(X_t)dt + \sigma(X_t)dB_t,
\ee
with drift $b:\R^d\to\R^d$ and diffusion matrix $\sigma:\R^d\to\Mr_d(\R)$, it is defined as
\bes
\lll = \frac12\sum_{i,j=1}^da_{i,j}(x)\partial_i\partial_j + \sum_{i=1}^db_i(x)\partial_i,
\ees
acting on functions on $\R^d$, denoting $(a_{i,j})_{i,j} = \sigma\sigma^T$. The study of $\lll$ is important, because under mild hypotheses, we have that the solution to the partial differential equation (PDE)
\be\label{eq:FeynmanKac}
\left\{
\begin{aligned}
&\partial_tu+\lll u = 0,\\
&u_{|t=0}=u_0,
\end{aligned}
\right.
\ee
is given by $u(t,x) = \E(u_0(X_t)\ |\ X_0=x)$, with $(X_t)_t$ a solution to \eqref{eq:SDEgeneral}.

Therefore, the generator of \eqref{eq:SDE0} is
\bes
\lll_{KFP} = v\cdot \partial_x - \partial_x V\cdot \partial_v + h\Delta_v - 2v\cdot\partial_v,
\ees
where $\Delta_v$ denotes the Laplacian in the $v$ variables only. Instead of working with $\lll_{KFP}$, it is very convenient to work with a conjugation of this operator. In the following, we will denote $Q^*$ the (formal) adjoint of any operator $Q$. Observing that (formally) $\lll_{KFP}1=0$ and $\lll_{KFP}^*(e^{-2(V(x) + \frac{|v|^2}{2})/h})=0$, we therefore introduce
\bes
P_{KFP} = -e^{(V(x) + \frac{|v|^2}{2})/h}\circ h\lll_{KFP}^*\circ e^{-(V(x) + \frac{|v|^2}{2})/h}.
\ees
This operator is called the (semiclassical) Kramers-Fokker-Planck (KFP) operator and takes the form
\be\label{eq:standardKFP}
P_{KFP} = v\cdot h\partial_x - \partial_x V\cdot h\partial_v - h^2\Delta_v + |v|^2 - hd.
\ee
It has the nice form of the sum of a (formally) skew-adjoint operator and a (formally) self-adjoint one. The skew part being a transport term, and the self-adjoint part being a harmonic oscillator in velocity. Moreover, we observe that (formally)
\bes
P_{KFP}(e^{-(V(x) + \frac{|v|^2}{2})/h}) = P_{KFP}^*(e^{-(V(x) + \frac{|v|^2}{2})/h}) = 0.
\ees

The operator $P_{KFP}$ has been extensively studied in the literature, in particular when trying to compare it with the Witten Laplacian associated with $V$, that is
\bes
\Delta_V = -h^2\Delta + |\nabla V|^2 - h\Delta V.
\ees
A general theory regrouping the study of $P_{KFP}$ and $\Delta_V$ (and their extensions to forms) which gives a profound motivation to the link between the two can be found in the works of J.-M. Bismut and his hypoelliptic Laplacian, we can mention for example \cite{Bi05,BiLe05} for an overview of these concepts. We can also mention the book \cite{HeNi06} which contains a deep study of both these operators.

One of the pioneer work studying $P_{KFP}$ is \cite{HeNi04} in which the authors determined estimates regarding the return to equilibrium for the semigroup $e^{-tP_{KFP}}$ for certain potentials that are homogeneous near infinity. This result has then been generalized for example in \cite{Li09,Be22}. In \cite{HeSjSt05} the authors derived a rough localization of the small eigenvalues of $P_{KFP}$ in the semiclassical limit $h\to0$. More precise asymptotics have been obtained in \cite{HeHiSj08-2,HeHiSj11_01}. Note that these work are not restricted to the KFP operator, the authors considered a wider class of supersymmetric non-self-adjoint operators. Then, this supersymmetric assumption has been relaxed in \cite{BoLePMi22}. The goal of this paper is to obtain a similar result for a KFP operator with general coefficients that does not fit the assumptions of \cite{HeHiSj08-2,HeHiSj11_01,BoLePMi22}.

\subsection{Main results}
In this paper, we consider the following stochastic differential equation
\be\label{eq:SDE}
\left\{
\begin{aligned}
&dx_t = \alpha(x_t,v_t)dt,\\
&dv_t = \beta(x_t,v_t)dt - 4\Sigma^T\Sigma v_tdt + \sqrt{2h}dB_t,
\end{aligned}
\right.
\ee
where $\alpha$ and $\beta$ are smooth functions of both $x\in\R^d$ and $v\in\R^{d'}$, and $\Sigma\in GL_{d'}(\R)$ is a fixed invertible matrix.

Its generator is
\bes
\lll = \alpha\cdot\p_x + \beta\cdot\p_v - 4\Sigma^T\Sigma v\cdot\p_v + h\Delta_v.
\ees
Therefore we have
\bes
-\lll^* = \alpha\cdot\p_x + \beta\cdot\p_v + \div_x\alpha + \div_v\beta - (h\Delta_v + 4\Sigma^T\Sigma v\cdot\p_v + 4\Tr(\Sigma^T\Sigma)).
\ees

We notice that $(h\Delta_v + 4\Sigma^T\Sigma v\cdot\p_v + 4\Tr(\Sigma^T\Sigma))e^{-2|\Sigma v|^2/h} = 0$ and thus we make the following assumption
\begin{assumption}\label{ass.skewadj}
    There exists $V\in\ccc^\infty(\R^d,\R)$ such that denoting
    \be\label{eq:f}
    f(x,v) = V(x) + |\Sigma v|^2,
    \ee
    we have
    \bes
    \lll^*(e^{-2f/h}) = 0.
    \ees
\end{assumption}
In other words, this assumption is
\be\label{eq:assskew}
\alpha\cdot\p_xV + 2\beta\cdot\Sigma^T\Sigma v - \frac h2(\div_x\alpha+\div_v\beta) = 0.
\ee
Noticing also that $\lll1=0$ we can consider $P = -e^{f/h}\circ h\lll^*\circ e^{-f/h}$, this way we formally have $P(e^{-f/h}) = P^*(e^{-f/h}) = 0$ and we have the explicit formula $P = X + N$, with
\be\label{eq:generalP}
\left\{
\begin{aligned}
&X=\alpha(x,v)\cdot h\partial_x + \beta(x,v)\cdot h\partial_v + \frac h2(\div_x\alpha+\div_v\beta),\\
&N=-h^2\Delta_v + 4|\Sigma^T\Sigma v|^2 - 2h\Tr(\Sigma^T\Sigma).
\end{aligned}
\right.
\ee
Where we denoted $\Delta_v$ the Laplacian acting on $v$ only. We observe that we have the formal algebraic relations:
\be\label{eq:adj0}
X^* = -X,\;\;N^* = N.
\ee

Considering the following operators
\be\label{eq:dV}
\d_V = h\p_x + \p_x V
\ee
and
\be\label{eq:dSigmav}
\d_{|\Sigma v|^2} = h\p_v + 2\Sigma^T\Sigma v,
\ee
we observe that thanks to Assumption \ref{ass.skewadj}, \eqref{eq:generalP} can be rewritten
\be\label{eq:generalP2}
\left\{
\begin{aligned}
&P = X+N,\\
&X = \alpha\cdot \d_V + \beta\cdot\d_{|\Sigma v|^2},\\
&N = \d_{|\Sigma v|^2}^*\d_{|\Sigma v|^2}.
\end{aligned}
\right.
\ee

We want to obtain similar precise asymptotics of the bottom of the spectrum of $P$ to the ones obtained in \cite{HeHiSj11_01,BoLePMi22}. One of the main building block of \cite{HeHiSj08-2} (upon which \cite{HeHiSj11_01,BoLePMi22} relies) is the determination of resolvent estimates. To that extent, \cite{HeHiSj08-2} crucially needs that the symbol of the operator they consider is locally quadratic. We can easily observe that for general $\alpha,\beta$ and non-Morse potential, this cannot be true.

However, the kinetic structure of $P$ is very important and yields powerful tools to palliate \cite{HeHiSj08-2}. We can make use of the notion of hypocoercivity. This notion, was initiated in \cite{He06}, and has been further developed in \cite{Vi09,CaDoHeMiMoSc24} and references therein. Here we are in a setting which is very convenient to use a semiclassical version of \cite{DoMoSc15} and we shall follow this strategy to obtain the desired resolvent estimates.

Throughout the paper, we shall consider a potential satisfying the following assumption.

\begin{assumption}\label{ass.confin}
There exist $C>0$ and a compact  set $K\subset\R^d$ such that
\bes
V(x)\ \geq\  -C,\;\;\;
 \vert\nabla V(x)\vert \geq \frac 1 C\;\;\;\text{and}\;\;\; \Vert\Hess V(x)\Vert_\infty \leq C,
\ees
for all $x\in\R^d\setminus K$.
\end{assumption}

\begin{lemma}\label{lem:soussolaffine}
    Let $V$ satisfying Assumption \ref{ass.confin}, then there exists $b\in\R$ such that
    \bes
    \forall x\in\R^d,\ \ V(x) \geq \frac1C|x| + b
    \ees
    with $C>0$ given by Assumption \ref{ass.confin}. 
\end{lemma}
With this lemma (which proof is postponed to Subsection \ref{sec:proofsoussolaffine}), we hence have $e^{-V/h}\in L^2(\R^d)$ and thus $\Ker\Delta_V = \C e^{-V/h}$ knowing that $\d_V = e^{-V/h}\circ h\p_x\circ\ e^{V/h}$. Then, thanks to Assumption \ref{ass.confin}, $e^{-f/h}\in L^2(\R^{d+d'})$ and $N(e^{-f/h}) = 0$. Moreover Assumption \ref{ass.skewadj} is equivalent to $X(e^{-f/h}) = 0$ and together with the last assumption, it ensures that $e^{-f/h}\in D(P)$, thus, the formal equality $P(e^{-f/h}) = P^*(e^{-f/h}) = 0$ is now true in $L^2(\R^{d+d'})$.

We also consider the following assumption.
\begin{assumption}\label{ass.accretive}
    There exist $C>0$ and a compact set $K\subset\R^{d+d'}$ such that outside $K$,
    \bes
    |\div_x \alpha + \div_v\beta|\leq Cf.
    \ees
\end{assumption}
It is essential to prove the following proposition.
\begin{proposition}\label{prop:existence}
    Suppose Assumptions \ref{ass.skewadj}, \ref{ass.confin}, and \ref{ass.accretive} hold. For any initial condition, the Cauchy problem associated with \eqref{eq:SDE} admits a solution. It is almost surely unique and determines an absolutely continuous stochastic process for all time $t\geq0$.
\end{proposition}
The goal of this work is to study the long time behavior of solutions to \eqref{eq:SDE}. Therefore it is crucial to have a global existence result for this process. We postpone the proof of this proposition to Subsection \ref{ssec:existence}.

In the following we consider hypoelliptic operators. More precisely, we shall assume they satisfy the standard hypoellipticity theorem for Hörmander operators \cite[Theorem 1.1]{Ho67}
\begin{assumption}\label{ass.Hormander}
    The operators $Y_{i_1},\ [Y_{i_1},Y_{i_3}],\ [Y_{i_4},[Y_{i_5},[Y_{i_6},\ldots,Y_{i_j}]]]\ldots$ where the $Y_{i_k}$ are in $\{X,\p_{v_1},\ldots,\p_{v_{d'}}\}$ span the whole tangent space at any point of $\R^{d+d'}$.
\end{assumption}

\begin{proposition}\label{prop:accretive}
The operator $P$ initially defined on $\ccc^\infty_c(\R^{d+d'})$ admits a unique maximally accretive extension that we still denote by $(P,D(P))$.
\end{proposition}
We postpone the proof of this Proposition to Subsection \ref{ssec:accretive}.

We denote by $\uuu$ the set of critical points of $V$.
\begin{assumption}\label{ass.morse}
    For any critical point $x^*\in\uuu$, there exists a neighborhood $\vvv\ni x^*$, $(t_i^{x^*})_{1\leq i\leq d}\subset \R^{*}$, $(\nu_i^{x^*})_{1\leq i\leq d}\subset \N\setminus\{0,1\}$, a $\ccc^\infty$ change of variable $U^{x^*}$ defined on $\vvv$ such that $U^{x^*}(x^*)=x^*$, $U^{x^*}$ and $d_{x^*}U^{x^*}$ are invertible and
    \be\label{eq:morse}
    \forall x\in\vvv,\ \ V\circ U^{x^*}(x) - V(x^*) = \sum_{i=1}^dt_i^{x^*}(x_i-x_i^*)^{\nu_i^{x^*}}.
    \ee
\end{assumption}
\begin{remark}
    We notice from \eqref{eq:morse} that any $x^*\in\uuu$ is an isolated critical point, but Assumption \ref{ass.confin} implies that $\uuu\subset K$ which we recall is compact, therefore the set $\uuu$ is finite. Furthermore, if $V$ is a Morse function, $V$ satisfies Assumption \ref{ass.morse} through the Morse Lemma with $\nu_i^{x^*} = 2$ for all $i$.
\end{remark}
\begin{remark}[\cite{De25_1}]\label{rem:relax}
    We can relax \eqref{eq:morse} to
    \be\label{eq:morsealt}
    \forall x\in\vvv,\ \ V\circ U^{x^*}(x) - V(x^*) = \sum_{i=1}^dt_i^{x^*}(x_i-x_i^*)^{\nu_i^{x^*}}(1+r_i(x-x^*)),
    \ee
    with $r_i(x) = O(x)$ smooth.
    Indeed, let us consider
    \bes
    \phi : (x_1,\ldots,x_d)\to(x_1(1+r_1(x))^{1/\nu_1^{x^*}},\ldots,x_d(1+r_d(x))^{1/\nu_d^{x^*}}),
    \ees
    then $d_{0}\phi = \Id$ hence it is a $\ccc^\infty$-diffeomorphism in a neighborhood of $0$. Considering now $U^{x^*}$ that satisfies \eqref{eq:morsealt}, then $\hat U^{x^*} = U^{x^*}\circ\tau_{-x^*}\circ\phi^{-1}\circ\tau_{x^*}$ satisfies \eqref{eq:morse}, where $\tau_a(x) = x-a$ and moreover we have $d_{x^*}\hat U^{x^*} = d_{x^*}U^{x^*}$.
\end{remark}
Most of the time when the reference is clear, we will just write $U,t_i$ and $\nu_i$ instead of $U^{x^*},t_i^{x^*}$ and $\nu_i^{x^*}$ in order to lighten the notations.

We consider the partition $\uuu = \uuu^{odd}\sqcup\uuu^{even}$ where $x^*\in\uuu^{even} \iff \forall i,\ \nu_i^{x^*}\in 2\N$. Now we shall say that $x^*\in\uuu^{even}$ is of index $j\in\lb0,d\rb$ if $\sharp\{i\ |\ t_i^{x^*}<0\}=j$ and therefore $\uuu^{even}=\bigsqcup_{j=0}^d\uuu^{(j)}$ where $\uuu^{(j)}$ is the set of critical points of $V$ with even order in each direction and of index $j$, we also denote $n_0=\sharp\uuu^{(0)}$ the number of minima of $V$.

Notice that the critical points of $f$ are exactly the $(x^*,0)$, with $x^*$ a critical point of $V$. Moreover, because $|\Sigma v|^2$ is quadratic and convex, we have the same partition $\uuu = \uuu^{odd}\sqcup\uuu^{even}$ and the critical points have the same index when defined. In the following we will identify those two and use $x^*$ instead of $(x^*,0)$ where it is clear which one we are really talking about ($x^*$ will mostly be denoted either $\m$ if of index $0$ or $\s$ if of index $1$).

We define the function $\rho(v) = (C_\Sigma h)^{-\frac{d'}{4}} e^{-\frac{|\Sigma v|^2}{h}}$, where $C_\Sigma > 0$ is a normalization constant so that $\vvert{\rho}_{L^2(\R^{d'}_v)} = 1$ (so we have $C_\Sigma = \frac{\pi}{2}\big(\det\Sigma^{-1}\big)^{\frac{2}{d'}}$). We can then consider the crucial matrix that will intervene in all our work
\be\label{eq:defG}
G = \<\alpha^T\alpha\rho,\rho\>_{L^2(\R^{d'}_v)} \in \ccc^\infty(\R^d,\Mr_d(\R)).
\ee 

In the following, given two operators $A$ and $B$, we shall say that $A\lesssim B$ when for all $u$, $\<Au,u\>\lesssim\<Bu,u\>$, in the formal sense, without giving too much importance to the domains for now.

\begin{assumption}\label{ass.G}
    There exist $g_1,g_2>0$ such that for all $x\in\R^d$,
    \begin{itemize}
        \item[$i)$] $g_1(h)I_{d} \leq G(x) \leq g_2(h)I_{d}$,
        \item[$ii)$] $\forall i\in\lb1,d\rb,\ \p_{x_i}G(x) \lesssim G(x)$,
        \item[$iii)$] $\forall i,j\in\lb1,d\rb,\ \<\alpha_i^2\alpha_j^2\rho,\rho\>_{L^2(\R_v^{d'})}(x) \leq g_2(h)^2$.
    \end{itemize}
\end{assumption}

In the following we denote
\be\label{eq:nubar}
\overline \nu = \max_{i,x^*}\nu_i^{x^*}.
\ee
It allows us to introduce the function that will control the spectral gap of $P$,
\be\label{eq:g}
g(h) = \frac{h}{1 + h^{\frac{4}{\overline\nu}-2}\Big(\frac{g_2(h)}{g_1(h)}\Big)^{3} + h^{\frac1{\overline\nu}}g_1(h)^{-\frac12}}.
\ee

With these assumptions, we can state a first result concerning the rough localization of the spectrum of $P$.

\begin{theorem}\label{thm:1}
Suppose Assumptions \ref{ass.skewadj}, \ref{ass.confin}, \ref{ass.morse}, \ref{ass.G} and \ref{ass.hypocoer} hold true. There exist $h_0>0$, $c_0,c_1,c>0$, such that for all $h\in]0,h_0]$, there exists $G_h$ subspace of $L^2(\R^{d+d'})$ of dimension $n_0$ such that
$$
\forall u\in D(P)\cap G_h^\bot\ \ \Vert (P-z)u\Vert_{L^2}\geq c_1g(h)\Vert u\Vert_{L^2}
$$
for every $z\in\C$ such that $|\Re z|\leq c_0g(h)$, with $g(h)$ defined in \eqref{eq:g}. Moreover, there exists an explicit constant $c_f>0$ depending only on $f$ such that if $g(h)$ satisfies
\be\label{eq:gthm2}
g(h)\geq e^{-\frac{\tilde c}{2h}}\mbox{ for any } \tilde c<c_f,
\ee
then there exist $\lambda_\m(h)\in\C$ for all $\m\in\uuu^{(0)}$ such that $\sigma(P)\cap\{|\Re z|\leq c_0g(h)\}=\{\lambda_\m(h),\m\in\uuu^{(0)}\}$ counted with multiplicity, and for all $\m\in\uuu^{(0)},\ |\lambda_\m(h)|\leq e^{-c/h}$. Finally (still under \eqref{eq:gthm2}), for all $0<c_0'<c_1$,
\bes
\forall |z|>c_0'g(h),\mbox{ such that }|\Re z|\leq c_0g(h),\ \vvert{(P-z)^{-1}}_{L^2}\leq\frac{2}{c_0'g(h)}.
\ees
In addition, under Assumptions \ref{ass.accretive} and \ref{ass.Hormander}, we have Proposition \ref{prop:accretive} and thus we can extend all these results for $\Re z \leq -c_0g(h)$.
\end{theorem}

\begin{remark}
    Observe from \eqref{eq:g} that if $\alpha$ and $\beta$ have low degeneracy leading to $g_1,g_2\sim h$ in Assumption \ref{ass.G}, then $\vvert{(P-z)^{-1}} = O(h^{\frac{4}{\overline{\nu}}-3})$ in the region $\{|z|>c_0'g(h)\}\cap\{\Re z\leq c_0g(h)\}$. We notice that this resolvent estimate is not the same as the one for the Witten Laplacian obtained in \cite{De25_1} which we recall is of order $h^{\frac{2}{\overline{\nu}}-2}$. Although the estimates are different, when taking $\overline{\nu} = 2$, that is $V$ is Morse, we obtain an $h^{-1}$ for both the Fokker-Planck operator and the Witten Laplacian (which we know is sharp), therefore, even if our result may not be optimal, it still remains relevant.
\end{remark}

For the rough localization of Theorem \ref{thm:1}, we only need vague constructions around the critical points of $V$ in order to have Proposition \ref{prop:gapWitten}. In order to obtain sharp asymptotics, we have to refine our constructions. It requires the introduction of the topological definitions we recall in Subsection \ref{sec:labeling}.

We recall the generic assumption \eqref{ass.gener} which is useful in order to lighten the result and the proof.
\be\label{ass.gener}\tag{Gener}
\begin{array}{l}
    (\ast)\mbox{ for any }\m\in\uuu^{(0)},\m \mbox{ is the unique global minimum of } V_{|E(\m)}\\
    (\ast)\mbox{ for all }\m\neq\m'\in\uuu^{(0)},\j(\m)\cap\j(\m')=\emptyset.
\end{array}
\ee
In particular, \eqref{ass.gener} implies that $V$ uniquely attains its global minimum at $\um$. This assumption allows us to avoid some heavy constructions regarding the set $\uuu$ and lighten the definition \ref{def:quasimodes} of the quasimodes. But it seems that its not a true obstruction and that we can pursue the computations without this assumption as described in \cite[Section 6]{BoLePMi22}, \cite{No24}, in the spirit of \cite{Mi19}.

This assumption comes from \cite[(1.7)]{BoGaKl05_01} and \cite[Assumption 3.8]{HeKlNi04_01}. Then it changed to become \cite[Hypothesis 5.1]{HeHiSj11_01} when finally reaching the form of \cite[Assumption 4]{LePMi20}.

One can show that \eqref{ass.gener} is weaker than \cite{BoGaKl05_01}'s, \cite{HeKlNi04_01}'s and \cite{HeHiSj11_01}'s assumptions. More precisely, they supposed that $\j(\m)$ is a singleton while here we have no restriction on the size of $\j(\m)$.

The hypotheses of Theorem \ref{thm:1} hold for very general coefficients $\alpha,\beta$ and potential $V$. Then, the objective is to refine this rough localization of the eigenvalues of $P$. To that end, we follow the method of \cite[Section 3]{BoLePMi22} which is an adaptation of the WKB method. However, we were not able to make it work under the hypotheses of Theorem \ref{thm:1}. In Section \ref{sec:geoloc}, we study several situations in which we are able to give results. These are summarized under Proposition \ref{prop:ell} which was the target of our study of the WKB method. To ease the study of the aforementioned situations, we shall assume that the potential $V$ is a Morse function.
\begin{assumption}\label{ass:vraimentMorse}
    The potential $V$ is a Morse function.
\end{assumption}
We can observe what behavior arises when this assumption is not satisfied for the Witten Laplacian in \cite{De25_1}. In addition, we also impose a stronger assumption than \eqref{eq:gthm2}.
\begin{assumption}\label{ass.polyg}
    The function $g$ defined in \eqref{eq:g} is no worse than polynomial. That is, there exist $h_0,c>0$ such that for all $h\in(0,h_0]$,
    \bes
    g(h)\geq h^c.
    \ees
\end{assumption}
From \eqref{eq:g}, $g(h) = O(h)$, therefore we in fact have $c\geq1$.

\begin{theorem}\label{thm:2}
    Suppose the hypotheses of Theorem \ref{thm:1}, Assumptions \ref{ass:vraimentMorse}, \ref{ass.polyg}, and \eqref{ass.gener} hold true. Suppose moreover that Proposition \ref{prop:ell} apply. There exist $h_0,\gamma > 0$ such that for all $h\in\ ]0,h_0]$, one has $\lambda(\um,h) = 0$ and for all $\m\neq\um$, $\lambda(\m,h)$ satisfies the following Eyring–Kramers type formula
    \bes
    \lambda(\m,h) = v(\m)h^{\mu(\m)}e^{-2S(\m)/h}(1+O(h^\gamma)),
    \ees
    where $v(\m)$ and $\mu(\m)$ are defined in \eqref{eq:defvmu} and depend explicitly on $V$, and $S$ is the standard height function defined by \eqref{eq:defS}.
\end{theorem}

As mentioned in the paragraph after \eqref{ass.gener}, it seems that this assumption is not required in order to obtain sharp results although we did no computation without it. For Assumption \ref{ass.polyg}, we observe that it is very commonly satisfied and all examples in this article satisfy it.

In \cite{De25_1}, we obtained a very similar result for the Witten Laplacian, but we can mention some differences on the prefactor $v(\m)$ and $\mu(\m)$. The difference between the prefactor $v(\m)h^{\mu(\m)}$ of the small eigenvalues of the standard KFP operator \eqref{eq:standardKFP} and the ones of the Witten Laplacian, both with Morse potential lies in the negative eigenvalue of some non-degenerate matrix (see for example \cite{BoLePMi22} for a formula in both cases). Here it is a little bit more subtle. The degeneracies in $\alpha$ and $\beta$ induce a degeneracy in the resolution of the equations obtained by the WKB method leading to another factor when doing the Laplace method. The exponent $\mu(\m)$ in \cite{De25_1} was completely determined by the degeneracy of the potential $V$. While here, we proved Theorem \ref{thm:2} only for Morse potentials and yet, we can have various powers of $h$. This is due to the fact that degeneracies on $\alpha$ and $\beta$ affect $\mu(\m)$.

Note that the Arrhenius law $\lim_{h\to0}h\ln\lambda(\m,h) = -2S(\m)$ remains the same as usual. It is coherent with \cite[Corollary 1.2.4]{NiSaWh24} in which the authors proved this law for the standard KFP operator with very general potential. It is reasonable to assume that this very robust law (continuous in $V$ for the $C^0$ topology) would hold for a wide class of coefficients $\alpha$ and $\beta$.

\subsection{Metastability}
The results presented here, proved in Subsection \ref{ssec:proofcor}, are adaptations of \cite[Corollary 1.5, 1.6]{BoLePMi22} to our settings. In generality, a result of the form of Theorem \ref{thm:2} allows to write such corollaries without much assumptions.

From Proposition \ref{prop:accretive}, $P$ is maximally accretive, therefore, for all $u_0\in L^2(\R^d)$, the following Cauchy problem
\be\label{eq:evolutionP}
\left\{
\begin{aligned}
    &h\partial_tu + Pu = 0,\\
    &u_{|t=0} = u_0,
\end{aligned}
\right.
\ee
admits a unique solution $u\in\ccc^0([0,+\infty),L^2(\R^d))\cap\ccc^1((0,+\infty),L^2(\R^d))$ denoted $u(t) = e^{-tP/h}u_0$. Theorem \ref{thm:2} gives the following result on the long time behavior of that solution.

\begin{corollary}\label{cor:return}
    In the setting of Theorem \ref{thm:2}, there exist $C,\varepsilon>0$ such that, for all $u_0\in L^2(\R^d)$ and $h$ small enough, there exists $(u_{\m,n})_{\m,n}\subset\C$ such that the solution $u(t)$ of \eqref{eq:evolutionP} satisfies
    \be\label{eq:return1}
    \forall t\geq0,\ \ \vvert{u(t) - \sum_{\m\in\uuu^{(0)}}\sum_{n=0}^{n_0-1}u_{\m,n}t^ne^{-\lambda(\m,h)t/h}} \leq Ce^{-\varepsilon h^{c-1}t}\vvert{u_0},
    \ee
    with $c\geq1$ given by Assumption \ref{ass.polyg}. Moreover, there exists $C>0$ such that, for all $u_0\in L^2(\R^d)$ and $h$ small enough, the solution $u(t)$ of \eqref{eq:evolutionP} satisfies
    \be\label{eq:return2}
    \forall t\geq 0,\ \ \vvert{u(t) - \frac{\<e^{-f/h},u_0\>}{\vvert{e^{-f/h}}^2}e^{-f/h}} \leq Ce^{-t\underset{\m\neq\um}{\min}\Re(\lambda(\m,h))(1-Ch)/h}\vvert{u_0}.
    \ee
\end{corollary}

Another way to write \eqref{eq:return1} is
\bes
u(t) = e^{-tP/h}\Pi_\Cr u_0 + O(e^{-\varepsilon t})\vvert{u_0},
\ees
while \eqref{eq:return2} is
\bes
u(t) = e^{-tP/h}\Pi_0 u_0 + O(e^{-t\underset{\m\neq\um}{\min}\Re(\lambda(\m,h))(1-Ch)/h})\vvert{u_0},
\ees
with the $O$ being uniform in $t$ and $h$. Here, $\Pi_0$ denotes the orthogonal projector on the kernel of $P$, and $\Pi_\Cr$ denotes the spectral projector of $P$ associated with its $n_0$ exponentially small eigenvalues (recalling $n_0$ is the number of minima of $V$). It is defined as
\bes
\Pi_\Cr = \frac{1}{2i\pi}\int_{\Cr}(z-P)^{-1}dz,
\ees
where $\Cr = \partial D(0,\frac{c_0}2g(h))$, with $c_0$ and $g(h)$ given by Theorem \ref{thm:1}.

Furthermore, we can describe the metastable behavior of the solutions of \eqref{eq:evolutionP}.

\begin{corollary}\label{cor:times}
    In the setting of Theorem \ref{thm:2}, let $S_1\leq\cdots\leq S_{p+1}=+\infty$ denote the non-decreasing sequence of the $S(\m)$'s defined in \eqref{eq:defS} such that if $S_k = S_{k+1}$, then $\mu_k < \mu_{k+1}$ and let $\Pi^\leq_k$ be the spectral projector of $P$ associated with its eigenvalues of modulus of order less than $h^{\mu_k}e^{-2S_k/h}$. For two positive functions $t_\pm(h)$ such that $t_-(h) = O(h^\infty)$ and $t_+^{-1}(h) = o(h^{c-1}|\ln h|^{-1})$, we define the times
    \bes
    t_0^+ = t_+(h)\ \text{ and }\ \forall1\leq k\leq p+1,\ t^\pm_k = t_\pm(h)h^{1-\mu_k}e^{2S_k/h}
    \ees
    (in particular $t^-_{p+1} = +\infty$). Then, for every $h$ small enough, the solution $u(t)$ of \eqref{eq:evolutionP} satisfies
    \bes
    \forall t^+_{k-1}\leq t\leq t^-_k,\ u(t) = \Pi^\leq_ku_0 + O(h^\infty)\vvert{u_0},
    \ees
    uniformly with respect to $t$, $1\leq k\leq p+1$, and $u_0\in L^2(\R^d)$.
\end{corollary}

In other words, $e^{-tP/h}$ is approximately constant equal to $\Pi^\leq_k$ on the time interval $[t^+_{k-1},t^-_k]$, with fast transition around the times $t_k = h^{\mu_k}e^{2S_k/h}\in(t^-_k,t^+_k)$. In this corollary, one can take $t_-(h) = e^{-\delta/h}$ for some $\delta>0$ and $t_+(h) = \frac{|\ln h|^2}{h^{c-1}}$.

We can link this with the stochastic process solving \eqref{eq:SDE}. Under mild hypotheses, one can show that the probability density $\rho(t,\cdot)$ of the process $(X_t)_t$ solution of \eqref{eq:SDE} is solution to the problem
\be\label{eq:probadensity}
\partial_t \rho = \lll^*\rho
\ee
Therefore, recalling that $P = -e^{f/h}\circ(h\lll^*)\circ e^{-f/h}$, $u$ is solution to \eqref{eq:evolutionP} if and only if $e^{-f/h}u$ is solution to \eqref{eq:probadensity} (with adapted initial conditions). We observe that Corollary \ref{cor:return} gives
\bes
\forall t\geq 0,\ \ \vvert{\rho(t) - \frac{\<e^{-f/h},u_0\>}{\vvert{e^{-f/h}}^2}e^{-2f/h}}_{TV} \leq Ce^{-t\underset{\m\neq\um}{\min}\Re(\lambda(\m,h))(1-Ch)/h}\vvert{e^{-f/h}}\vvert{u_0},
\ees
using that for absolutely continuous measure $\mu,\nu$, $\vvert{\mu-\nu}_{TV} = \vvert{\mu-\nu}_1$ and the Cauchy-Schwarz inequality.

\subsection*{Acknowledgements}
The author is grateful to Laurent Michel for his advice through this work and to Jean-François Bony for helpful discussions. This work is supported by the ANR project QuAMProcs 19-CE40-0010-01.

The rest of the paper is organized as follows. In the next section, we develop the hypocoercive estimates in order to have a rough localization of the eigenvalues of $P$ as well as resolvent estimates resulting in Theorem \ref{thm:1}. In Section \ref{sec:examples}, we give a short list of examples of coefficients $\alpha$ and $\beta$ for which Theorem \ref{thm:1} apply.

Then, we want to obtain precise Eyring-Kramers laws for our operator. We did not manage to obtain such results in broad generality, but we could prove these in several prescribed situations. Section \ref{sec:geoloc} focuses on the local constructions of the WKB method in the spirit of \cite{BoLePMi22}, while the purpose of Section \ref{sec:global} is to glue these local constructions to obtain globally defined cutoffs. This leads to the proof of Theorem \ref{thm:2}.

\section{Hypocoercive estimates}\label{sec:Hypocoer}

Let $\chi_\m$, $\m\in \uuu^{(0)}$ be some cutoffs in $\ccc_c^\infty(\R^d)$ such that $\chi_\m$ is supported in $B(\m,r)$ for some $r>0$ to be chosen small enough and $\chi_\m=1$ near $\m$. We then consider
\bes
V_\m(x) = \chi_\m(x) e^{-(V(x) - V(\m))/h}
\ees
and their space
\bes
E_h = \span\{V_\m,\ \m\in\uuu^{(0)}\}.
\ees
For $r$ small enough, the $V_\m$ have disjoint support hence $\dim E_h = n_0$.

\begin{proposition}\label{prop:gapWitten}\cite[Proposition 2.1]{De25_1}
    Recalling $\overline{\nu}=\max_{i,x^*}\nu_i^{x^*}$, we have
    \bes
    \exists C>0,\ \forall u\in D(\Delta_V)\cap E_h^\bot\ \ \<\Delta_Vu,u\> \geq Ch^{2-\frac2{\overline{\nu}}}\vvert{u}^2.
    \ees
\end{proposition}

As the critical points of $f$ are the $(x^*,0)$ for $x^*\in\uuu$, with the same index (when it is defined), we will identify those two and use $x^*$ instead of $(x^*,0)$ where it is clear which one we are really talking about ($x^*$ will mostly be denoted either $\m$ if of index $0$ or $\s$ if of index $1$). We also define the global quasimodes
\bes
f_\m(x,v)=\chi_\m(x) e^{-(f(x,v)-f(\m))/h} = V_\m(x)e^{-\frac{|\Sigma v|^2}{h}},
\ees
\bes
F_h=\span\{f_\m,\;\m\in\uuu^{(0)}\}.
\ees
Both $E_h$ and $F_h$ have dimension $n_0$ but we must notice that $E_h\subset L^2(\R^d)$ while $F_h\subset L^2(\R^{d+d'})$.

We recall the function $\rho(v)=(C_\Sigma h)^{-\frac {d'}4}e^{-\frac{|\Sigma v|^2}{h}}$ and we introduce the projector onto the kernel of $N$ defined on $L^2(\R^{d+d'})$ by
\bes
\Pi u(x,v)=\int_{\R^{d'}}u(x,v')\rho(v')dv' \rho(v)=u_\rho(x)\rho(v),
\ees
where we denoted
\be\label{eq:rho}
u_\rho = \< u, \rho\>_{L^2(\R_v^{d'})}.
\ee
We observe that $\Pi (E_h\otimes L^2(\R^{d'})) = F_h$.

\begin{lemma}\label{lem:XPi}
    Under Assumption \ref{ass.skewadj},
    \be\label{eq:XPi}
    X\Pi = \alpha\cdot\d_V\Pi
    \ee
    and hence recalling $G = \<\alpha\alpha^T\rho,\rho\>_{L^2(\R_v^{d'})} \in \ccc^\infty(\R^d,\Mr_d(\R))$, we obtain
    \be\label{eq:XPistarXPi}
    (X\Pi)^*(X\Pi) = \d_V^*G\d_V\Pi.
    \ee
\end{lemma}

\bp Using \eqref{eq:generalP2}, and having that $\Pi$ is a projector on the kernel of $\d_{|\Sigma v|^2}$, we immediately get that $X\Pi = \alpha\cdot\d_V\Pi$. We then obtain
\bes
\begin{aligned}
    (X\Pi)^*(X\Pi) &= \Pi(\alpha\cdot\d_V)^*(\alpha\cdot\d_V)\Pi = \d_V^*\Pi\alpha\alpha^T\Pi\d_V\\
    \Pi\alpha\alpha^T\Pi u &= \Pi\alpha\alpha^Tu_\rho(x)\rho = u_\rho(x)\<\alpha\alpha^T\rho,\rho\>_{L^2(\R_v^{d'})}\rho = G\Pi u
\end{aligned}
\ees
and hence $(X\Pi)^*(X\Pi) = \d_V^*G\d_V\Pi$.

\ep

Under Assumption \ref{ass.G}, we have a straightforward corollary of Proposition \ref{prop:gapWitten},

\begin{corollary}\label{cor:lowbound}
    $\exists C>0,\ \forall u\in D(\Delta_V)\cap F_h^\bot\ \<\d_V^*G\d_Vu,u\> \geq Ch^{2-\frac2{\overline{\nu}}}g_1(h)\vvert{u}^2$.
\end{corollary}

We now define the following auxiliary operator
\be\label{eq:A}
A = (h^{2-\frac2{\overline{\nu}}}g_1(h) + (X\Pi)^*(X\Pi))^{-1}(X\Pi)^* = (h^{2-\frac2{\overline{\nu}}}g_1(h) + \d_V^*G\d_V)^{-1}(X\Pi)^*.
\ee
This auxiliary operator is introduced in \cite{DoMoSc15} and used in \cite{LeSaSt20} in order to ease the calculus in the proof of Theorem \ref{thm:1}. This kind of method to compute hypocoercivity was mainly introduced and used at first in \cite{Vi09}, \cite{HeNi04} and \cite{He06}.
\begin{lemma}\label{lem:Aborne}
The operator $A$ is bounded on $L^2(\R^{d+d'})$, it satisfies $A = \Pi A$ and one has the estimate
\bes
\Vert A\Vert_{L^2}\leq O(h^{\frac1{\overline{\nu}}-1}g_1(h)^{-\frac12})
\ees
\end{lemma}
\bp The bound is easily seen using Lemmas \ref{lem:1} and \ref{lem:2}.

\ep

The following assumption will help us prove some bounds on $A$ for the hypocoercivity result.

\begin{assumption}\label{ass.hypocoer}
    We consider coefficients $\alpha$ and $\beta$ such that
    \begin{itemize}
        \item[$i)$] $\forall (x,v)\in\R^{d+d'}\ \int \alpha(x,\sqrt{h}v)e^{-2|\Sigma v|^2}\dv = 0$.
        \item[$ii)$] For all $q\in\{hJ_x\alpha\alpha,hJ_v\alpha\beta,hJ_v\alpha\Sigma^T\Sigma v,h^2\Delta_v\alpha\}$, for all $i,j\in\lb1,d\rb$
        \bes
        \Pi q_iq_j\Pi \lesssim g_2(h)^2(|\nabla V|^2+h)\Pi,
        \ees
        where we denote $J_x\alpha$ the Jacobian matrix of $\alpha$ with respect to the variable $x$ (and likewise for $J_v\alpha$) and $\Delta_v\alpha$ the vector $(\Delta_v\alpha_i)_i$.
    \end{itemize}
\end{assumption}

\begin{remark}
    Let us notice that $i)$ implies that $\Pi X\Pi = \Pi\alpha\d_V\Pi = 0$ and thus $A = A(1-\Pi)$.
\end{remark}


This leads to the intermediate Lemma

\begin{lemma}\label{lem:bounds}
    Under Assumption \ref{ass.hypocoer}, there exists $C,h_0>0$ such that for all $h\in]0,h_0]$, for all $u\in L^2(\R^{d+d'})$, one has
    \be\label{eq:maj1}
    |\<AX (1-\Pi)u,u\>|\leq Ch^{\frac{2}{\overline\nu}-1}\Big(\frac{g_2(h)}{g_1(h)}\Big)^{\frac32}\Vert \Pi u\Vert\, \Vert(1- \Pi) u\Vert,
    \ee
    \be\label{eq:maj2}
    |\<AN u,u\>|\leq Ch^{\frac{2}{\overline\nu}-1}\Big(\frac{g_2(h)}{g_1(h)}\Big)^{\frac32}\Vert \Pi u\Vert\, \Vert(1- \Pi) u\Vert,
    \ee
    \be\label{eq:maj3}
    |\<X u,A u\>|\leq C\Vert(1- \Pi) u\Vert^2.
    \ee
\end{lemma}
\bp
Within this proof, $C$ will denote a positive constant that may only depends on the dimension $d$ and $\Sigma$ and can change from line to line. Let us denote
\be\label{eq:defR}
R = (h^{2-\frac2{\overline{\nu}}}g_1(h) + \d_V^*G\d_V)^{-1},
\ee
this will makes the computations clearer. This way, $A = R(X\Pi)^*$.

Let us start with the proof of \eqref{eq:maj1}. Since $A=\Pi A$, by the Cauchy-Schwarz inequality it is sufficient to show that the operator $AX$ (or equivalently its adjoint) is bounded on $L^2$. But $X^*A^* = -X^2\Pi R$ with
\bes
X^2\Pi = X(\alpha\cdot\d_V)\Pi = (\alpha\cdot\d_V)^2\Pi + \beta\cdot\d_{|\Sigma v|^2}\alpha\cdot\d_V\Pi
\ees
thanks to Lemma \ref{lem:XPi} and \eqref{eq:generalP2}. Using that $\d_{|\Sigma v|^2}\Pi = 0$,
\be\label{eq:X2Pi}
X^2\Pi = (\alpha\cdot\d_V)^2\Pi + h\beta\cdot\p_v(\alpha\cdot\d_V)\Pi = (\alpha\cdot\d_V)^2\Pi + hJ_v\alpha\beta\cdot\d_V\Pi.
\ee
Moreover, denoting $\d_{V,i} = h\p_{x_i} + \p_{x_i}V$, we can write
\be\label{eq:alphadv2}
\begin{aligned}
    (\alpha\cdot\d_V)^2 &= \sum_{i,j}\alpha_i\d_{V,i}\alpha_j\d_{V,j}\\
    &= \sum_{i,j}\alpha_i\alpha_j\d_{V,i}\d_{V,j} + \sum_{i,j}\alpha_i[\d_{V,i},\alpha_j]\d_{V,j}\\
    &= -\sum_{i,j}\alpha_i\alpha_j\d_{V,i}^*\d_{V,j} + 2\sum_{i,j}\alpha_i\alpha_j\p_{x_i}V\d_{V,j} + hJ_x\alpha\alpha\cdot\d_V.
\end{aligned}
\ee

Using Assumption \ref{ass.G} $iii)$, for $i,j\in\lb1,d\rb$ and $u\in L^2(\R^{d+d'})$, we have
\bes
\begin{aligned}
    \vvert{\alpha_i\alpha_j\d_{V,i}^*\d_{V,j}\Pi Ru}^2 &= \<\Pi \alpha_i^2\alpha_j^2\Pi\d_{V,i}^*\d_{V,j}Ru,\d_{V,i}^*\d_{V,j}Ru\>\\
    &\leq Cg_2(h)^2\vvert{\d_{V,i}^*\d_{V,j}Ru}^2\\
    &\leq C\big(h^{\frac{2}{\overline{\nu}} - 1}g_1(h)^{-3/2}g_2(h)^{3/2}\big)^2
\end{aligned}
\ees
with Lemma \ref{lem:4}. Using the same arguments, we obtain
\bes
\vvert{\alpha_i\alpha_j\p_{x_i}V\d_{V,j}\Pi Ru} \lesssim g_2(h)\vvert{\p_{x_i}V\d_{V,j}\Pi Ru}.
\ees
Moreover, with Assumption \ref{ass.hypocoer}, we have
\bes
\vvert{hJ_x\alpha\alpha\cdot\d_V\Pi Ru} \lesssim g_2(h)(\vvert{|\nabla V|\d_V\Pi Ru} + \sqrt{h}\vvert{\d_V\Pi Ru}),
\ees
hence there just remains to control terms of the form $\p_{x_i}V\d_{V,j}R$ using Lemma \ref{lem:3}. Noticing that the non-negativity of the Laplacian implies that $|\nabla V|^2\leq \Delta_V + h\Delta V$, we have
\be\label{eq:nablaVleqDeltaV}
|\nabla V|^2\leq\Delta_V + C'h
\ee
for some $C'>0$ using Assumption \ref{ass.confin}. Thus we just have to estimate
\bes
\<\d_V^*\Delta_V\d_VRu,Ru\>.
\ees
Notice that the $\Delta_V$ appearing is actually a scalar matrix $\Delta_VI_d$ where $I_d$ denotes the identity matrix of size $d$. Consider now the relation $\Delta_V^{(1)} = \Delta_VI_d + 2h\Hess V$ where $\Delta_V^{(1)}$ denotes the Witten Laplacian acting on $1$-forms (that we identify with $\R^d$-valued functions, we refer to \cite{HeSj85_01} for more details about Witten Laplacians on $p$-forms). Using Assumption \ref{ass.confin}, and the commutation rule (\cite[(2.1)]{HeKlNi04_01})
\bes
\Delta_V^{(1)}\d_V = \d_V\Delta_V,
\ees
we obtain
\bes
\begin{aligned}
    \<\d_V^*\Delta_V\d_VRu,Ru\> &= \<\d_V^*\Delta_V^{(1)}\d_VRu,Ru\> - 2h\<\d_V^*\Hess V\d_VRu,Ru\>\\
    &\leq \<\d_V^*\d_V\Delta_VRu,Ru\> + 2Ch\vvert{\d_VRu}^2\\
    &= \vvert{\Delta_VRu}^2 + 2Ch\vvert{\d_VRu}^2.
\end{aligned}
\ees
This shows that
\be\label{eq:pVdVR}
\vvert{|\nabla V|\d_VR} \lesssim \vvert{\Delta_VR} + \sqrt{h}\vvert{\d_VR},
\ee
therefore we obtain with
\bes
\vvert{X^2\Pi R} \lesssim h^{\frac{2}{\overline{\nu}} - 1}g_1(h)^{-3/2}g_2(h)^{3/2} + g_2(h)\vvert{\Delta_V\Pi R} + \sqrt hg_2(h)\vvert{\d_V \Pi R}.
\ees
Hence \eqref{eq:maj1} using Lemmas \ref{lem:1}, \ref{lem:3} and \ref{lem:4}.

Now for \eqref{eq:maj2}, because $\Pi$ is the projection onto the kernel of $N$, we have $N = N(1-\Pi)$ and we recall that $A = \Pi A$. Hence by the Cauchy-Schwarz inequality, we only need a bound on $AN$ or $NA^*$ (since $N$ is self-adjoint). Recall that $X\Pi = \alpha\cdot\d_V\Pi$ and thus
\bes
\begin{aligned}
    NA^* &= N\alpha\cdot\d_V\Pi R = [N,\alpha\cdot\d_V]\Pi R = [N,\alpha]\cdot\d_V\Pi R\\
    &= -h^2(\Delta_v\alpha + 2J_v\alpha\p_v)\Pi\cdot\d_V R\\
    &= (-h^2\Delta_v\alpha + 4hJ_v\alpha\Sigma^T\Sigma v)\cdot\d_V\Pi R
\end{aligned}
\ees
where we recall we denoted $\Delta_v\alpha$ the vector $(\Delta_v\alpha_i)_i$. Therefore thanks to Assumption \ref{ass.hypocoer} and \eqref{eq:pVdVR},
\bes
\vvert{NA^*} \lesssim g_2(h)\big(\vvert{\Delta_V\Pi R} + \sqrt{h}\vvert{\d_V\Pi R}\big),
\ees
which proves \eqref{eq:maj2} using the same estimates as for $X^*A^*$.

And finally for \eqref{eq:maj3}, using that $A = \Pi A (1-\Pi)$ and $\Pi X \Pi = 0$ we quickly obtain
\bes
\begin{aligned}
    \<Xu,Au\> &= \<Xu,\Pi A(1-\Pi)u\> = \<\Pi X(1-\Pi)u,A(1-\Pi)u\>\\
    &\leq \vvert{(\Pi X)^*A}\vvert{(1-\Pi)u}^2\\
    &= \vvert{X\Pi R(X\Pi)^*}\vvert{(1-\Pi)u}^2
\end{aligned}
\ees
And we recognize the norm of $QQ^*$ with $Q = X\Pi R^{1/2}$ which is bounded by Lemmas \ref{lem:XPi} and \ref{lem:2}, hence the result.

\ep

\begin{proposition}\label{prop:hypocoer}
There exists $C,\delta_0,h_0>0$ such that for all $h\in]0,h_0]$, and for all $u\in D(P)\cap F_h^\bot$, one has
$$
\Re\,\<Pu,(1+\delta(h)(A+A^*))u\>\geq C\delta(h) \Vert u\Vert^2,
$$
where $\delta(h)=\delta_0 \frac{h}{1 + h^{\frac{4}{\overline\nu}-2}\Big(\frac{g_2(h)}{g_1(h)}\Big)^{3}}$.
\end{proposition}
\bp
For all $\delta>0$ and $u\in D(P)\cap F_h^\bot$, let us define 
\bes
I_\delta = \Re\,\<P u,(1+\delta(A+A^*)u\>
\ees
Using the decomposition $P = X + N$, and the skew-adjointness of $X$ coming from \eqref{eq:adj0}, one gets
\bes
I_\delta = \<N u,u\> + \delta\Re\<Pu,(A+A^*)u\>
\ees
From the spectral properties of $N$, it follows that
\be\label{eq:coer1}
I_\delta \geq h\Vert (1-\Pi)u\Vert^2 + \delta\Re\<Pu,(A+A^*)u\>.
\ee
Denoting $J = \<Pu,(A+A^*)u\>$, one has
\bes
J = \<AXu,u\> + \<ANu,u\> + \<Xu,Au\> + \<Nu,Au\>
\ees
and since $A = \Pi A$ and $\Pi N = 0$ it follows that 
\be\label{eq:J}
J = \<AX\Pi u,u\> + J'
\ee
with 
\be\label{eq:J'}
J' = \<AX(1-\Pi)u,u\> + \<ANu,u\> + \<Xu,Au\>.
\ee
Moreover, by definition of $A$ and Lemma \ref{lem:XPi}, on $F_h^\bot$
\be\label{eq:AXPi}
AX\Pi = (h^{2-\frac2{\overline{\nu}}}g_1(h) + \d_V^*G\d_V)^{-1}\d_V^*G\d_V\Pi \geq c_0\Pi
\ee
for some $c_0>0$ using functional calculus and Corollary \ref{cor:lowbound}. Therefore, combining \eqref{eq:coer1}, \eqref{eq:J}, \eqref{eq:J'}, \eqref{eq:AXPi} and Lemma \ref{lem:bounds}, we obtain
\bes
\begin{aligned}
    \forall u\in D(P)\cap F_h^\bot,\ \ I_\delta &\geq h\vvert{(1-\Pi)u}^2 + \delta c_0\vvert{\Pi u}^2 - C\delta\vvert{(1-\Pi)u}^2\\&\phantom{****} - C\delta h^{\frac{2}{\overline\nu}-1}\Big(\frac{g_2(h)}{g_1(h)}\Big)^{\frac32}\Vert \Pi u\Vert\, \Vert(1- \Pi) u\Vert,
\end{aligned}
\ees
thus with Young's inequality, we have
\bes
I_\delta \geq \Big(h-C\delta-\frac{C^2}{2c_0}\delta h^{\frac{4}{\overline\nu}-2}\Big(\frac{g_2(h)}{g_1(h)}\Big)^{3}\Big)\vvert{(1-\Pi)u}^2 + \delta\frac{c_0}{2}\vvert{\Pi u}^2.
\ees
Optimizing the right hand side by taking
\bes
\delta = \frac{2c_0h}{c_0^2 + 2c_0C + C^2h^{\frac{4}{\overline\nu}-2}\Big(\frac{g_2(h)}{g_1(h)}\Big)^{3}},
\ees
we obtain
\bes
I_\delta \geq \delta\frac{c_0}{2}\vvert{u}^2.
\ees
Noticing that there exists $C>0$ such that for $h$ small enough,
\bes
\frac1C\delta(h)\leq\frac{h}{1 + h^{\frac{4}{\overline\nu}-2}\Big(\frac{g_2(h)}{g_1(h)}\Big)^{3}}\leq C\delta(h),
\ees
we proved the proposition.

\ep

\subsection{Proof of Theorem \ref{thm:1}}
Let $z\in\C$, $u\in D(P)\cap F_h^\bot$ and $\delta(h)$ as in Proposition \ref{prop:hypocoer}. First with the Cauchy-Schwarz inequality, we have
\be\label{eq:thm1}
\Re\< (P-z)u, (1+\delta(h)(A+A^*))u \> \leq \vvert{(P-z)u}\vvert{1+\delta(h)(A+A^*)}\vvert{u},
\ee
then thanks to Proposition \ref{prop:hypocoer}, 
\bes
\Re\< (P-z)u, (1+\delta(h)(A+A^*))u \> \geq C\delta(h)\vvert{u}^2-\Re(z\< u,(1+\delta(h)(A+A^*))u\>).
\ees
Because $(1+\delta(h)(A+A^*))$ is symmetric, this leads to
\bes
\Re\< (P-z)u, (1+\delta(h)(A+A^*))u \> \geq C\delta(h)\vvert{u}^2-|\Re z|\vvert{u}^2\vvert{1+\delta(h)(A+A^*)}.
\ees
Using that $\vvert{1+\delta(h)(A+A^*)} \leq 1 + 2\delta(h)\vvert{A}$, this and \eqref{eq:thm1} lead to
\bes
\vvert{(P-z)u} \geq C\frac{\delta(h)}{1 + 2\delta(h)\vvert{A}}\vvert{u} - |\Re z|\vvert{u}.
\ees
With the expression of $\delta(h)$ and Lemma \ref{lem:Aborne}, let us recall \eqref{eq:g}, for $h>0$, $g(h)$ is defined as
\bes
g(h) = \frac{h}{1 + h^{\frac{4}{\overline\nu}-2}\Big(\frac{g_2(h)}{g_1(h)}\Big)^{3} + h^{\frac1{\overline\nu}}g_1(h)^{-\frac12}}.
\ees
We notice that $g(h) = O(h)$. We then obtain that there exists $c_0,c_1>0$ such that for $|\Re z| \leq c_0g(h)$,\be\label{eq:hypoPi}
\forall u\in D(P)\cap F_h^\bot,\ \vvert{(P-z)u} \geq c_1g(h)\vvert{u}.
\ee

And we can now deduce the second part of Theorem \ref{thm:1} from that, following the same sketch of proof as in \cite{No23}.

Let $\m\in\uuu^{(0)}$, by recalling $f_\m(x,v)=\chi_\m(x) e^{-(f(x,v)-f(\m))/h}$, and because $e^{-f/h}\in\Ker P$, we obtain
\bes
P(f_\m)= [P,\chi_\m]e^{-(f-f(\m))/h} = h\alpha\cdot\nabla\chi_\m e^{-(f-f(\m))/h} = O(e^{-c_\m/h})
\ees
with $c_\m=\inf_{\supp\nabla\chi_\m} f-f(\m)>0$ (because $\chi\equiv1$ near $\m$).

Moreover, with a change of variable $(x,v)\mapsto(h^{\frac{1}{\nu_1}}x_1,\ldots,h^{\frac{1}{\nu_d}}x_d,h^{\frac12}v)$,
\bes
\vvert{f_\m} = Ch^{\frac{d'}{4} + \sum_{i=1}^d\frac{1}{2\nu_i^\m}}(1+O(h^\frac{1}{\overline\nu}))
\ees
for some $C>0$. Thus, since the $(f_\m)_{\m\in\uuu^{(0)}}$ are orthogonal, we actually have :
\be\label{eq:Pfm}
\forall u\in F_h,\ \vvert{Pu}=O(e^{-c/h})\vvert{u}
\ee
for all $c<c_f$ where $c_f=\min_{\m\in\uuu^{(0)}}c_\m>0$. Furthermore, \eqref{eq:Pfm} is also true replacing $P$ by $P^*$ because $P^*(f_\m)=-P(f_\m)$. And because
\bes
P^*Pf_\m = -X(Pf_\m) = (\alpha\cdot h\p_x + \beta\cdot h\p_v)(h\alpha\cdot\nabla\chi_\m)e^{-(f-f(\m))/h},
\ees
\eqref{eq:Pfm} is still valid replacing $P$ by $P^*P$.

We denote by $\Pi_F$ the projector on $F_h$. Let $0<c_0'\leq c_1$, $u\in D(P)$ and $z$ such that $|\Re z|\leq c_0g(h)$ and $|z|\geq c_0'g(h)$,
\bes
\begin{aligned}
    \vvert{(P-z)u}^2 &= \vvert{(P-z)(\Pi_F+\Id-\Pi_F)u}^2\\
    &= \vvert{(P-z)(\Id-\Pi_F)u}^2+\vvert{(P-z)\Pi_F u}^2\\&\phantom{*****}+2\Re\<(P-z)(\Id-\Pi_F)u,(P-z)\Pi_F u\>,
\end{aligned}
\ees
but one has
\bes
\vvert{(P-z)(\Id-\Pi_F)u} \geq c_1g(h)\vvert{(\Id-\Pi_F)u}
\ees
thanks to \eqref{eq:hypoPi}, and 
\bes
\vvert{(P-z)\Pi_F u}^2 \geq (\vvert{P\Pi_F u}-\vvert{z\Pi_F u})^2 \geq \vvert{z\Pi_F u}(\vvert{z\Pi_F u}-2\vvert{P\Pi_F u}).
\ees
Assume now that $g(h)\geq e^{-c/(2h)}$ for some $c<c_f$ (this is \eqref{eq:gthm2}),
\bes
|z|\geq c_0'g(h)\geq c_0'e^{-c/(2h)}\geq c_0'e^{-c/h},
\ees
thus using \eqref{eq:Pfm}, we get
\bes
\vvert{(P-z)\Pi_F u}^2 \geq \frac{|z|^2}{2}\vvert{\Pi_F u}^2.
\ees
Studying each term in the scalar product, there exists $c>0$ such that
\bes
\begin{aligned}
    (\ast) :&= \Re\<(P-z)(\Id-\Pi_F)u,(P-z)\Pi_F u\>\\
    &= \Re\big(\<P(\Id-\Pi_F)u,P\Pi_F u\>-z\<(\Id-\Pi_F)u,P\Pi_F u\> - \bar{z}\<P(\Id-\Pi_F)u,\Pi_F u\>\big)\\
    &= (1+|z|)\vvert{(\Id-\Pi_F)u}\vvert{\Pi_F u}O(e^{-c/h})\\
    &= \Big(\vvert{u}^2+|z|^2\vvert{\Pi_F u}^2+\vvert{(\Id-\Pi_F)u}^2\Big)O(e^{-c/h})
\end{aligned}
\ees
hence
\bes
\begin{aligned}
    \vvert{(P-z)u}^2 &\geq (c_1g(h))^2\vvert{(\Id-\Pi_F)u}^2 + \frac{|z|^2}3\vvert{\Pi_F u}^2\\
    &\phantom{*******}+(\vvert{u}^2+\vvert{(\Id-\Pi_F)u}^2)O(e^{-c/h})\\
    &\geq \frac13(c_0'g(h))^2\vvert{u}^2+(\vvert{u}^2+\vvert{(\Id-\Pi_F)u}^2)O(e^{-c/h})\\
    &\geq \frac14(c_0'g(h))^2\vvert{u}^2
\end{aligned}
\ees
for $h$ small enough, using \eqref{eq:gthm2}. It leads to
\be\label{eq:injectivePz}
\forall u\in D(P),\quad\vvert{(P-z)u} \geq \frac{c_0'}{2}g(h)\vvert{u}.
\ee

By using the same arguments for $P^*$ we have the same result for it (the key point is that $e^{-f/h}$ is in the kernel of $X$ and $N$ hence it also is in $P^*$'s one). It just remains to show that $P-z$ is surjective in order to obtain the resolvent estimate, we show it the classical way, by showing that $\Ran(P-z)$ is closed and dense.

Let $u_n\in D(P)$ and $w\in L^2$ such that $(P-z)u_n\to w$ therefore $((P-z)u_n)_{n\in\N}$ is Cauchy and so is $(u_n)_{n\in\N}$ thanks to \eqref{eq:injectivePz}, hence there exists $u\in L^2$ such that $u_n\to u$. Because the convergence is also true in $\ddd'$, $(P-z)u = w$ in $\ddd'$, and since $w\in L^2$, so is $(P-z)u$, thus $u\in D(P)$ and $\Ran(P-z)$ is closed. Now to show that $\Ran(P-z)$ is dense, we use \eqref{eq:injectivePz} for $P^*$ and so $\Ker(P^*-\overline{z})=\{0\}$.

All this leads to the resolvent estimate
\be\label{eq:resol}
\vvert{(P-z)^{-1}}\leq\frac{2}{c_0'g(h)}.
\ee
Hence, $P$ has no spectrum in
\bes
\{|\Re z|\leq c_0g(h)\}\cap\{|z|\geq c_0'g(h)\}.
\ees
Suppose now that Assumption \ref{ass.Hormander} and \ref{ass.accretive} hold true, from Proposition \ref{prop:accretive} we know that $P$ is maximally accretive and therefore $P - z$ is invertible for all $\Re z < 0$. Moreover we easily see that
\bes
\vvert{(P-z)u}\vvert{u} \geq \Re \<(P-z)u,u\> \geq -\Re z\vvert{u}^2,
\ees
thus for all $\Re z < 0$,
\bes
\vvert{(P-z)^{-1}} \leq \frac{1}{-\Re z}.
\ees
which extends \eqref{eq:resol}
\bes
\forall z\in\{\Re z\leq c_0g(h)\}\cap\{|z|\geq c_0'g(h)\},\ \ \vvert{(P-z)^{-1}}\leq\frac{2}{c_0'g(h)}.
\ees

There is left to show that the spectrum within $\{|z|\leq c_0'g(h)\}$ is composed of $n_0$ eigenvalues exponentially small compared to $h^{-1}$. By denoting $D =D(0,c_0' g(h))$ the disk in $\C$ centered at $0$ of radius $c_0' g(h)$, we consider the spectral projector
\bes
\Pi_D=\frac{1}{2i\pi}\int_{\p\negthinspace D}(z-P)^{-1}dz
\ees
the projector on the small eigenvalues. We start by proving the following lemma
\begin{lemma}\label{lem:Pi_0}
There exists $C>0$ such that $\vvert{P\Pi_D}\leq Cc_0'g(h)$.
\end{lemma}

\bp
\bes
P\Pi_D = \frac{1}{2i\pi}\int_{\p\negthinspace D}P(z-P)^{-1}dz = \frac{1}{2i\pi}\int_{\p\negthinspace D}z(z-P)^{-1}dz,
\ees
hence
\bes
\vvert{P\Pi_D}\leq C(c_0' g(h))^2\frac{2}{c_0'g(h)}
\ees
thanks to \eqref{eq:resol}.

\ep

Let us now prove that $\dim\Ran\Pi_0 = n_0$. We first show that $\dim\Ran\Pi_0\leq n_0$. By contradiction, let us suppose $F_h^\bot\cap\Ran\Pi_D\neq\emptyset$ and so let us take $u\in F_h^\bot\cap\Ran\Pi_D$ of norm one. Since $u\in\Ran\Pi_D$, by Lemma \ref{lem:Pi_0}, $\vvert{Pu}\leq Cc_0' g(h)$, but because $u\in F_h^\bot$ we can use \eqref{eq:hypoPi} and so $\vvert{Pu}\geq c_1g(h)$. Taking $c_0'$ low enough, we have the contradiction we aimed for and thus, $\dim\Ran\Pi_D\leq n_0$.

For the converse inequality, we have
\bes
\Pi_D-\Id=\frac{1}{2i\pi}\int_{\p\negthinspace D}z^{-1}(z-P)^{-1}Pdz
\ees
and therefore
\be\label{quasiortho}
\begin{aligned}
    \eps_\m&=\Pi_Df_\m-f_\m=\frac{1}{2i\pi}\int_{\p\negthinspace D}(z-P)^{-1}P(f_\m)\frac{dz}z\\
    &= O(g(h)^{-1}e^{-c/h})=O(e^{-\frac{c}{2h}})
\end{aligned}
\ee
for some $c>0$, using \eqref{eq:resol}, \eqref{eq:Pfm} and the hypothesis \eqref{eq:gthm2}.

Let us suppose $\sum_{\m\in\uuu^{(0)}} a_\m\Pi_D f_\m = 0$ with $\sum_{\m\in\uuu^{(0)}}  |a_\m|^2 = 1$, since $\Pi_D f_\m = f_\m + \eps_\m$, we have for all $\m'\in\uuu^{(0)}\sum_{\m\in\uuu^{(0)}}  a_\m(\delta_{\m,\m'} + \<\eps_\m,f_{\m'}\>) = 0$ and thus for all $\m, a_\m = O(e^{-c/h})$ for some $c>0$, which is in contradiction with $\sum |a_\m|^2 = 1$. We deduce that $\dim\Ran\Pi_D\geq n_0$ and hence $\dim\Ran\Pi_D = n_0$.

This leads to
\bes
\sigma(P)\cap\{|\Re z|\leq c_0g(h)\}=\{\lambda_\m(h),\m\in\uuu^{(0)}\}\subset D(0,\frac{c_1}{2} g(h)).
\ees

It only remains to show that $\lambda_\m(h) = O(e^{-c/h})$. Noticing that $\Ran\Pi_D$ is $P$-stable and that $(\Pi_Df_\m)_{\m\in\uuu^{(0)}}$ is one of its basis, there exists $C>0$,
\bes
\vvert{P\Pi_D f_\m} = \vvert{\Pi_D Pf_\m} \leq C\vvert{Pf_\m} = O(e^{-c/h}).
\ees
Having that $\vvert{f_\m} = 1 + o(1)$, this yields $P_{|\Ran\Pi_D} = O(e^{-c/h})$m hence $\sigma(P_{|\Ran\Pi_D})\subset D(0,e^{-c/h})$ for some $c>0$.

\section{Examples of hypocoercive operators}\label{sec:examples}

\subsection{Example 1 : a generalization of the adaptive Langevin dynamics.}\label{ssec:exampleAdaptive}

We can try to generalize the process considered in \cite{LeSaSt20}. The starting point is to model a lack of knowledge on the gradient of $V$ by a drift proportional to another random process, in other words, we consider the following SDE
\be\label{eq:SDEAdaptLang}
\left\{
\begin{aligned}
&dx'_t = 2\Sigma^T\Sigma v_tdt,\\
&dv_t = -\p_{x'}V(x'_t)dt - M_t\Sigma^T\Sigma v_tdt - 4\Sigma^T\Sigma v_tdt + \sqrt{2h}dB_t,\\
&dM_t = \Gamma(x'_t,M_t,v_t)dt.
\end{aligned}
\right.
\ee
Where $x'_t,v_t\in\R^d$, $M_t\in\Mr_d(\R)$ is the new variable and $\Gamma$ is to be set so that the system admits an invariant probability measure with the same form as the standard Langevin process. Therefore, denoting $f(x',v,M) = V(x') + |\Sigma v|^2 + W(M)$ with $W$ another smooth real valued function to be determined, having Assumption \ref{ass.skewadj} satisfied, that is $-h\lll^*(e^{-2f/h}) = 0$ is equivalent to \eqref{eq:assskew}, which in this case is
\be\label{eq:skewAdapt1}
4\<M\Sigma^T\Sigma v,\Sigma^T\Sigma v\> - h\div_v(M\Sigma^T\Sigma v) -\Gamma\cdot\p_M W + h\div_M\Gamma = 0,
\ee
where
\bes
\Gamma \cdot \p_MW = \sum_{i,j}\Gamma_{i,j}\p_{M_{i,j}}W \text{ and likewise, }\div_M\Gamma = \sum_{i,j}\p_{M_{i,j}}\Gamma_{i,j}.
\ees
To have a more workable framework, we can consider the case where there exists a fixed real unitary matrix $Q$ such that $Q^T\Gamma Q$ is diagonal. Let us denote $\gamma_i$ such that the $i$-th diagonal block of $Q^T\Gamma Q$ is $\gamma_iI_{r_i}$, $r_i\in\N^*$ such that $\sum r_i = d$. In the following if $A$ is a matrix, then $\diag(A) = (A_{i,i})_i$ and if $u$ is a vector, then $\diag(u) = (u_i\delta_{i,j})_{i,j}$ using $\delta$ the Kronecker symbol. In the next, we will use $\diag$ to switch from vector to matrix and vice versa. Thus we denote $\gamma = \diag(Q^T\Gamma Q)$. This special form for $\Gamma$ induces the same for $M$, let us denote $D$ diagonal by block, which $i$-th block is $y_iI_{r_i}$, such that $M = QDQ^T$. Therefore, the $y_i$ are the new variables. That way, we denote $y = \diag(D)$.

Thus, from the last equation of \eqref{eq:SDEAdaptLang} and the equality $M = QDQ^T$, $\gamma$ has the same form as $D$ and denoting $y$ the new variable,
\bes
dy_{i,t} = \gamma_i(x'_t,QD_tQ^T,v_t)dt.
\ees
Assuming $\Gamma$ does not depend on $y$, then \eqref{eq:skewAdapt1} becomes
\be\label{eq:skewAdapt2}
4\<QDQ^T\Sigma^T\Sigma v,\Sigma^T\Sigma v\> - h\Tr(QDQ^T\Sigma^T\Sigma) - \sum_i\<\gamma_i, \p_{y_i}W\> = 0.
\ee
Now, a natural $W$ to consider is $W(y) = |y|^2$, hence $\p_{y_i}W = 2y_i$.

Noticing now that for any $d\times d$ matrix $A$ and vector $u$ of size $d$ we have
\bes
\Tr(A\diag(u)) = \<\diag(A),u\>,
\ees
thus $\Tr(QDQ^T\Sigma^T\Sigma) = \<\diag(Q^T\Sigma^T\Sigma Q),y\>$. Moreover, using that for any vector $u$, $Du = \diag(u) y$, we obtain
\bes
\<QDQ^T\Sigma^T\Sigma v,\Sigma^T\Sigma v\> = \<\diag(Q^T\Sigma^T\Sigma v)y,Q^T\Sigma^T\Sigma v\>.
\ees

From \eqref{eq:skewAdapt2}, we thus obtain
\be\label{eq:skewAdapt3}
\<4\diag((\Sigma Q)^T\Sigma v)(\Sigma Q)^T\Sigma v - h\diag((\Sigma Q)^T\Sigma Q), y\> = 2\sum_i\<\gamma_i, y_i\>.
\ee

We can consider the case where $Q = I_d$ to simplify the expression of $\gamma$. Because \eqref{eq:skewAdapt3} must be true for any $y$, we get
\bes
\gamma_i = 2\sum_{k = 1+\sum_{n\leq i-1}r_{n}}^{\sum_{n\leq i}r_{n}}(\sum_{j=1}^d(\Sigma^T\Sigma)_{k,j}v_j)^2 - \frac h2\sum_{k = 1+\sum_{n\leq i-1}r_{n}}^{\sum_{n\leq i}r_{n}}(\Sigma^T\Sigma)_{k,k},
\ees
with the convention $r_0 = 0$. When taking $r_1=d$ we have the following
\bes
\gamma_1 = 2|\Sigma^T\Sigma v|^2 - \frac h2\Tr(\Sigma^T\Sigma)
\ees
where we recognize the last equation of \cite[(1.2)]{LeSaSt20} (up to the change $\Sigma\mapsto\frac12\Sigma$).

In the following we will assume $Q = I_d$ and all the multiplicities of the variables $y_i$ are simple, in other words, $\forall i, r_i = 1$, this way we choose $\Gamma = \diag(\gamma)$ and for $i\in\lb1,d\rb$,
\bes
\gamma_i = 2(\sum_{j=1}^d(\Sigma^T\Sigma)_{i,j}v_j)^2 - \frac h2(\Sigma^T\Sigma)_{i,i}.
\ees

Using the decomposition $x = (x',y)$, this leads to the coefficients $\beta(x',y,v) = -\p_{x'}V(x') - \diag(y)\Sigma^T\Sigma v$ and
\bes
\alpha(x',y,v) = \begin{pmatrix}2\Sigma^T\Sigma v\\2\diag(\Sigma^T\Sigma v)\Sigma^T\Sigma v - \frac h2\diag(\Sigma^T\Sigma)\end{pmatrix},
\ees
and we combine the variables to be in the settings of the previous sections: $x = (x',y)$. We now need to check if these coefficients satisfy the different assumptions we made.

\textbullet\ Assumption \ref{ass.skewadj}: We constructed $\alpha$ and $\beta$ from this starting point, so it is satisfied.

\textbullet\ Assumption \ref{ass.confin} and \ref{ass.morse}: Since $W$ is quadratic, as long as $V$ satisfies both of them, the global potential $\hat V(x',y) = V(x') + W(y)$ does.

\textbullet\ Assumption \ref{ass.G}: There, we need to compute the matrix $G = (G_{i,j})_{1\leq i,j\leq 2d}$,
\bes
G_{i,j} = (C_\Sigma h)^{-\frac d2}\int_{\R^d}\alpha_i\alpha_je^{-2|\Sigma v|^2/h}dv.
\ees
For $i,j\leq d$, $G_{i,j} = 0$ for $i\neq j$ using the change of variable $v\mapsto(\Sigma^T\Sigma)^{-1} v$ and a parity argument. Then
\bes
\begin{aligned}
    G_{i,j} &= 4\delta_{i,j}(C_\Sigma h)^{-\frac d2}\int_{\R^d}(\Sigma^T\Sigma v)_i^2e^{-2|\Sigma v|^2/h}dv\\
    &= 4\delta_{i,j}(C_\Sigma h)^{-\frac d2}\det(\Sigma^{-1})\sum_{k,n}\int_{\R^d}\Sigma_{k,i}\Sigma_{n,i}v_kv_ne^{-2|v|^2/h}dv\\
    &= 4\delta_{i,j}(C_\Sigma h)^{-\frac d2}\det(\Sigma^{-1})\sum_{k}\Sigma_{k,i}^2\int_{\R^d}v_k^2e^{-2|v|^2/h}dv\\
    &= 4\delta_{i,j}(C_\Sigma h)^{-\frac d2}\det(\Sigma^{-1})\sum_{k}\Sigma_{k,i}^2\frac h4\int_{\R^d}e^{-2|v|^2/h}dv\\
    &= \delta_{i,j}h(\Sigma^T\Sigma)_{i,i}
\end{aligned}
\ees
using the expression of $\rho$ (see \eqref{eq:rho}).

If $i\leq d$ and $j>d$ with another parity argument, we again have that $G_{i,j} = 0$. For index greater than $d$, we have from the previous computations
\bes
\begin{aligned}
    G_{i+d,j+d} &= (C_\Sigma h)^{-\frac d2}\int_{\R^d}(2(\Sigma^T\Sigma v)_i^2 - \frac h2(\Sigma^T\Sigma)_{i,i})(2(\Sigma^T\Sigma v)_j^2 - \frac h2(\Sigma^T\Sigma)_{j,j})e^{-2|\Sigma v|^2/h}dv\\
    &= 4(C_\Sigma h)^{-\frac d2}\int_{\R^d}(\Sigma^T\Sigma v)_i^2(\Sigma^T\Sigma v)_j^2e^{-2|\Sigma v|^2/h}dv\\
    &\phantom{********} - h(C_\Sigma h)^{-\frac d2}(\Sigma^T\Sigma)_{j,j}\int_{\R^d}(\Sigma^T\Sigma v)_i^2e^{-2|\Sigma v|^2/h}dv\\
    &\phantom{********} - h(C_\Sigma h)^{-\frac d2}(\Sigma^T\Sigma)_{i,i}\int_{\R^d}(\Sigma^T\Sigma v)_j^2e^{-2|\Sigma v|^2/h}dv\\
    &\phantom{********} + \frac{h^2}{4}(C_\Sigma h)^{-\frac d2}(\Sigma^T\Sigma)_{j,j}(\Sigma^T\Sigma)_{i,i}\int_{\R^d}e^{-2|\Sigma v|^2/h}dv\\
    &= 4(C_\Sigma h)^{-\frac d2}\int_{\R^d}(\Sigma^T\Sigma v)_i^2(\Sigma^T\Sigma v)_j^2e^{-2|\Sigma v|^2/h}dv - \frac{h^2}{4}(\Sigma^T\Sigma)_{j,j}(\Sigma^T\Sigma)_{i,i}\\
\end{aligned}
\ees
which, with parity arguments leads to 
\bes
\begin{aligned}
    (\ast_1) :&= 4(C_\Sigma h)^{-\frac d2}\det(\Sigma^{-1})\int_{\R^d}(\Sigma^Tv)_i^2(\Sigma^T v)_j^2e^{-2|v|^2/h}dv\\
    &= 4(C_\Sigma h)^{-\frac d2}\det(\Sigma^{-1})\int_{\R^d}\sum_{k_1,k_2,k_3,k_4}\Sigma_{k_1,i}\Sigma_{k_2,i}\Sigma_{k_3,j}\Sigma_{k_4,j}v_{k_1}v_{k_2}v_{k_3}v_{k_4}e^{-2|v|^2/h}dv\\
    &= 4(C_\Sigma h)^{-\frac d2}\det(\Sigma^{-1})\int_{\R^d}\sum_{k\neq n}\Sigma_{k,i}^2\Sigma_{n,j}^2v_k^2v_n^2e^{-2|v|^2/h}dv\\
    &\phantom{********} + 4(C_\Sigma h)^{-\frac d2}\det(\Sigma^{-1})\int_{\R^d}2\sum_{k\neq n}\Sigma_{k,i}\Sigma_{n,i}\Sigma_{k,j}\Sigma_{n,j}v_k^2v_n^2e^{-2|v|^2/h}dv\\
    &\phantom{********} + 4(C_\Sigma h)^{-\frac d2}\det(\Sigma^{-1})\int_{\R^d}\sum_{k}\Sigma_{k,i}^2\Sigma_{k,j}^2v_k^4e^{-2|v|^2/h}dv\\
    &= \frac{h^2}{4}\big(\sum_{k\neq n}\Sigma_{k,i}^2\Sigma_{n,j}^2 + 2\sum_{k\neq n}\Sigma_{k,i}\Sigma_{n,i}\Sigma_{k,j}\Sigma_{n,j} + 3\sum_{k}\Sigma_{k,i}^2\Sigma_{k,j}^2\big)\\
    &= \frac{h^2}{4}\big((\Sigma^T\Sigma)_{i,i}(\Sigma^T\Sigma)_{j,j} + 2(\Sigma^T\Sigma)_{i,j}^2\big).
\end{aligned}
\ees
Hence $G_{i+d,j+d} = \frac{h^2}{2}(\Sigma^T\Sigma)_{i,j}^2$. We can summarize these computations in
\bes
G = \frac h2\left(\begin{array}{c|c}
    2\diag(\diag(\Sigma^T\Sigma)) & 0 \\
    \hline
    0 & h(\Sigma^T\Sigma)^{\odot 2}
\end{array}\right),
\ees
where $\odot$ is the Hadamard product. We then know using Schur's product Theorem \cite[Theorem VII]{Sc11} that $(\Sigma^T\Sigma)^{\odot 2}$ is definite positive since $\Sigma^T\Sigma$ is, which therefore proves Assumption \ref{ass.G} $i)$ and $ii)$ with $g_1(h)\asymp h^2$ and $g_2(h)\asymp h$. This also makes \eqref{eq:gthm2} true. For $iii)$:

let $i,j\in\lb1,d\rb$, with the scaling $v\mapsto\sqrt{h}v$, we observe that
\bes
\begin{aligned}
    &\<\alpha_i^2\alpha_j^2\rho,\rho\>_{L^2(\R^d_v)} = O(h^2)\\
    &\<\alpha_{i+d}^2\alpha_j^2\rho,\rho\>_{L^2(\R^d_v)} = O(h^3)\\
    &\<\alpha_{i+d}^2\alpha_{j+d}^2\rho,\rho\>_{L^2(\R^d_v)} = O(h^4)
\end{aligned}
\ees

\textbullet\ Assumption \ref{ass.hypocoer}: Using that $\alpha$ is a polynomial of order $2$ in $v$, does not depend on the other variables, and recalling that $\beta(x',y,v) = -\p_{x'}V(x') - \diag(y)\Sigma^T\Sigma v - 4\Sigma^T\Sigma v$, $ii)$ of this assumption is easily satisfied.

Indeed, take $q\in h\{J_x\alpha\alpha,J_v\alpha\beta,J_v\alpha\Sigma^T\Sigma v,h\Delta_v\alpha\}$, we have for all $i,j\in\lb1,d\rb$,
\bes
\Pi q_iq_j\Pi \lesssim h^2(|\nabla \hat V|^2 + h)\Pi
\ees
recalling that $\hat V$ denotes the potential in both $x'$ and $y$, using that $\Pi v_i^k\Pi = Ch^{k/2}\Pi$ for some $C>0$. Noticing that we took $g_2(h)\asymp h^2$, we have the result.

Now for $i)$, taking $j\in\lb1,d\rb$, we have $\<\alpha_j\rho,\rho\>_{L^2(\R^d_v)} = 0$ using parity and from the first set of equation of the previous point, we obtain
\bes
\<\alpha_{j+d}\rho,\rho\>_{L^2(\R^d_v)} = 2(C_\Sigma h)^{-\frac d2}\int_{\R^d}(\Sigma^T\Sigma v)_j^2e^{-2|\Sigma v|^2/h}dv - \frac{h}{2}(\Sigma^T\Sigma)_{j,j} = 0.
\ees

\textbullet Assumption \ref{ass.Hormander}: Because $\Sigma^T\Sigma$ is invertible, the family $\p_{v_i}\alpha_j$ for $i,j\in\lb1,d\rb$ will generate $\p_{x'}$. With the same argument, we obtain $\p_y$ from $\p_{v_k}\p_{v_i}\alpha_{j+d}$, with $k\in\lb1,d\rb$.

\textbullet Assumption \ref{ass.accretive}: We first notice that $\alpha$ does not depend on $x$, then $\div_{x}\alpha = 0$. While for $\beta$, $\div_v\beta = -\Tr(\diag(y)\Sigma^T\Sigma)$ hence $|\div_v\beta| \lesssim |y| \leq f^{1/2}\leq f$ near infinity.

\subsection{Example 2 : a Langevin dynamic rescaled in time.}\label{ssec:exampleRescaled}
Another example can be found in \cite{LeLeLaHi24}, the Langevin dynamics is written with the form

\be\label{eq:SDEg}
\left\{
\begin{aligned}
&dx_t = g(x_t)v_tdt,\\
&dv_t = -g(x_t)\p_xVdt + h\p_xgdt - g(x_t)v_tdt + \sqrt{2hg(x_t)}dB_t.
\end{aligned}
\right.
\ee
And considering $f(x,v) = \frac{V(x)}{2} + \frac{|v|^2}{4}$, $P_g = -e^{f/h}\circ h\lll^*\circ e^{-f/h}$ has the form $P_g = X_g + g\Delta_{\frac{|v|^2}{4}}$ where
\bes
X_g = hg(x)(v\cdot\p_x - \p_xV\cdot\p_v) + h\p_xg\cdot(\frac{v}{2} + h\p_v)
\ees
satisfies
\bes
P_g(e^{-f/h}) = P_g^*(e^{-f/h}) = 0.
\ees
Let us notice that taking $g\equiv 1$ we recover the usual Fokker-Planck operator.

Looking at the details in Section \ref{sec:Hypocoer}, we see that we mostly manipulate $X$ and rather rarely $N$. So we will consider a generalization of the operator extracted from \eqref{eq:SDEg}, taht is $x,v\in\R^d$, $P = X + gN$, with
\bes
\left\{
\begin{aligned}
&X = 2g(x)\Sigma^T\Sigma v\cdot h\p_x - g(x)\p_xV\cdot h\p_v + \frac h2\p_xg\cdot(2\Sigma^T\Sigma v + h\p_v),\\
&N = -h^2\Delta_v + 4|\Sigma^T\Sigma v|^2 - 2h\Tr(\Sigma^T\Sigma).
\end{aligned}
\right.
\ees
From these expressions, we thus consider
\bes
\alpha(x,v) = 2g(x)\Sigma^T\Sigma v\quad \beta(x,v) = - g(x)\p_xV(x) + \frac h2\p_xg(x).
\ees

Like in \cite{LeLeLaHi24}, we assume $g$ satisfies
\be\label{eq:assgExample3}
\exists m,M,\ \forall x\in\R^d,\ 0 < m \leq |\p_xg(x)| \lesssim g(x) \leq M.
\ee
We will see that this condition ensures that all the assumptions for Theorem \ref{thm:1} hold.

Because we have $g$ in factor of $N$, we must be careful to the change it would induce in the proofs of the theorems. First, observe that if $\{g = 0\}$ is of measure $0$, then the kernel of $gN$ and $N$ are the same, and thus $\Pi$ still is the orthogonal projector on the kernel of $gN$. Therefore there is no need to change the definition of $A$, and within the proofs of the lemmas of Section \ref{sec:Hypocoer}, we only use $N$ once, when we bound $NA^*$, which means the other results of that section hold. If $g\in L^\infty$, then $\vvert{gNA^*} \leq \vvert{g}_\infty \vvert{NA^*}$. This lead us to consider both these hypothesis to be true. Now let us look at the other assumptions.

\textbullet\ Assumption \ref{ass.skewadj}: The form of $P$ is made so that this is satisfied.

\textbullet\ Assumption \ref{ass.confin} and \ref{ass.morse}: We just need $V$ to satisfy them.

\textbullet\ Assumption \ref{ass.G}: We observe from the computations of the previous example that $G = hg(x)\diag(\diag(\Sigma^T\Sigma))$. Hence with the assumption that there exists $m,M$ such that for all $x\in\R^d$, we have $0<m\leq g(x)\leq M$ and $\p_ig \lesssim g$ for all $i\in\lb1,d\rb$, then Assumption \ref{ass.G} is satisfied with both $g_1(h)\asymp h$ and $g_2(h)\asymp h$. This also makes \eqref{eq:gthm2} true.

\textbullet\ Assumption \ref{ass.hypocoer}: For $i)$, this is immediate noticing $\alpha$ is odd in $v$. For $ii)$, we need to study each term individually

\quad $\star$ Remark that $hJ_x\alpha\alpha = hg(x)\p_xg|\Sigma^T\Sigma v|^2$, hence we must assume that
\bes
|\p_x(g^2)| \lesssim |\p_x V| + 1,
\ees
which is ensured by \eqref{eq:assgExample3}.

\quad $\star$ We then have $hJ_v\alpha\beta = -2hg(x)^2\Sigma^T\Sigma\p_xV + \frac{h^2}2\Sigma^T\Sigma\p_x(g^2)$. Both term are controlled by $h(|\p_xV| + 1)$ because $g$ is bounded for the first one, and because of the previous point for the other one.

\quad $\star$ Using the two previous points, we also have that $hJ_v\alpha\Sigma^T\Sigma v = 2hg(x)(\Sigma^T\Sigma)^2v$ is well controlled.

\quad $\star$ And lastly $h^2\Delta_v\alpha = 0$.

\textbullet\ Assumption \ref{ass.Hormander}: We have that $\p_{v_i}\alpha_j = 2g(x)(\Sigma^T\Sigma)_{j,i}$. Using that $g$ is bounded from below and $\Sigma^T\Sigma$ is non-degenerate, we can recover all the $\p_{x_i}$.

\textbullet\ Assumption \ref{ass.accretive}: We have that $\div_x\alpha + \div_v\beta = 2\p_xg\cdot\Sigma^T\Sigma$, thus we need to assume $|\p_xg|\lesssim V$ at infinity.

\subsection{Example 3 : adding a magnetic field in the usual KFP operator.}\label{ssec:exampleMagnetic}
The usual KFP equation consist in considering the coefficients
\bes
\alpha(x,v) = 2\Sigma^T\Sigma v;\quad \beta(x,v) = - \p_xV(x).
\ees
In $d = d' = 3$, adding a magnetic field means to add a part of the form $b\wedge v\cdot\p_v$ within the stochastic equation, which means that this term will appear in $\beta$, this leads to
\bes
\alpha(x,v) = 2\Sigma^T\Sigma v;\quad \beta(x,v) = - \p_xV(x) + b(x)\wedge\Sigma^T\Sigma v.
\ees

We now need to check if those coefficients satisfy the assumptions of the previous sections. The part in $2\Sigma^T\Sigma v\cdot h\p_x - \p_xV\cdot h\p_v$ in $P$ is the usual diffusion term of the well-known Fokker-Planck equation, so we will focus on the magnetic term $b(x)\wedge\Sigma^T\Sigma v\cdot \p_v$.

\textbullet\ Assumption \ref{ass.skewadj}: Using that for any $a,u\in\R^3$, $a\wedge u\cdot u = 0$ and since $\alpha$ does not depend on $x$, we have that this assumption is satisfied if and only if $\div_v(b\wedge\Sigma^T\Sigma v) = 0$. If we write $b = (b_i)_{1\leq i\leq3}$ and $M = (M_{i,j})_{i\leq i,j\leq3}$, we observe that
\bes
\div_v(b\wedge Mv) = b_1(M_{2,3} - M_{3,2}) + b_2(M_{3,1} - M_{1,3}) + b_3(M_{1,2} - M_{2,1}),
\ees
hence taking $M = \Sigma^T\Sigma$, which is symmetric, this assumption is satisfied.

We can then write the full operator considered, $P = X + N$, with
\bes
\left\{
\begin{aligned}
&X = 2\Sigma^T\Sigma v\cdot h\p_x - \p_xV\cdot h\p_v + b(x)\wedge\Sigma^T\Sigma v\cdot h\p_v,\\
&N = -h^2\Delta_v + 4|\Sigma^T\Sigma v|^2 - 2h\Tr(\Sigma^T\Sigma).
\end{aligned}
\right.
\ees

\textbullet\ Assumption \ref{ass.confin} and \ref{ass.morse}: We just need to have that $V$ satisfy those.

\textbullet\ Assumption \ref{ass.G}: Since the change in $P$ from the standard KFP operator lies in $\beta$, we have that $G = h\diag(\diag(\Sigma^T\Sigma))$, hence this assumption is satisfied with $g_1(h) = g_2(h) \propto h$.

\textbullet\ Assumption \ref{ass.hypocoer}: For $i)$ this works just like for example \ref{ssec:exampleAdaptive}, and for $ii)$, we only need to check the term $hJ_v\alpha\beta\Pi$, which is $2h\Sigma^T\Sigma\beta\Pi$. The first term appearing is $2h\Sigma^T\Sigma\p_xV\Pi$ which is directly controlled by $g_2(h)|\p_xV|\Pi$. For the other one, we need to work a bit, let $u\in L^2(\R^6)$ and $i\in\lb1,3\rb$,
\bes
\begin{aligned}
    \vvert{h(\Sigma^T\Sigma (b\wedge\Sigma^T\Sigma v))_i\Pi u}^2 &= h^2(C_\Sigma h)^{-\frac d2}\int_{\R^6}u_\rho^2(x)(\Sigma^T\Sigma (b\wedge\Sigma^T\Sigma v))_i^2e^{-2|\Sigma v|^2/h}dvdx\\
    &= h^2(C_\Sigma h)^{-\frac d2}\det(\Sigma^{-1})\\
    &\phantom{********}\times\int_{\R^6}u_\rho^2(x)(\Sigma^T\Sigma (b\wedge\Sigma^T v))_i^2e^{-2|v|^2/h}dvdx\\
    (\Sigma^T\Sigma (b\wedge\Sigma^T v))_i &= \sum_{j=1}^3(\Sigma^T\Sigma)_{i,j}(b\wedge\Sigma^T v)_j.
\end{aligned}
\ees
For $j\in\lb1,3\rb$, we denote $k,\ell\in\lb1,3\rb$ such that $k=j+1[3]$ and $\ell=k+1[3]$, that way we get
\bes
\begin{aligned}
    (b\wedge\Sigma^T v)_j &= b_k(\Sigma^T v)_\ell - b_\ell(\Sigma^T v)_k\\
    &= \sum_{n=1}^3(b_k\Sigma_{n,\ell} - b_\ell\Sigma_{n,k})v_n.
\end{aligned}
\ees
Using a parity argument to get rid of terms of the form $v_nv_m$ where $n\neq m$, we obtain
\bes
\begin{aligned}
    \vvert{h(\Sigma^T\Sigma b\wedge\Sigma^T\Sigma v)_i\Pi u}^2 &=h^2(C_\Sigma h)^{-\frac d2}\det(\Sigma^{-1})\sum_{n=1}^3\Bigg(\int_{\R^3}v_n^2e^{-2|v|^2/h}dv\\
    &\phantom{****}\times\int_{\R^3}u_\rho^2(x)\Big(\sum_{j=1}^3(\Sigma^T\Sigma)_{i,j}(b_k(x)\Sigma_{n,\ell} - b_\ell(x)\Sigma_{n,k})\Big)^2dx\Bigg)\\
    &= \frac{h^3}{4}\int_{\R^3}u_\rho^2(x)\sum_{n=1}^3\Big(\sum_{j=1}^3(\Sigma^T\Sigma)_{i,j}(b_k(x)\Sigma_{n,\ell} - b_\ell(x)\Sigma_{n,k})\Big)^2dx\\
    &= \frac{h^3}{4}\int_{\R^3}u_\rho^2(x)\sum_{n=1}^3(\Sigma^T\Sigma (b(x)\wedge\Sigma^T e_n))_i^2dx\\
\end{aligned}
\ees
denoting $e_n$ the $n$-th element of the canonical basis of $\R^3$. Assuming we have $|b(x)| \lesssim |\p_xV(x)| + 1$ on all $\R^3$, we will then have the assumption verified. Notice that constant or bounded magnetic fields satisfy this bound.

\textbullet\ Assumption \ref{ass.Hormander}: Using that $\p_{v_j}\alpha_i\p_{x_i} = 2(\Sigma^T\Sigma)_{i,j}\p_{x_i}$ this assumption is directly satisfied.

\textbullet\ Assumption \ref{ass.accretive}: Since we have no divergence, this assumption is trivially verified.

\section{Geometric construction of the quasimodes}\label{sec:geoloc}

We now want to have a better view on the small eigenvalues of $P$. For this purpose, we are going to build sharp quasimodes, and so we are following the steps of \cite[section 3\&4]{BoLePMi22}. Their theorem does not apply here because $P$ does not satisfy either (1.9), (Harmo), (Hypo) or (Morse) from their paper for general $\alpha,\beta$ and $V$ we consider. Therefore, the work here is to find a way to obtain a similar result without these assumptions.

Due to several inconvenience, we will consider $V$ Morse from this section onward. However, the other assumptions from \cite{BoLePMi22} still remain false in the general case. To have an overview of what behavior can arise when $V$ is not a Morse function, we refer to \cite{De25_1} where we study the Witten Laplacian $\Delta_V$ with degenerate potentials a similar way to the one we present here.

As in \cite{BoLePMi22}, given $\s\in\vvv^{(1)}$ we look for an approximate solution to the equation $Pu = 0$ in a neighborhood $\www = \www_x\times\www_v$ of $\s$.

We look for an approximate solution of $P u = 0$ of the form
\bes
u=\chi e^{-(f-f(\m))/h},
\ees
and we set
\bes
\chi(x,v)=\int_0^{\ell(x,v,h)}\zeta(s/\tau)e^{-\frac{s^2}{2h}}ds
\ees
where the function $\ell\in\ccc^\infty(\www)$ has a formal classical expansion $\ell\sim\sum_{j\geq0}h^j\ell_j$. Here, $\zeta$ denotes a fixed smooth even function equal to $1$ on $[-1,1]$ and supported in $[-2,2]$, and $\tau > 0$ is a small parameter which will be fixed later.

The object of this section is to construct the function $\ell$. Here we adapt the construction made in \cite[section 3]{BoLePMi22}.

Since $P(e^{-f/h}) = 0$ we have that $P(\chi e^{-f/h}) = [P,\chi](e^{-f/h})$, but
\bes
\begin{aligned}~
    [P,\chi] &= \alpha\cdot h\p_x\chi + \beta\cdot h\p_v\chi - h^2[\Delta_v,\chi]\\
    &= \alpha\cdot h\p_x\chi + \beta\cdot h\p_v\chi - h^2(\Delta_v \chi + 2\p_v \chi\cdot\p_v),\\
    \p_v \chi &= \p_v\ell\zeta(\ell/\tau)e^{-\frac{\ell^{2}}{2h}},\\
    \p_x \chi &= \p_x\ell\zeta(\ell/\tau)e^{-\frac{\ell^{2}}{2h}},\\
    \Delta_v \chi &= \Big(\Delta_v\ell\zeta(\ell/\tau) + \frac1\tau|\p_v\ell|^2\zeta'(\ell/\tau)-\zeta(\ell/\tau)|\p_v\ell|^2\frac{\ell}{h}\Big)e^{-\frac{\ell^{2}}{2h}}.
\end{aligned}
\ees
Therefore, if we set that $\ell_0(\s) = 0$, there exists $r$ smooth such that $r\equiv0$ near $\s$, $r$ and its derivatives are locally uniformly bounded with respect to $h$ and
\be\label{eq:wplusr}
P(\chi e^{-f/h}) = h(w+r)e^{-(f+\frac{\ell^{2}}{2})/h},
\ee
where
\bes
w = \alpha\cdot \p_x\ell + \beta\cdot \p_v\ell + 4\Sigma^T\Sigma v\cdot\p_v\ell + \ell|\p_v\ell|^2 - h\Delta_v\ell.
\ees

If $\alpha$ and $\beta$ admits a similar classical expansion $\ell$'s : $\alpha\sim\sum_{j\geq0}h^j\alpha^j$, $\beta\sim\sum_{j\geq0}h^j\beta^j$, we can see that $w$ admits a formal classical expansion $w \sim \sum_{j\geq0} h^jw_j$. Solving formally $w=0$ and identifying the powers of $h$ leads to a system of equations
\be\label{eq:eik}\tag{eik}
\alpha^0\cdot \p_x\ell_0 + \beta^0\cdot \p_v\ell_0 + 4\Sigma^T\Sigma v\cdot\p_v\ell_0 + \ell_0|\p_v\ell_0|^2 = 0
\ee
and for $j\geq1$,
\be\label{eq:transp}\tag{T$_j$}
\big(\alpha^0\cdot \p_x + (\beta^0+ 4\Sigma^T\Sigma v + 2\ell_0\p_v\ell_0)\cdot\p_v + |\p_v\ell_0|^2\big) \ell_j + R_j = 0
\ee
where $R_j$ is a smooth polynomial of the $\p^\gamma\ell_k$ for $|\gamma|\leq2$ and $k<j$. As in \cite{BoLePMi22} and by analogy with the WKB method, we call eikonal equation the first one and transport equations the next ones.
\begin{proposition}\label{prop:ell}There exists $\ell\in\ccc^\infty(\www)$ satisfying the following
    \begin{itemize}
        \item[$i)$] The eikonal and transport equations are solved up to any order leading to
        \be\label{eq:Pchiell}
            P(\chi e^{-f/h}) = hO(X^\infty + h^\infty)e^{-(f+\frac{\ell^{2}}{2})/h}
        \ee
        denoting $X=(x,v)$ and we have $|\nabla\ell_{0}(\s)|^2\neq0$.
        \item[$ii)$] Moreover, either $|\p_v\ell_{0}(\s)|^2\neq0$ or there exists a non-zero, $h$-independent, positive semidefinite matrix $A$, such that $|\p_v\ell_{0}|^2 = 4\<A\nabla f,\nabla f\>(1+O(X-\s))$, and $\<A\nabla f,\nabla f\>$ is not identically $0$ near $\s$.
        \item[$iii)$] Furthermore, $\ell$ elliptizes $f$ around $\s$ in the following sense
        \be\label{eq:dethess}
        \det\Hess_\s\big(f + \frac{\ell_0^{2}}{2}\big) = -\det\Hess_\s f.
        \ee
        \item[$iv)$] Finally,
        \bes
        \forall X\in\www\setminus\{\s\},\ \ \ X-\s\in\eta(\s)^\bot = 0 \Rightarrow f(X) > f(\s),
        \ees
        where $\eta(\s) = \nabla\ell_{0}(\s)$.
    \end{itemize}
\end{proposition}

\begin{remark}\label{rem:signell}
    Note that $-\ell$ solves the equations the same way $\ell$ does. For now, the sign does not matter, but we will fix it in the next section when properly constructing the quasimodes on a global setting.
\end{remark}

In the following, we will not try to solve \eqref{eq:eik} and \eqref{eq:transp} in their utmost generality due to many different difficulties. We will rather show what behavior can arise when treating with some of the degeneracy one can encounter. We manage to prove Proposition \ref{prop:ell} in all three situations we describe thereafter. Up to translations, we can also assume without loss of generality that $\s = 0$ and $V(\s) = 0$.

\subsection{Situation 1 : partially degenerate $\alpha$ and $\beta$.}\label{sub:ex1geo}

In this first example, we consider a case generalizing a bit the Section \ref{ssec:exampleAdaptive}. We assume the differential $d_{x,v}\alpha^0$ is non-zero (or equivalently from the next equation $d_{x,v}\beta^0\neq0$). We use the relation \eqref{eq:assskew}, which gives when identifying the coefficients of $h^0$:
\bes
\alpha^0\cdot\p_xV + 2\beta^0\cdot\Sigma^T\Sigma v = 0.
\ees
Therefore, assuming $\alpha^0$ has a non-zero linear part, we can write $\alpha^0 = Mv + o(v)$ for some $M\in\Mr_{d',d}(\R)$. Hence we have $\beta^0 = M_\beta x + o(x)$, with $M_\beta\in\Mr_{d,d'}(\R)$, and at its principal order, the previous equation becomes
\bes
(M^TH + 2\Sigma^T\Sigma M_\beta)x\cdot v = 0
\ees
denoting $H$ the Hessian of $V$ at $\s$. This gives the relation $M_\beta = -\frac12(\Sigma^T\Sigma )^{-1}M^TH$. We now denote $X = (x,v)$ and consider $\ell_0 = \xi\cdot X + O(X^2)$ for some $\xi=(\xi_x,\xi_v)\in\R^{d+d'}$. Therefore, at its principal order, \eqref{eq:eik} becomes
\be\label{eq:eik1}
(\Lambda\xi + |\xi_v|^2\xi)\cdot X = 0,
\ee
where $\Lambda = \begin{pmatrix}
    0&-\frac12HM(\Sigma^T\Sigma)^{-1}\\
    M^T&4\Sigma^T\Sigma
\end{pmatrix}$.

We therefore have that \eqref{eq:eik1} is satisfied, for a non-zero $\xi$, if and only if $\xi$ is an eigenvector of $\Lambda$ associated with the eigenvalue $-|\xi_v|^2$. Let $\mu\geq0$,
\be\label{eq:Lambdamu}
\Lambda\xi = -\mu\xi \iff \left\{\begin{array}{l}
    \frac12HM(\Sigma^T\Sigma)^{-1}\xi_v = \mu\xi_x,\\
    M^T\xi_x+ 4\Sigma^T\Sigma\xi_v = -\mu\xi_v.
\end{array}\right.
\ee

Provided that this system has a solution with $\mu\neq0$, we can solve the eikonal equation up to the first order. Note that $\mu\neq0$ implies $\xi_v\neq0$ because we want $\xi\neq0$. We now consider the following assumption
\be\label{ass.simple}\tag{Simple}
    \exists \mu>0,\ \ -\mu\in\sigma(\Lambda)\subset\{-\mu\}\sqcup\{\Re z\geq0\}\text{ and }-\mu\text{ is simple.}
\ee

We observe that this case is not always contained in \cite{BoLePMi22}. If $M^T$ has a non-trivial kernel, then the Kalman-type criterion \cite[Remark 2.5]{BoLePMi22} is not verified, which means that their assumption (Harmo) cannot be satisfied. To see this, let $\Sigma^T\Sigma = I_{d'}$ and $\eta\in\Ker(M^T)$ then with their notations $A^0 = \begin{pmatrix}0&0\\0&I_{d'}\end{pmatrix}$ and $B^T = \begin{pmatrix}0&M_\beta^T\\M^T&0\end{pmatrix}$, hence
\bes
(\eta,0)\in\Ker(A^0)\cap\Ker(B^T)\subset\bigcap_{n=0}^{d+d'-1}\Ker(A^0(B^T)^n).
\ees

Note that as soon as $d>d'$, $M^T$ is not injective, this gives a setting where we have some information that is not present in \cite{BoLePMi22}.

Notice that when $\Sigma = \Id$, we have from \eqref{eq:Lambdamu}
\bes
\Lambda\xi = -\mu\xi \iff \left\{\begin{array}{l}
    \xi_v = \frac{-1}{4+\mu}M^T\xi_x,\\
    HMM^T\xi_x = -2\mu(4+\mu)\xi_x.
\end{array}\right.
\ees
Hence \eqref{ass.simple} is equivalent to the same assumption with $\Lambda$ replaced by $HMM^T$ which may be easier to verify in practice. For example with Section \ref{ssec:exampleAdaptive}, we observe that in this case
\bes
M = \begin{pmatrix}2\Sigma^T\Sigma&0\\0&0\end{pmatrix},
\ees
hence taking $\Sigma^T\Sigma = I_{d'}$ leads to
\bes
HMM^T = \begin{pmatrix}4\Hess_\s V&0\\0&0\end{pmatrix},
\ees
which satisfies \eqref{ass.simple} by definition of $\s$ being a saddle point of $V$.

We now look for $\ell_0$ that admits a decomposition $\ell_0\sim\sum_{j\geq0}\ell_{0,j}$ with $\ell_{0,j}\in\ppp^j_{hom}(X)$, where $\ppp^j_{hom}(X)$ denotes the set of homogeneous polynomials in the $X$ variables of degree $j$. Upon assuming $\alpha^0$ and $\beta^0$ admit similar decompositions, so does $w_0$. Moreover we have, recalling $\nu_1=2$, for all $j\geq1$
\bes
w_{0,j} = ((\Lambda^T + 2A\Pi_\xi)X\cdot\nabla + \mu)\ell_{0,j} + R_{0,j}
\ees
where $A$ is the projector on the velocity: $A:(x,v)\in\R^{d+d'}\mapsto v\in\R^{d'}$, $\Pi_\xi$ the one on $\xi$, from now on we denote $\mu = |\xi_v|^2 = |\p_v\ell_{0,1}|^2$, and $R_{0,j}$ is a smooth polynomial of the $\p^\gamma\ell_{0,k}$ for $|\gamma|\leq1$ and $k<j$.

Under Assumption \eqref{ass.simple}, we have that $\lll_0 = (\Lambda^T + 2A\Pi_\xi)X\cdot\nabla + \mu$ is an automorphism of $\ppp^j_{hom}(X)$ for all $j$. In order to prove this, denoting $\Upsilon = \Lambda^T + 2A\Pi_\xi$, it is sufficient to prove that $\sigma(\Upsilon)\subset\{\Re z\geq0\}$ thanks to \cite[Lemma A.1]{BoLePMi22} (this lemma gives the result for $\{\Re z>0\}$ but we can easily extend it to $\{\Re z\geq0\}$).

In a basis of $\C^{d+d'}$ adapted to $\xi$ in which $\Lambda$ is upper triangular, we have that only the first entry of its diagonal has negative real part, being $-\mu$. Moreover in that same basis, $2\Pi_\xi A$ has zeros outside its first row, and the first element of that row is $2\mu$. Hence all the eigenvalues of $\Upsilon^T$, and thus of $\Upsilon$, have non-negative real part.

This way we solved \eqref{eq:eik} at infinite order: there exists $\ell_0\sim\sum_{j\geq0}\ell_{0,j}$ such that
\be\label{eq:resoleikinf}
\alpha^0\cdot \p_x\ell_0 + \beta^0\cdot \p_v\ell_0 + 4\Sigma^T\Sigma v\cdot\p_v\ell_0 + \ell_0|\p_v\ell_0|^2 = O(X^\infty).
\ee

In order to give a proper definition of $\ell_0\sim\sum_{j\geq0}\ell_{0,j}$, we can use a Borel procedure, which makes the sum converge in $\ccc^\infty(\www)$ and such that $\ell_0$ still solves \eqref{eq:resoleikinf}.

Let us now solve the transport equations \eqref{eq:transp}. We observe that these equations are of the form
\bes
\lll\ell_j + R_j = 0
\ees
with $\lll = \lll_0 + \lll_>$ where $\lll_>(p) = O(X^{j+1})$ for $p\in\ppp^j_{hom}(X)$. Using the same argument as for \eqref{eq:eik}, we can thus solve all the transport equations. And with another Borel procedure in the $h$ variable, we can make sense of $\ell\sim\sum_{j\geq0}h^j\ell_j$ in $\ccc^\infty(\www)$. This leads to
\bes
P(\chi e^{-f/h}) = hO(X^\infty + h^\infty)e^{-(f+\frac{\ell^2}{2})/h}.
\ees

The proof of Proposition \ref{prop:ell} $iii)$, works just like for \cite[Lemma 3.3]{BoLePMi22}. Denoting $\hat H = \begin{pmatrix}H&0\\0&2\Sigma^T\Sigma\end{pmatrix}$, and $E = 1 + \hat H^{-1}\Pi_\xi$, we have that $\Hess_\s(f+\frac{\ell_0^2}{2}) = \hat HE$. Moreover $E \equiv1$ on the hyperplane $\xi^\bot$ and $\<E\xi,\xi\> = -\vvert{\xi}^2$, hence $\det(\hat HE) = -\det \hat H >0$. This also proves $iv)$.

\subsection{Situation 2 : $\alpha$ and $\beta$ with no linear part.}\label{ssec:nolin}

Let us now consider a case where $\alpha$ and $\beta$ are degenerate in all directions. For this purpose we study the operator
\bes
P =  (v^2-h)\cdot h\p_x - 2v\p_xV\cdot h\p_v - h\p_xV + \Delta_{\frac{|v|^2}{4}}
\ees
acting on $(x,v)\in\R^{1+1}$, and with $\Sigma = \frac12$. Although this model does not quite fit in any example of hypocoercive operators given in the previous section, one can easily check that this operator is hypocoercive, with $g_1(h) \asymp g_2(h) \asymp h^2$ (where $g_1$ and $g_2$ are defined in Assumption \ref{ass.G}).

Since we are one-dimensional, we can write $\p_xV = V'$, we will use both notations but mostly the first one, to emphasize the dependence in $x$ and to see how it can be generalized to a multi dimensional case.

This operator obviously does not fall in the framework of \cite{BoLePMi22}, in particular the Kalman-type condition \cite[Remark 2.5]{BoLePMi22} is not satisfied since in their notations, $A^0 = \begin{pmatrix}0&0\\0&1\end{pmatrix}$ and $B = 0$. Recall the equation $w=0$, being
\be\label{eq:w=0situ2}
(v^2-h)\cdot \p_x\ell - 2v\p_xV\cdot \p_v\ell + v\cdot\p_v\ell + \ell|\p_v\ell|^2 - h\Delta_v\ell = 0.
\ee
And the eikonal equation is
\bes
v^2\cdot \p_x\ell_0 - 2v\p_xV\cdot \p_v\ell_0 + v\cdot\p_v\ell_0 + \ell_0|\p_v\ell_0|^2 = 0.
\ees

Working with homogeneous polynomials just like in the previous example we have the principal order of the eikonal equation that is
\bes
v\xi_v + \xi\cdot X|\xi_v|^2 = 0
\ees
recalling $\ell_{0,1}(X) = \xi\cdot X$, $\xi\in\R^2$. We notice that if $\xi_v\neq0$, then we must have $\xi_x=0$. But another purpose of these constructions is to change the nature of the saddle point $\s$ of $V$ into a minimum of $V + \frac{\ell^2}{2}$, which cannot be achieved if $\xi_x=0$. Therefore, we must set $\xi_v=0$, and thus $\xi_x$ remains free.

Before going further in the resolution of the equations, we recall we have $\ell\sim\sum_{j\geq0}h^j\ell_j$, $\ell_j$ $h$-independent, $\ell_j\sim\sum_{k\geq0}\ell_{j,k}$, $\ell_{j,k}\in\ppp^k_{hom}(X)$. We now also consider $\ell_{j,k} = \sum_{a+b=k}\ell_{j,a,b}$ with $\ell_{j,a,b}\in\R x^av^b$. Note that this decomposition is unique. This leads to the same decomposition for $w$ and identifying the monomials that leads to $w_{j,a,b}=0$ gives
\be\label{eq:decompjab}
\begin{aligned}
    &v^2\p_x\ell_{j,a+1,b-2} - \p_x\ell_{j-1,a+1,b} - 2v\sum_{k=1}^a\p_x^{k+1}V(0)\frac{x^k}{k!}\p_v\ell_{j,a-k,b}\\
    &\phantom{********}+ v\p_v\ell_{j,a,b} - \p_v^2\ell_{j-1,a,b+2} + p_{j,a,b} = 0
\end{aligned}
\ee
where $p_{j,a,b}$ contains the terms coming from $|\p_v\ell|^2\ell$, that is
\bes
p_{j,a,b} = \sum\ell_{j_1,a_1,b_1}\p_v\ell_{j_2,a_2,b_2}\p_v\ell_{j_3,a_3,b_3}
\ees
with the sum being on the indexes $j_i,a_i,b_i$, $i\in\lb1,3\rb$ such that
\bes
j_1 + j_2 + j_3 = j;\ a_1 + a_2 + a_3 = a;\ b_1 + b_2 + b_3 = b.
\ees
All this with the convention that $\ell_{j',a',b'} = 0$ if $(j',a',b')\notin\N^3$. The goal now is to find monomials $\ell_{j,a,b}$ that solves \eqref{eq:decompjab} for all $j,a,b\in\N$.

From the structure of the eikonal equation, we have the following lemma.
\begin{lemma}\label{lem:situ2}
    Any function $\ell$ solution of \eqref{eq:w=0situ2} is odd with respect to $v$.
\end{lemma}
We are the one building an $\ell$ solution so we could just look for one that is odd with respect to $v$, but we can actually prove it must be so. Therefore it is not needed to prove this result to continue, this is why we postpone the proof in the appendix \ref{ssec:lemma}.

This property in particular implies that $\p_v\ell = O(v)$ hence plugging $v=0$ in the definition of $w$ we get $\p_x\ell_{|v=0} + \p_v^2\ell_{|v=0} = 0$, which gives for all $j,a\in\N$
\bes
\p_x\ell_{j,a+1,0} + \p_v^2\ell_{j,a,2} = 0.
\ees
We observe that integrating with respect to $v$ this equation gives
\be\label{eq:linkell2}
v\p_x\ell_{j,a+1,0} + \p_v\ell_{j,a,2} = 0,
\ee
which will be useful later.

Note also that thanks to the structure of $w$, having that $\ell$ is even in $v$ directly implies that the odd terms in $w$ are $0$ without further assumptions. Moreover, using that $\p_v\ell = O(v)$, we therefore have that for $b=0$ \eqref{eq:decompjab} is automatically satisfied.

In order to determine the rest of the terms in the expansion of $\ell$, the aim is to solve \eqref{eq:decompjab} in increasing order of $j$, namely solve the eikonal equation first and then the series of the transport ones.

For the eikonal equation, which correspond to $j=0$, we have
\be\label{eq:decomp0ab}
v^2\p_x\ell_{0,a+1,b-2} - v\sum_{k=1}^a\theta_kx^k\p_v\ell_{0,a-k,b} + v\p_v\ell_{0,a,b} + p_{0,a,b}(x,v) = 0,
\ee
where we denote $\theta_k = \frac{2}{k!}\p_x^{k+1}V(0)$. Let us now solve \eqref{eq:decomp0ab} by induction on $b$. We already know that $b=0$ gives a trivial equation and so does $b\in2\N+1$, hence there remains to consider $b\in2\N^*$. For $b=2$ we have
\be\label{eq:decomp0a2-1}
v(v\p_x\ell_{0,a+1,0} + \p_v\ell_{0,a,2}) -v\sum_{k=1}^a\theta_kx^k\p_v\ell_{0,a-k,2} + p_{0,a,2}(x,v) = 0.
\ee
Observe that the first term is zero thanks to \eqref{eq:linkell2}. Note also that because $\ell_{0,a,b} = 0$ for all odd $b$ thanks to Lemma \ref{lem:situ2}, and $\ell_{0,0,0} = 0$, we have $a=0\Rightarrow p_{0,a,2}(x,v)=0$, so for $a=0$, \eqref{eq:decomp0a2-1} is automatically satisfied. And for $a=1$ we get
\bes
- v\theta_{1}x\p_v\ell_{0,0,2} + \ell_{0,1,0}|\p_v\ell_{0,0,2}|^2 = 0,
\ees
using again \eqref{eq:linkell2} which gives $\ell_{0,1,0} = \xi_x x$ and $\p_v\ell_{0,0,2} = -\xi_xv$, we obtain
\bes
v^2\theta_{1}x\xi_x + \xi_xx|\xi_xv|^2 = 0
\ees
leading to $|\xi_x|^2 = -\theta_{1}$ (note that because $\s$ is a saddle point of $V$ according to our definition, $\theta_1 = 2\p_x^{2}V(\s)<0$, hence the negative sign makes sense). Noticing that $\ell$ solves $w=0$ if and only if $-\ell$ do so, we can thus choose $\xi_x =\sqrt{-\theta_1}$ for now.

For $a\geq2$, writing
\bes
p_{0,a,2}(x,v) = \sum_{c + d + e = a}\ell_{0,c,0}\p_v\ell_{0,d,2}\p_v\ell_{0,e,2},
\ees
we observe that (since $\ell_{0,0,0} = 0$)
\bes
p_{0,a,2}(x,v) = \ell_{0,a,0}|\p_v\ell_{0,0,2}|^2 + 2\ell_{0,1,0}\p_v\ell_{0,0,2}\p_v\ell_{0,a-1,2} + q_{0,a,2}(x,v)
\ees
where $q_{0,a,2}$ only contains terms of the form $\p_v\ell_{0,a',2}$ with $a'<a-1$ and $\ell_{0,a'',0}$ with $a''<a$. Therefore we can write \eqref{eq:decomp0a2-1} as
\bes
-v\theta_1x\p_v\ell_{0,a-1,2} + \ell_{0,a,0}|\p_v\ell_{0,0,2}|^2 + 2\ell_{0,1,0}\p_v\ell_{0,0,2}\p_v\ell_{0,a-1,2} = R_{0,a,2},
\ees
with $R_{0,a,2}$ a smooth polynomial of $\p_v\ell_{0,a',2}$ with $a'<a-1$ and $\ell_{0,a'',0}$ with $a''<a$. Using another time \eqref{eq:linkell2}, we can say that $R_{0,a,2}$ is a smooth polynomial of $\ell_{0,a'',0}$ with $a''<a$ exclusively, and that the previous equation can be rewritten
\be\label{eq:decomp0a2}
\lll \ell_{0,a,0} = R_{0,a,2}.
\ee
where $\lll = -v^2\theta_1(x\p_x+1)$ which is invertible over $\ppp^{a}_{hom}(X)$ for all $a\geq2$, then we can solve \eqref{eq:decomp0a2} for all $a\geq2$ by a direct induction. It determines all the $\ell_{0,a,0}$ as well as the $\ell_{0,a,2}$ thanks to \eqref{eq:linkell2}.

In conclusion, solving \eqref{eq:decomp0ab} for all $a\geq0$ and $b\leq2$ determines $\ell_{0,a,b}$ for all $a\geq0$, $b\leq2$.

Now, let $b\geq4$ and assume we have constructed the $\ell_{0,a,b'}$ for all $a\in\N$ and $b'<b$, then \eqref{eq:decomp0ab} can be rewritten
\bes
v\p_v\ell_{0,a,b} = R_{0,a,b}
\ees
with $R_{0,a,b}$ a smooth polynomial of the $\p^\delta\ell_{0,a',b'}$ with $|\delta|\leq1$ and $b'\leq b$ or $b'=b$ but $a'<a$. Because $b\neq0$, $v\p_v$ is invertible over $\R x^av^b$ and we can construct by a direct induction on $a\in\N$ all the $\ell_{0,a,b}$.

For the transport equations, we can observe they have a structure very similar to the eikonal one. Let us solve them by induction on $j$. Let $j\geq1$ and assume we have constructed $\ell_{k,a,b}$ for all $k<j$, $a,b\in\N$. Recall that choosing $b=0$ in \eqref{eq:decompjab} gives $0=0$. For $b=2$ we have
\bes
v(v\p_x\ell_{j,a+1,0} + \p_v\ell_{j,a,2}) - v\sum_{k=1}^a\theta_kx^k\p_v\ell_{j,a-k,2} + p_{j,a,2}(x,v) = R_{j,a,2}
\ees
with $R_{j,a,2} = \p_x\ell_{j-1,a+1,2} + \p_v^2\ell_{j-1,a,4}$. Here again, the first term is zero thanks to \eqref{eq:linkell2}. Consider $a=0$, we obtain
\bes
\ell_{j,0,0}|\p_v\ell_{0,0,2}|^2 = \tilde R_{j,0,2}
\ees
with $\tilde R_{j,0,2}$ a smooth polynomial of the $\p^\delta\ell_{k,a',b'}$ for $|\delta|\leq1$ and $k<j$. Because $|\p_v\ell_{0,0,2}|^2 \neq 0$ this determines $\ell_{j,0,0}$. Let now $a\in\N^*$, assume that the $\ell_{j,a',0}$ are constructed for all $a'<a$. Then by \eqref{eq:linkell2}, the $\ell_{j,a',2}$ are also determined for $a'<a-1$, and we have
\be\label{eq:decompja2}
-v\theta_1x\p_v\ell_{j,a-1,2} + \ell_{j,a,0}|\p_v\ell_{0,0,2}|^2 + 2\ell_{0,1,0}\p_v\ell_{0,0,2}\p_v\ell_{j,a-1,2} = \tilde R_{j,a,2}
\ee
where $\tilde R_{j,a,2}$ is a smooth polynomial of the $\p^\delta\ell_{k,a',b'}$ for $|\delta|\leq1$ and $k<j$ or $k=j$ and either $a'<a$ with $b=0$, or $a'<a-1$. Just like for the eikonal equation, using \eqref{eq:linkell2}, we obtain
\bes
\lll\ell_{j,a,0} = \tilde R_{j,a,2}
\ees
recalling $\lll = -v^2\theta_1(x\p_x+1)$ which is invertible over $\ppp^a_{hom}(X)$. By induction we solve \eqref{eq:decompja2} for all $a\in\N$.

Finally for $b\geq4$, assume we have constructed the $\ell_{j,a,b'}$ for all $a\in\N$ and $b'<b$, then \eqref{eq:decompjab} can be rewritten
\bes
v\p_v\ell_{j,a,b} = R_{j,a,b}
\ees
with $R_{j,a,b}$ a smooth polynomial of the $\p^\delta\ell_{k,a',b'}$ with $|\delta|\leq1$, $k<j$ and either $b'\leq b$ or $b'=b$ but $a'<a$. Because $b\neq0$, $v\p_v$ is invertible over $\R x^av^b$ and we can construct by a direct induction on $a\in\N$ all the $\ell_{j,a,b}$.

Using Borel procedures to give sense to $\ell_j\sim\sum_{k\geq0}\ell_{j,k}$ and $\ell\sim\sum_{j\geq0}h^j\ell_j$ in $\ccc^\infty(\www)$, we have solved $w = O(X^\infty + h^\infty)$ and thus we obtain
\bes
P(\chi e^{-f/h}) = hO(X^\infty + h^\infty)e^{-(f+\frac{\ell^2}{2})/h}.
\ees

Note that for this example, $|\p_v\ell_0(\s)|^2 = 0$, but Proposition \ref{prop:ell} $ii)$ is satisfied because $|\p_v\ell_0|^2 = - 2v^2\p_x^2V(\s)(1+O(X))$. For the last two items, we have
\bes
(f+\frac{\ell_0^2}{2})(\s) = f(\s) + \frac12\p_x^2V(\s)x^2 + \frac{v^2}{4} + \frac{\ell_{0,1,0}^2}{2} + \ell_{0,1,0}\ell_{0,0,2} + O(v^4 + x^4).
\ees
Note that $\ell_{0,1,0}\ell_{0,0,2} = O(|x|v^2) = O(|x|^3 + |v|^3)$, and recall that $\ell_{0,1,0} = \xi_xx$, $|\xi_x|^2 = -\theta_1$ and $\theta_1 = -2\p_x^2V(\s)$. Combined together we have
\bes
(f+\frac{\ell_0^2}{2})(\s) = f(\s) + \frac12|\p_x^2V(\s)|x^2 + \frac{v^2}{4} + O(|(x,v)|^3),
\ees
which proves $iii)$ and $iv)$.

\subsection{Situation 3 : adding a magnetic source.}

Here we consider the operator
\bes
P = 4\Sigma^T\Sigma v\cdot h\p_x - 2\p_xV\cdot h\p_v + b(x)\wedge \Sigma^T\Sigma v\cdot h\p_v + \Delta_{|\Sigma v|^2}
\ees
acting on $(x,v)\in\R^{3+3}$. We will show that this operator satisfies all the hypotheses made in \cite{BoLePMi22} except possibly for $(1.9)$ when $b$ is not bounded. First observe that (Confin), (Gibbs) and (Morse) are already implied by our assumptions, hence there only remains to look at (Harmo) and (Hypo).

For (Harmo), we use \cite[Corollary 2.4 and Remark 2.5]{BoLePMi22}: we have to check that
\bes
\bigcap_{n=0}^{d-1}\Ker(A^0(B^T)^n) = \{0\}
\ees
where with their notation we have $A^0 = \begin{pmatrix}0&0\\0&I_3\end{pmatrix}$ and $B = \begin{pmatrix}0&4\Sigma^T\Sigma\\-2H&\tilde B\end{pmatrix}$, we recall $H$ denotes the Hessian of $V$ at $\s=0$ and here $\tilde B = \p_v(b\wedge\Sigma^T\Sigma v)(0)$. Now, let $(x,v)\in\Ker(A^0)\cap\Ker(A^0B^T)$. From the first kernel, we have that $v=0$, hence the second one gives us $4\Sigma^T\Sigma x = 0$ which implies $x=0$ since $\Sigma$ is invertible.

For (Hypo), we need to look at the dynamic they denoted $c^0(e^{tb^0\cdot\nabla}(x,v))$ for $x$ outside a neighborhood of $E = \{b^0 = 0\}\cap\{c^0 = 0\}$. In our case, we have
\bes
b^0 = \begin{pmatrix}
    4\Sigma^T\Sigma v\\
     - 2\p_xV + b(x)\wedge \Sigma^T\Sigma v
\end{pmatrix}
\ees
and $c^0 = 4|\Sigma^T\Sigma v|^2$.

Let $(x_0,v_0)$ outside a neighborhood of $E$. We consider two cases. Either $v_0\neq0$, then around $(x_0,v_0)$, $c^0$ is uniformly bounded from below, and we obtain that (Hypo) is satisfied. Or $v_0 = 0$ and necessarily, $x_0$ is far enough from a critical point of $V$, therefore $e^{tb^0\cdot\nabla}(x_0,v_0) = e^{- 2t\p_xV(x_0)\cdot\p_v}(x_0,v_0)$. Thus for small time $t>0$, this flow pushes $(x_0,v_0)$ outside a neighborhood of $\R^3_x\times\{0\}$. We can then use the uniform lower bound on $c^0$ to conclude.

This allows us to resolve \eqref{eq:eik} and \eqref{eq:transp} exactly the same way as in \cite{BoLePMi22} since they do not require any assumption at infinity for these local constructions. Their construction of $\ell$ gives the proof of Proposition \ref{prop:ell} $i)$, $ii)$ and $iii)$. For $iv)$ see \cite[Lemma 4.1]{BoLePMi22}.

\begin{remark}
    In \cite{BoLePMi22}, the slow growth assumptions of $(1.9)$ are useful to determine resolvent estimates and rough spectrum localization. If one manage to prove such results without their method, for example using the hypocoercive estimates of Section \ref{sec:Hypocoer}, then satisfying only (Harmo) and (Hypo) of \cite{BoLePMi22} is enough to apply their sharp construction of Gaussian quasimodes, resulting in a precise description and Eyring-Kramers law for the bottom of the spectrum of the operator.
\end{remark}

\section{Global construction}\label{sec:global}

To construct proper quasimodes, we need the notions introduced in Definition \ref{def:labeling} and the labeling given after this. We shall also suppose the generic assumption \eqref{ass.gener} holds true, and we recall it
\bes\tag{Gener}
\begin{array}{l}
    (\ast)\mbox{ for any }\m\in\uuu^{(0)},\m \mbox{ is the unique global minimum of } V_{|E(\m)},\\
    (\ast)\mbox{ for all }\m\neq\m'\in\uuu^{(0)},\j(\m)\cap\j(\m')=\emptyset.
\end{array}
\ees
Given $\m\in\uuu^{(0)}\setminus\{\um\}$, one has $\bsigma(\m) = \sigma_i$ for a certain $i\geq2$. Hence, since $\sigma_{i-1} > \sigma_i$, there exists a unique connected component of $\{V<\sigma_{i-1}\}$ containing $\m$, we denote $E_-(\m)$ that set.

Then we follow the construction of \cite[Section 4]{BoLePMi22}. Notice that in our setting of a kinetic operator, we have two ways of constructing the geometric setup, either by defining the objects on $\R^d_x$ and then tensorizing them by $\R^{d'}_v$, or directly defining the objects on $\R^{d+d'}$. Here we describe the construction on $\R^{d+d'}$. We refer to \cite[Section 3.1]{No23} for a very thorough description of these constructions and a justification of the equivalence between this method and the tensorization one, justifying Definition \ref{def:labeling} can extend from $V$ to $f$.

We recall $\um$ is the unique global minimum of $f$ (the uniqueness is implied by \eqref{ass.gener}). Now let us consider some arbitrary $\m\in\uuu^{(0)}\setminus\{\um\}$. For every $\s\in\j(\m)$, for any $\tau,\delta>0$, we define the sets $\bbb_{\s,\tau,\delta}$, $\ccc_{\s,\tau,\delta}$ and $E_{\m,\tau,\delta}$ by
\bes
\bbb_{\s,\tau,\delta} = \{f\leq f(\s)+\delta\}\cap\{X\in\R^{d+d'},|\eta(\s)\cdot(X - \s)|\leq\tau\},
\ees
\be\label{eq:ccc}
\ccc_{\s,\tau,\delta} \mbox{ the connected component of } \bbb_{\s,\tau,\delta} \mbox{ containing } \s
\ee
and
\bes
E_{\m,\tau,\delta} = \big(E_-(\m)\cap\{f<f(\j(\m))+\delta\}\big)\setminus\bigcup_{\s\in\j(\m)}\ccc_{\s,\tau,\delta},
\ees
where $\eta(\s) = \nabla\ell_{\s,0}(\s)\neq0$ with $\ell_{\s}$ defined in Proposition \ref{prop:ell}. We recall Proposition \ref{prop:ell} $iv)$,
\bes
\forall X\in\www\setminus\{\s\},\ \ \ X-\s\in\eta(s)^\bot = 0 \Rightarrow f(X) > f(\s).
\ees
For $\tau,\delta>0$ small enough, this leads to a partition of $E_{\m,\tau,\delta} = E^+_{\m,\tau,\delta}\sqcup E^-_{\m,\tau,\delta}$ where $E^+_{\m,\tau,\delta}$ is defined as the connected component of $E^+_{\m,\tau,\delta}$ containing $\m$. Let us notice that any path connecting $E^+_{\m,\tau,\delta}$ to $E^-_{\m,\tau,\delta}$ within $\{f<f(\j(\m))+\delta\}$ shall cross at least one $\ccc_{\s,\tau,\delta}$ because of Definition \ref{def:labeling} of the separating saddle points.

We can now define, for $h>0$ and $\tau,\delta$ small enough, the function $\theta_\m$ on the sublevel set $E_-(\m)\cap\{f<f(\j(\m)) + 3\delta\}$ as follows. On the disjoint open sets $E^+_{\m,3\tau,3\delta}$ and $E^-_{\m,3\tau,3\delta}$, we define
\bes
\theta_\m(X)=\left\{\begin{aligned}
    1\mbox{ for } X\in E^+_{\m,3\tau,3\delta},\\
    -1\mbox{ for } X\in E^-_{\m,3\tau,3\delta}.
\end{aligned}\right.
\ees
In addition, for $X\in\ccc_{\s,3\tau,3\delta}$ we set
\bes
\theta_\m(X) = C_{\s,h}^{-1}\int_0^{\ell_\s(X)}\zeta(r/\tau)e^{-\frac{r^{2}}{2 h}}dr,
\ees
where the function $\ell_\s$ is the one constructed in the previous section and its sign (see Remark \ref{rem:signell}) is chosen so that there exists a neighborhood $\www$ of $\s$ such that $E(\m)\cap\www$ is included in the half plane $\{\eta(\s)\cdot(X - \s)>0\}$. We recall that $\zeta\in\ccc_c^\infty(\R,[0,1])$ is even and satisfies $\zeta = 1$ on $[-1,1]$ and $\zeta=0$ outside $[-2,2]$, and we set the normalizing constant
\bes
C_{\s,h}=\frac 12\int_{-\infty}^{+\infty}\zeta(r/\tau)e^{-\frac{r^{2}}{2h}}dr.
\ees
Therefore, since, for every $\tau>0$ and then $\delta>0$ small enough, the sets $E^+_{\m,3\tau,3\delta}$, $E^-_{\m,3\tau,3\delta}$ and $\ccc_{\s,3\tau,3\delta}$ for $\s\in\j(\m)$ are mutually disjoint, $\theta_\m$ is well defined on their reunion which forms $E_-(\m)\cap\{V<V(\j(\m))+3\delta\}$. Moreover, on a small neighborhood of the common boundary between $\ccc_{\s,3\tau,3\delta}$ and $E_{\m,3\tau,3\delta}$, we have that $|\eta(\s)\cdot(X - \s)|\geq\frac52\tau$ or in other words, $|\ell_{\s,0,1}|\geq\frac52\tau$. Using that $\ell_\s = \ell_{\s,0,1} + O(|x-\s|^2 + h)$ and having that $|x-\s| = O(\delta)$ in this neighborhood, we thus obtain that for every $\tau > 0$ small and then $\delta,h>0$ small enough, $|\ell_\s| \geq 2\tau$ in a neighborhood of the boundary between $\ccc_{\s,3\tau,3\delta}$ and $E_{\m,3\tau,3\delta}$. This shows that $\theta_\m$ is $\ccc^\infty$ on $E_-(\m)\cap\{V<V(\j(\m))+3\delta\}$.

Note also that there exists $\gamma,\eps>0$ such that
\bes
\begin{aligned}
    C_{\s,h} &= \frac 12\int_{-\infty}^{+\infty}\zeta(r/\tau)e^{-\frac{r^{2}}{2h}}dr = \int_0^{+\infty}\zeta(r/\tau)e^{-\frac{r^{2}}{2h}}dr\\
    &= \int_0^{+\infty}e^{-\frac{r^{2}}{2h}}dr + \int_0^{+\infty}(\zeta(r/\tau)-1)e^{-\frac{r^{2}}{2h}}dr\\
    &= \sqrt{\frac{\pi h}{2}} + O\big(\int_\gamma^{+\infty}e^{-\frac{r^{2}}{2h}}dr\big)\\
    &= \sqrt{\frac{\pi h}{2}}(1+O(e^{-\eps/h})).
\end{aligned}
\ees
Hence
\be\label{eq:ch}
\exists\eps>0,\ C_{\s,h}^{-1} = \sqrt{\frac{2}{\pi h}}(1+O(e^{-\eps/h})).
\ee

We now want to extend $\theta_\m$ to a cutoff defined on $\R^{d+d'}$. Considering a smooth function $\chi_\m$ such that
\bes
\chi_\m(X) = \left\{\begin{aligned}
    &1\mbox{ for } X\in E_-(\m)\cap\{f \leq f(\j(\m)) + 2\delta\},\\
    &0\mbox{ for } X\in\R^{d+d'}\setminus\big(E_-(\m)\cap\{f < f(\j(\m)) + 3\delta\}\big),
\end{aligned}\right.
\ees
we have that $\chi_\m\theta_\m$ belongs to $\ccc_c^\infty(\R^{d+d'},[-1,1])$ and
\bes
\supp(\chi_\m\theta_\m)\subset E_-(\m)\cap\{f < f(\j(\m)) + 3\delta\}.
\ees

\begin{defin}\label{def:quasimodes}
For $\tau>0$ and then $\delta,h>0$ small enough, we define the quasimodes 
\bes
\left\{\begin{array}{l}
    \psi_{\um}(X)=2e^{-\frac{f(X)-f(\um)}{h}}\\
    \psi_{\m}(X)=\chi_\m(X)(\theta_\m(X)+1)e^{-\frac{f(X)-f(\m)}{h}}\quad\mbox{ for }\m\in\uuu^{(0)}\setminus\{\um\}.
\end{array}\right.
\ees
And at the same time, we define the normalized quasimodes for $\m\in\uuu^{(0)}$ by
\bes
\phii_{\m}=\frac{\psi_{\m}}{\vvert{\psi_{\m}}}.
\ees
\end{defin}

From this definition, we have the following Lemma which gives us a first relationship between the quasimodes
\begin{lemma}\label{lem:supppsi}Let $\m\neq\m'\in\uuu^{(0)}$,

    If $\bsigma(\m) = \bsigma(\m')$ and $\j(\m)\cap\j(\m')=\emptyset$ then $\supp(\psi_\m)\cap\supp(\psi_{\m'}) = \emptyset$.

    If $\bsigma(\m) > \bsigma(\m')$, then
    \begin{itemize}
        \item[$\star$] either $\supp(\psi_\m)\cap\supp(\psi_{\m'}) = \emptyset$,
        \item[$\star$] or $\psi_\m = 2e^{-(f-f(\m))/h}$ on $\supp(\psi_{\m'})$.
    \end{itemize}
\end{lemma}

The proof is the same as for \cite[Lemma 4.4]{BoLePMi22}. We recall it for the reader's convenience.

\bp Because $E(\m)$ is a connected component of $\{f<\bsigma(\m)\}$, the boundary of $\overline{E(\m)}$ is made of non-critical points of $f$ (points where $\nabla f\neq0$) and separating saddle points $\s\in\j(\m)$ by the definition of a separating saddle point. Therefore, from the definition of $\psi_\m$, we see that for all $\eps>0$,
\bes
\supp\psi_\m\subset \overline{E^+_{\m,3\tau,3\delta}\cup\bigcup_{\s\in\j(\m)}\ccc_{\s,3\tau,3\delta}}\subset\overline{E(\m)}+B(0,\eps)
\ees
for $\tau,\delta$ small enough. Now if $\bsigma(\m) = \bsigma(\m')$ then necessarily $E(\m)\cap E(\m') = \emptyset$. If it was not the case, because they are critical components of $\{f\leq\bsigma(\m)\}$, they would be the same which is in a contradiction with the construction of $E$. When in addition $\j(\m)\cap\j(\m')=\emptyset$, then
\bes
\overline{E(\m)}\cap\overline{E(\m')} = \p E(\m)\cap\p E(\m') = \j(\m)\cap\j(\m') = \emptyset,
\ees
hence with $\eps$ sufficiently small we have $\supp(\psi_\m)\cap\supp(\psi_{\m'}) = \emptyset$. If $\bsigma(\m) > \bsigma(\m')$ then either $\m'\notin E(\m)$ in which case we have with the above that $\supp(\psi_\m)\cap\supp(\psi_{\m'}) = \emptyset$. And if $\m'\in E(\m)$ then $\overline{E(\m')}\subset E_-(\m')\subset E(\m)$ but $\chi_\m\theta_\m\equiv1$ on a neighborhood of $\overline{E(\m)}\setminus\bigcup_{\s\in\j(\m)}\ccc_{\s,3\tau,3\delta}$ and $\j(\m)\cap\supp\psi_{\m'}=\emptyset$ hence the last result.

\ep

We denote
\be\label{eq:defmu}
\mu(\m) = 1 + \inf_{\s\in\j(\m)}b_{\ell_\s},
\ee
where $b_{\ell_\s}$ is defined in \eqref{eq:defab}. It is a power that will appear in the following, and we denote $\tilde\j(\m)\subset\j(\m)$ the set of saddle points that satisfy the infimum, in other words
\be\label{eq:deftildej}
\s_0\in\tilde\j(\m) \iff b_{\ell_\s} = \min_{\s\in\j(\m)}b_{\ell_\s}.
\ee

\begin{proposition}\label{prop:interactionmatrix}
Suppose that Proposition \ref{prop:ell} holds, under Assumption \eqref{ass.gener}, there exists $C>0$ such that for $\tau>0$ and then $\delta,h>0$ small enough, for every $\m,\m'\in\uuu^{(0)}$,
\begin{itemize}
    \item[$i)$] $\<\phii_\m,\phii_{\m'}\> = \delta_{\m,\m'}+O(e^{-C/h}),$
    \item[$ii)$] $\<P\phii_\m,\phii_\m\> = \frac1{2\pi}\sum_{\s\in\tilde\j(\m)}\frac{(\det\Hess_\m f)^{\frac12}}{|\det\Hess_\s f|^{\frac12}}a_{\ell_\s}h^{\mu(\m)}e^{-2S(\m)/h}(1+O(h))$, with $\mu$ defined in \eqref{eq:defmu} and $a_{\ell_\s}$ defined in \eqref{eq:defab}. Where $S(\m) = f(\j(\m))-f(\m) = V(\j(\m))-V(\m)$ for $\m\neq\um$ and $S(\um)=+\infty$ as denoted in \eqref{eq:defS}.
    \item[$iii)$] $\vvert{P\phii_\m}^2 = O(h^\infty)\<P\phii_\m,\phii_\m\>$
    \item[$iv)$] $\vvert{P^*\phii_\m}^2 = O(h^{-c})\<P\phii_\m,\phii_\m\>$ for some $c\in\R$.
\end{itemize}
\end{proposition}

\bp We will follow the proof of \cite[Proposition 5.1]{BoLePMi22} and will not explain every arguments that did not change from their proof. In this proof we will use several Laplace method (LM), they are all justified thanks to Proposition \ref{prop:ell} $iii)$.

Noticing $f$ uniquely attains its global minimum at $\m$ on $\supp\psi_\m$, by using a LM applied to $2f$ we obtain for $\m\neq\um$
\bes
\begin{aligned}
    \vvert{\psi_\m}^2 &= \int_{\R^{d+d'}}\chi_\m^2(\theta_\m+1)^2e^{-2\frac{f-f(\m)}{h}}\\
    &= \chi_\m^2(\m)(\theta_\m(\m)+1)^2\frac{(h\pi)^{\frac {d+d'}2}}{(\det\Hess_\m f)^{\frac12}}(1+O(h)).
\end{aligned}
\ees
This leads to
\be\label{eq:eclpsi}
\vvert{\psi_\m} = 2\frac{(h\pi)^{\frac {d+d'}4}}{(\det\Hess_\m f)^{\frac14}}(1+O(h)).
\ee
Now, notice that this result also holds for $\m = \um$. Then, the proof of $i)$ is exactly the same as in \cite{BoLePMi22} using Lemma \ref{lem:supppsi}.

For $ii)$, Using the computation done after the proof of \cite[Proposition 5.1 $i)$]{BoLePMi22} of $\<P\psi_\m,\psi_\m\>$ we have 
\bes
\begin{aligned}
    \<P\psi_\m,\psi_\m\> &= h^2\sum_{\s\in\j(\m)}C_{\s,h}^{-2}\int_{\ccc_{\s,3\tau,3\delta}}\chi_\m^2\zeta(\ell_\s/\tau)^2|\p_v\ell_\s|^2e^{-2\big(f+\frac{\ell_\s^{2}}{2}-f(\m)\big)/h}\\&\phantom{***************}+O(e^{-2(S(\m)+2\delta)/h}).
\end{aligned}
\ees
Thanks to Proposition \ref{prop:ell} $iii)$, we have that $\s$ is the unique minima of $f+\frac{\ell^{2}_{\s,0}}{2s}-f(\m)$ on $\ccc_{\s,3\tau,3\delta}$ and
\bes
\big(f+\frac{\ell^{2}_{\s,0}}{2}-f(\m)\big)(\s) = f(\s) - f(\m) = f(\j(\m)) - f(\m) = S(\m).
\ees
Moreover, thanks to Remark \ref{rem:relax} and Proposition \ref{prop:ell} $iii)$, $f+\frac{\ell^{2}_{\s,0}}{2}$ satisfies Assumption \ref{ass.morse}, thus we can use a LM and we obtain
\bes
\begin{aligned}
    \<P\psi_\m,\psi_\m\> &= h^2\sum_{\s\in\j(\m)}C_{\s,h}^{-2}\chi_\m^2(\s)\zeta(\ell_{\s,0}(\s)/\tau)^2|\p_v\ell_{\s,0}(\s)|^2\frac{(h\pi)^{\frac {d+d'}2}}{(\det\Hess_\s f)^{\frac12}}\\
    &\phantom{************}\times e^{-2S(\m)/h}(1+O(h))\\
\end{aligned}
\ees
in the case where $|\p_v\ell_{\s,0}(\s)|^2\neq0$. Otherwise, thanks to Proposition \ref{prop:ell} and a LM, we replace $|\p_v\ell_{\s,0}(\s)|^2$ by $h\div(A\nabla f)(\s)$, with $A$ given by Proposition \ref{prop:ell} (we have an extra factor $2$ because we use the LM with $\phii = 2f$). For the example of Subsection \ref{ssec:nolin}, this gives $-{2}{\p_xV(\s)}h$ (recall $d=d'=1$ in this case). In the following, we will use
\be\label{eq:defab}
\begin{aligned}
    &a_{\ell_\s} = |\p_v\ell_{\s,0}(\s)|^2 \text{ and } b_{\ell_\s} = 0 \text{ if } |\p_v\ell_{\s,0}(\s)|^2 \text{ is non-zero, and}\\
    &a_{\ell_\s} = \div(A\nabla f)(\s) \text{ and } b_{\ell_\s} = 1 \text{ otherwise.}
\end{aligned}
\ee
Therefore we can write
\bes
\begin{aligned}
    \<P\psi_\m,\psi_\m\> &= h^2\sum_{\s\in\j(\m)}C_{\s,h}^{-2}\chi_\m^2(\s)\zeta(\ell_{\s,0}(\s)/\tau)^2a_{\ell_\s}h^{b_{\ell_\s}}\frac{(h\pi)^{\frac {d+d'}2}}{(\det\Hess_\s f)^{\frac12}}\\
    &\phantom{************}\times e^{-2S(\m)/h}(1+O(h)).\\
\end{aligned}
\ees
Thanks to \eqref{eq:ch}, this leads to 
\be\label{eq:Ppsipsi}
\begin{aligned}
    \<P\psi_\m,\psi_\m\> &= \sum_{\s\in\j(\m)}a_{\ell_\s}\frac{2(h\pi)^{\frac {d+d'}2}}{\pi(\det\Hess_\s f)^{\frac12}}h^{1+b_{\ell_\s}}e^{-2S(\m)/h}(1+O(h))\\
    &= \sum_{\s\in\tilde\j(\m)}a_{\ell_\s}\frac{2\pi^{\frac {d+d'}2-1}}{(\det\Hess_\s f)^{\frac12}}h^{\mu(\m) + \frac {d+d'}2}e^{-2S(\m)/h}(1+O(h)).
\end{aligned}
\ee
Combining \eqref{eq:Ppsipsi} and \eqref{eq:eclpsi}, we obtain
\bes
\<P\phii_\m,\phii_\m\> = \frac1{2\pi}\sum_{\s\in\tilde\j(\m)}\frac{(\det\Hess_\m f)^{\frac12}}{|\det\Hess_\s f|^{\frac12}}a_{\ell_\s}h^{\mu(\m)}e^{-2S(\m)/h}(1+O(h)),
\ees
this proves $ii)$.

Let us now prove $iii)$. Starting as in \cite{BoLePMi22}'s work, we have that
\be\label{eq:Dvpsim}
\vvert{P\psi_\m}^2 = \vvert{P(\theta_\m e^{-(f-f(\m))/h})}^2_{L^2(\supp\chi_\m)} + O(e^{-2(S(\m)+2\delta))/h}).
\ee
And on $\supp\chi_\m$, $P(\theta_\m e^{-(f-f(\m))/h})$ is supported in $\bigcup_{\s\in\j(\m)}\ccc_{\s,3\tau,3\delta}$, thus using \eqref{eq:Pchiell}, we obtain with a LM
\bes
\begin{aligned}
    \vvert{P(\theta_\m e^{-(f-f(\m))/h})}^2_{L^2(\ccc_{\s,3\tau,3\delta})} &= \int_{\ccc_{\s,3\tau,3\delta}} O(X^\infty + h^\infty)e^{-2\big(f-f(\m)+\frac{\ell_\s^{2}}{2}\big)/h}\\
    &= O(h^\infty)e^{-2S(\m)/h}\\
    &= O(h^\infty)\<P\psi_\m,\psi_\m\>_{L^2(\ccc_{\s,3\tau,3\delta})}.
\end{aligned}
\ees
Combining this results and \eqref{eq:Dvpsim}, we proved point $iii)$.

For $iv)$, notice that around $\s\in\j(\m)$, $P^*(\theta_\m e^{-f/h}) = h(\hat w + r)e^{-(f+\frac{\ell_\s^{2}}{2})/h}$ with
\bes
\hat w = - \alpha\cdot\p_x\ell_\s - \beta\cdot\p_v\ell_\s + 4\Sigma^T\Sigma v\cdot\p_v\ell_\s + \ell_\s|\p_v\ell_\s|^2 - h\Delta_v\ell_\s,
\ees
therefore with a LM and using the analog of \eqref{eq:Dvpsim} with $P^*$ instead of $P$, we obtain the result.

\ep

\subsection{Proof of Theorem \ref{thm:2} and graded matrices}\label{ssec:proofthm2}
The proof is the same as in \cite{BoLePMi22} we refer to it and to \cite{LePMi20} for the details, let us mention the main arguments. From now on, we relabel the minima by increasing saddle height:
\be\label{eq:decreasSmu}
\forall j\in\lb1,n_0-1\rb,\left\{\begin{aligned}
    &S(\m_j) \leq S(\m_{j+1}),\\
    &\mu(\m_{j+1}) \leq \mu(\m_j) \text{ if }S(\m_j) = S(\m_{j+1})
\end{aligned}\right.
\ee
Thus $\m_{n_0} = \um$ because $S(\um) = +\infty$ and $\um$ is the only minima that has this property by construction. We also denote for shortness
\bes
\forall j\in\lb1,n_0\rb,\ \ S_j = S(\m_j),\ \ \phii_j = \phii_{\m_j},\ \mbox{ and } \ \tilde{\lambda}_j = \<P\phii_j,\phii_j\>.
\ees
Therefore we have
\be\label{eq:Pmm'}
\forall j,k\in\lb1,n_0\rb,\ \ \<P\phii_j,\phii_k\> = \delta_{j,k}\tilde{\lambda}_j.
\ee
This statement is obvious when $j=k$. Now if $j > k$, from Lemma \ref{lem:supppsi} either $\supp\phii_j\cap\supp\phii_k = \emptyset$ or $\phii_j = c_he^{-(f-f(\m_j)/h}$ on $\supp\phii_k$, $c_h$ being a normalizing constant (or the same swapping $j$ and $k$). Using that $P(e^{-f/h}) = P^*(e^{-f/h}) = 0$, we see that \eqref{eq:Pmm'} is indeed true. Let $\Cr = \partial D(0,\frac{c_0}2g(h))$. We now define
\bes
\Pi_\Cr = \frac{1}{2i\pi}\int_{\Cr}(z-P)^{-1}dz
\ees
the spectral projector on the small eigenvalues of $P$, where $c_0$ is given by Theorem \ref{thm:1}, $g(h)$ is defined in \eqref{eq:g}. Observe that thanks to Theorem \ref{thm:1}, $\vvert{\Pi_\Cr}\leq \frac{c_0}{c_0'}$. Hence, denoting $u_j = \Pi_\Cr\phii_j$ for $j\in\lb1,n_0\rb$ (we notice that $u_{n_0} = \phii_{\um}$), we obtain the following proposition

\begin{proposition}
    Under Assumption \ref{ass.polyg}, there exists $c>0$ such that for every $j,k\in\lb1,n_0\rb$ and every $h>0$ small enough, one has
    \be\label{eq:quasiortho}
    \<u_j,u_k\> = \delta_{j,k} + O(e^{-c/h})
    \ee
    and
    \be\label{eq:interaction}
    \<Pu_j,u_k\> = \delta_{j,k}\tilde{\lambda}_j + O\Big(h^\infty\sqrt{\tilde{\lambda}_j\tilde{\lambda}_k}\Big).
    \ee
\end{proposition}

\bp Using that
\bes
\Pi_\Cr-1 = \frac{1}{2i\pi}\int_{\Cr}(z-P)^{-1}P\frac{dz}{z}
\ees
we have that for $u\in D(P)$, $\vvert{(\Pi_\Cr-1)u}\leq\sup_{\Cr}\Vert(z-P)^{-1}\Vert\vvert{Pu}$, and thus
\bes
\<u_j,u_k\> = \<\phii_j,\phii_k\> + g(h)^{-1}O\big(\Vert P\phii_j\Vert + \Vert P\phii_k\Vert\big)
\ees
using the resolvent estimate given by Theorem \ref{thm:1}. Using now Assumption \ref{ass.polyg}, we obtain $\<u_j,u_k\> = \delta_{j,k} + O(e^{-c/h})$ thanks to Proposition \ref{prop:interactionmatrix}. Therefore, thanks to Theorem \ref{thm:1} and \eqref{eq:Pmm'}, we then obtain
\bes
\begin{aligned}
    \<Pu_j,u_k\> &= \<P\phii_j,\phii_k\> + \<P(\Pi_\Cr-1)\phii_j,\phii_k\> + \<P\Pi_\Cr\phii_j,(\Pi_\Cr-1)\phii_k\>\\
    &= \delta_{j,k}\tilde{\lambda}_j + g(h)^{-1}O\big(\Vert P\phii_j\Vert\Vert P^*\phii_k\Vert + \Vert P\phii_j\Vert\Vert P\phii_k\Vert\big)\\
    &= \delta_{j,k}\tilde{\lambda}_j + O\Big(h^\infty\sqrt{\tilde{\lambda}_{j}\tilde{\lambda}_{k}}\Big).
\end{aligned}
\ees

\ep

This grants to the interaction matrix a graded structure we will develop below. Then we use the Gram-Schmidt process to transform the basis $(u_{n_0-j+1})_{1\leq j\leq n_0}$ into an orthonormal basis $(e_{n_0-j+1})_{1\leq j\leq n_0}$ of $\Ran\Pi_\Cr$. Moreover, thanks to \eqref{eq:quasiortho} we have that
\bes
\forall j\in\lb1,n_0\rb,\ \ e_j = u_j + O(e^{-c/h}),
\ees
see \cite[Lemma 4.11]{LePMi20} for the details. And thus, using \eqref{eq:interaction} along with \cite[Proposition 4.12]{LePMi20}, we have that
\be\label{eq:Pee}
\forall j\in\lb1,n_0\rb,\ \ \<Pe_j,e_k\> = \delta_{j,k}\tilde{\lambda}_j + O\Big(h^\infty\sqrt{\tilde{\lambda}_j\tilde{\lambda}_k}\Big).
\ee
Using its graded structure, we can now compute the eigenvalues of the matrix
\be\label{eq:matrix}
M := (\<Pe_j,e_k\>)_{j,k} = {P}_{|\Ran\Pi}.
\ee
We recall the results stated in \cite{LePMi20}:

We denote by $\Dr_0(E)$ the set of invertible and diagonalizable complex matrices of an Euclidean space $E$.

\begin{defin}\cite[Definition A.1]{LePMi20}
    Let $\Er = (E_j)_{1\leq j \leq p}$ be a sequence of vector spaces $E_j$ of (finite) dimension $r_j>0$, let $E = \bigoplus_{j=1}^pE_j$ and let $\tau = (\tau_i)_{2\leq i\leq p}\in(\R_+^*)^{p-1}$. Suppose that $(h,\tau)\mapsto \mmm_h(\tau)$ is a map from $(0,1]\times(\R_+^*)^{p-1}$ to the set of complex matrices on $E$.

    We say that $\mmm_h(\tau)$ is an $(\Er,\tau,h)$-graded matrix if there exists $\mmm'\in\Dr_0(E)$ independent of $(h,\tau)$ such that $\mmm_h(\tau) = \Omega(\tau)(\mmm' + O(h))\Omega(\tau)$ with $\Omega(\tau)$ and $\mmm'$ such that
    \begin{itemize}
        \item $\mmm' = \diag(M_j,1\leq j\leq p)$ with $M_j\in\Dr_0(E_j)$,
        \item $\Omega(\tau) = \diag(\eps_j(\tau)I_{r_j},1\leq j\leq p)$ with $\eps_1(\tau) = 1$ and $\eps_j(\tau) = \prod_{k=2}^j\tau_k$ for  $j\geq2$.
    \end{itemize}
\end{defin}

\begin{theorem}\cite[Theorem A.4]{LePMi20}\label{thm:graded}
    Suppose that $\mmm_h(\tau)$ is $(\Er,\tau,h)$-graded. Then, there exists $\tau_0, h_0>0$ such that for all $0<\tau_j<\tau_0$ and $h\in(0,h_0]$, one has
    \bes
    \sigma(\mmm_h(\tau))\subset\bigsqcup_{j=1}^p\eps_j(\tau)^2(\sigma(M_j) + O(h)).
    \ees
    Moreover, for any eigenvalue $\lambda$ of $M_j$ with multiplicity $m_j(\lambda)$, there exists $K>0$ such that, denoting $D_j(\lambda) = \{z\in\C,\ |z-\eps_j(\tau)^2\lambda| < \eps_j(\tau)^2Kh\}$, one has
    \bes
    n(D_j(\lambda);\mmm_h(\tau)) = m_j(\lambda),
    \ees
    where $n(D_j;\mmm_h(\tau))$ is the rank of the Riesz projector associated with $\mmm_h(\tau)$ with contour $D_j$. Finally, there exists $C>0$, such that for all $z\in\C\setminus\bigcup_{j=1}^p\bigcup_{\lambda\in\sigma(M_j)}D_j(\lambda)$,
    \be\label{eq:resolventGraded}
    \vvert{(\mmm_h(\tau) - z)^{-1}} \leq Cd(z,\sigma(\mmm_h(\tau)))^{-1}.
    \ee
\end{theorem}

We want to apply this theorem to \eqref{eq:matrix}, we therefore need to show that it is graded. Let us first notice that
\bes
M = \begin{pmatrix}M'&0\\0&0\end{pmatrix}
\ees
with $M'\in\Mr_{n_0-1}(\R)$ because by construction $e_{n_0}$ is the ground state of $P$. Now, for all $j\in\lb1,n_0-1\rb$, let us define
\be\label{eq:defvmu}
v_j = \frac1{2\pi}\sum_{\s\in\tilde\j(\m)}\frac{(\det\Hess_\m f)^{\frac12}}{|\det\Hess_\s f|^{\frac12}}a_{\ell_\s},\ \mbox{ and }\ \mu_j = 1 + b_{\ell_\s} = \mu(\m_j)
\ee
for some $\s\in\tilde\j(\m_j)$ (recalling $\tilde \j$ is defined in \eqref{eq:deftildej}), with $a_{\ell_\s}$ and $b_{\ell_\s}$ defined in \eqref{eq:defab}. Therefore we have that
\be\label{eq:lamdbav}
\tilde\lambda_j = v_jh^{\mu_j}e^{-2S_j/h}(1 + O(h)).
\ee
Because of \eqref{eq:decreasSmu}, there exists a partition $J_1\sqcup\ldots\sqcup J_p$ of $\lb1,n_0-1\rb$ such that for all $k\in\lb1,p\rb$, there exists $\iota(k)\in\lb1,n_0-1\rb$ such that
\bes
\forall j\in J_k,\ \ S_j = S_{\iota(k)},\ \mu_j = \mu_{\iota(k)}
\ees
and
\bes
\forall 1\leq k< k'\leq p,\ \ S_{\iota(k)} < S_{\iota(k')}\ \text{ or }\ S_{\iota(k)} = S_{\iota(k')}\ \text{ and }\  \mu_{\iota(k')} < \mu_{\iota(k)}
\ees
Hence, using \eqref{eq:Pee} and \eqref{eq:lamdbav} we have that $\big(h^{\mu_1}e^{-2S_1/h}\big)^{-1}M'$ is $(\Er,\tau,h)$-graded with $\Er = (\R^{J_k})_{1\leq k\leq p}$ and $\tau_k = h^{\frac{\mu_{\iota(k)}-\mu_{\iota(k-1)}}{2}} e^{-\frac{S_{\iota(k)}-S_{\iota(k-1)}}{h}}$ for $k\geq 2$. We can now apply Theorem \ref{thm:graded} and we obtain
\bes
\sigma(M') \subset \bigsqcup_{k=1}^ph^{\mu_1}e^{-2S_1/h}\eps_k^2(\sigma(M_k) + O(h))
\ees
with $M_k = \diag(v_j, j\in J_k)$, $h^{\mu_1}e^{-2S_1/h}\eps_k^2 = h^{\mu_{\iota(k)}}e^{-2S_{\iota(k)}/h}$ and corresponding multiplicities. All this leads to the announced result.

\subsection{Proof of Corollaries \ref{cor:return} and \ref{cor:times}}\label{ssec:proofcor}We follow the proof of \cite[Corollaries 1.5 and 1.6]{BoLePMi22}.

Recall that $\Cr = \partial D(0,\frac{c_0}2g(h))$ and
\bes
\Pi_\Cr = \frac{1}{2i\pi}\int_{\Cr}(z-P)^{-1}dz.
\ees
We already obtained that $\vvert{\Pi_\Cr} \leq \frac{c_0}{c_0'}$. Moreover, recall \eqref{eq:resol} that is
\bes
\vvert{(P-z)^{-1}}\leq\frac{2}{c_0'g(h)},
\ees
for $\{\Re z \leq c_0g(h)\}\cap\{|z|\geq c_0'g(h)\}$. Since the operator valued function $(P-z)^{-1}(1-\Pi_\Cr)$ is holomorphic on $\{\Re z \leq c_0g(h)\}$, the maximum principle yields
\be\label{eq:resol2}
\vvert{(P-z)^{-1}(1-\Pi_\Cr)}\leq\frac{2}{c_0'g(h)},
\ee
for $\Re z \leq c_0g(h)$

The solution of \eqref{eq:evolutionP} can be written
\be\label{eq:developEvol}
u(t) = e^{-tP/h}u_0 = e^{-tP/h}\Pi_\Cr u_0 + e^{-tP/h}(1-\Pi_\Cr)u_0.
\ee
Let $Q:\Im(1-\Pi_\Cr)\to\Im(1-\Pi_\Cr)$ be the operator $P$ restricted to the Hilbert space $\Im(1-\Pi_\Cr)$. Since $P$ is maximally accretive from Proposition \ref{prop:accretive}, so is $Q$ and thus $\vvert{e^{-tQ/h}}\leq1$. Moreover, \eqref{eq:resol2} shows that $\vvert{(Q-z)^{-1}}\leq Cg(h)^{-1}$ for some $C>0$ and $\Re z \leq c_0g(h)$. To estimate the last term in \eqref{eq:developEvol}, we use a Gearhardt-Pr\"uss type inequality with an explicit bound. More precisely, \cite[Proposition 2.1]{HeSj10_01} (see also \cite{HeSj21_01,HeSjVi24}) gives, for some $C>0$ and all $t\geq0$,
\bes
\vvert{e^{-tQ/h}} \leq \Big(1 + 2\frac{c_0}2g(h)\sup_{\Re z = \frac{c_0}2g(h)}\vvert{(Q-z)^{-1}}\Big)e^{-t\frac{c_0}2g(h)/h} \leq Ce^{-t\frac{c_0}2\frac{g(h)}{h}}.
\ees
Therefore, denoting $\varepsilon = c_0/2$ and using Assumption \ref{ass.polyg}, there exists $C>0$ such that for all $t\geq0$,
\be\label{eq:semigroupQ}
\vvert{e^{-tQ/h}(1-\Pi_\Cr)u_0} \leq Ce^{-\varepsilon h^{c-1} t}\vvert{u_0},
\ee
with $c\geq1$ given by Assumption \ref{ass.polyg}. Noticing now that $P_{|\Im\Pi_\Cr}$ is a matrix of size $n_0$ whose eigenvalues are the $\lambda(\m,h)$'s, using \eqref{eq:developEvol} and \eqref{eq:semigroupQ}, we obtain \eqref{eq:return1} from the usual formula for the exponential of a matrix applied to $e^{-tP/h}\Pi_\Cr u_0$.

Let us now prove \eqref{eq:return2}. Denoting $\Pi_0$ the orthogonal projector onto the kernel of $P$, we need to prove
\bes
e^{-tP/h} = \Pi_0 + O(e^{-t\underset{\m\neq\um}{\min}\Re(\lambda(\m,h))(1-Ch)/h}),
\ees
with a $O$ uniform in $t$ and $h$. Writing
\bes
e^{-tP/h} - \Pi_0 = e^{-tP/h}\Pi_\Cr - \Pi_0 + e^{-tP/h}(1-\Pi_\Cr),
\ees
it remains to bound $e^{-tP/h}\Pi_\Cr - \Pi_0$ thanks to \eqref{eq:semigroupQ}. Using that
\bes
\Pi_0u = \<\frac{e^{-f/h}}{\vvert{e^{-f/h}}},u\>\frac{e^{-f/h}}{\vvert{e^{-f/h}}},
\ees
which is a bounded operator, we have $\Pi_0 = e^{-tP/h}\Pi_0$ and it suffices to show that
\bes
\exists C>0,\ \ \vvert{e^{-tP/h}(\Pi_\Cr-\Pi_0)} \leq Ce^{-t\underset{\m\neq\um}{\min}\Re(\lambda(\m,h))(1-Ch)/h},
\ees
or in other words,
\bes
\exists C>0,\ \ \vvert{e^{-tM'/h}} \leq Ce^{-t\underset{\m\neq\um}{\min}\Re(\lambda(\m,h))(1-Ch)/h},
\ees
with $M'$ the matrix of $P$ in the orthonormal basis $(e_j)_{1\leq j\leq n_0-1}$ of $\Ran(\Pi_\Cr-\Pi_0)$ defined in the preceding subsection. We already showed, using Theorem \ref{thm:graded} that
\bes
\sigma(M') \subset \bigsqcup_{\m\in\uuu^{(0)}\setminus\{\um\}} D\big(v(\m)h^{\mu(\m)}e^{-2S(\m)/h},Kh^{\mu(\m)+1}e^{-2S(\m)/h}\big),
\ees
with $v$ and $\mu$ defined in \eqref{eq:defvmu} and some constant $K>0$. Therefore, using the functional calculus representation of $e^{-tM'/h}$ along with the resolvent estimate of Theorem \ref{thm:graded}, we obtain
\bes
\begin{split}
    \vvert{e^{-tM'/h}} &= O\Big(\sup_{\m\neq\um}e^{-tv(\m)h^{\mu(\m)}e^{-2S(\m)/h}(1-\frac K2h)/h}\Big)\\
    &= O\Big(e^{-t\underset{\m\neq\um}{\min}\Re(\lambda(\m,h))(1-\frac K2h)/h}\Big).
\end{split}
\ees

Let us now prove Corollary \ref{cor:times}. Recall that, in the proof of Theorem \ref{thm:2}, we introduced an application $\iota:\lb1,p+1\rb\to\lb1,n_0\rb$ such that for all $1\leq k<k'\leq p+1$,
\bes
h^{\mu_{\iota(k)}}e^{-S_{\iota(k)}/h} > h^{\mu_{\iota(k')}}e^{-S_{\iota(k')}/h}.
\ees
In the following we will omit the function $\iota$ and just write $\mu_k,S_k$ to avoid heavy expressions.

For $R>1$ we define the balls
\bes
\forall k\in\lb1,p\rb,\ \ \ D_k = D\big((R+R^{-1})h^{\mu_k}e^{-2S_k/h},Rh^{\mu_k}e^{-2S_k/h}\big)
\ees
and $D_{p+1} = D(0,R^{-1}h^{\mu_p}e^{-2S_{p}/h})$. For $R$ fixed large enough and every $h$ small enough, each exponentially small eigenvalue of $P$ belongs to exactly one of the disjoint sets $D_k$ from Theorem \ref{thm:2}. Moreover, $\partial D_k$ is at distance of order $h^{\mu_k}e^{-2S_k/h}$ (resp. $h^{\mu_p}e^{-2S_{p}/h}$) from the spectrum of $P$ for $k\in\lb1,p\rb$ (resp. $k=p+1$). Using the resolvent estimate of $P$ on the image of $\Pi_\Cr$ given by Theorem \ref{thm:graded} and \eqref{eq:resol2} to control the contribution on the image of $1-\Pi_\Cr$, we get
\be\label{eq:resol3}
\forall z\in \partial D_k,\ \ \vvert{(P-z)^{-1}}\leq \left\{
\begin{aligned}
    &Ch^{-\mu_k}e^{2S_k/h}&&\text{for }k\in\lb1,p\rb,\\
    &Ch^{-\mu_p}e^{2S_{p}/h}&&\text{for }k=p+1,
\end{aligned}
\right.
\ee
using Assumption \ref{ass.polyg}. In particular, the spectral projector associated with the eigenvalues of order $h^{\mu_k}e^{-2S_k/h}$,
\bes
\Pi_k = \frac{1}{2i\pi}\int_{\partial D_k}(z-P)^{-1}dz,
\ees
is well-defined and satisfies $\vvert{\Pi_k}\leq C$.

We can now decompose
\be\label{eq:decomProj}
e^{-tP/h}\Pi_\Cr = \sum_{k=1}^{p+1}e^{-tP/h}\Pi_k.
\ee
For $k\in\lb1,p\rb$ and $0\leq t\leq t_k^-$, \eqref{eq:resol3} and $t_k^-e^{-2S_k/h} = O(h^\infty)$ imply
\be\label{eq:estim1}
\begin{split}
    e^{-tP/h}\Pi_k &= \frac{1}{2i\pi}\int_{\partial D_k}e^{-tz/h}(z-P)^{-1}dz\\
    &= \frac{1}{2i\pi}\int_{\partial D_k}(z-P)^{-1}dz + \frac{1}{2i\pi}\int_{\partial D_k}(e^{-tz/h} - 1)(z-P)^{-1}dz\\
    &= \Pi_k + \int_{\partial D_k}O(t|z|/h)\vvert{(P-z)^{-1}}dz\\
    &= \Pi_k + O(th^{\mu_k-1}e^{-2S_k/h})\\
    &= \Pi_k + O(h^\infty).
\end{split}
\ee
using that $\mu_k\geq1$. On the contrary, for $t_k^+\leq t$, \eqref{eq:resol3} and $e^{-t_k^+h^{\mu_k-1}e^{-2S_k/h}/R} = O(h^\infty)$ give
\be\label{eq:estim2}
\begin{split}
    e^{-tP/h}\Pi_k &= \frac{1}{2i\pi}\int_{\partial D_k}e^{-tz/h}(z-P)^{-1}dz\\
    &= O\Big(\int_{\partial D_k}e^{-t\Re z/h}\vvert{(P-z)^{-1}}dz\Big)\\
    &= O\big(e^{-tR^{-1}h^{\mu_k}e^{-2S_k/h}/h}\big)\int_{\partial D_k}\vvert{(P-z)^{-1}}dz\\
    &= O(h^\infty).
\end{split}
\ee
Lastly, $e^{-tP/h}\Pi_{p+1} = \Pi_{p+1}$ since $\Pi_{p+1}$ is the rank-one spectral projector associated with the eigenvalue $0$. On the other hand, since $e^{-\varepsilon h^{c-1}t_0^+} = O(h^\infty)$, \eqref{eq:semigroupQ} becomes
\be\label{eq:estim3}
\vvert{e^{-tP/h}(1-\Pi_\Cr)} = O(h^\infty),
\ee
for $t\geq t_0^+$. Summing up, Corollary \ref{cor:times} is a direct consequence of the formulas \eqref{eq:resol3} and \eqref{eq:decomProj} with the relation $\Pi_k^\leq = \sum_{j=k}^{p+1}\Pi_j$ and the estimates \eqref{eq:estim1}, \eqref{eq:estim2} and \eqref{eq:estim3}.

\appendix
\section{Labeling of the critical points}

In this section, we consider a smooth function $W\in\ccc^\infty(\R^n,\R)$ satisfying Assumption \ref{ass.morse} that we recall thereafter. We denote by $\uuu$ the set of its critical points.
\begin{assumption}\label{ass:morsedegen2}
    For any critical point $x^*\in\uuu$, there exists a neighborhood $\vvv\ni x^*$, $(t_i^{x^*})_{1\leq i\leq n}\subset \R^{*}$, $(\nu_i^{x^*})_{1\leq i\leq n}\subset \N\setminus\{0,1\}$, a $\ccc^\infty$ change of variable $U^{x^*}$ defined on $\vvv$ such that $U^{x^*}(x^*)=x^*$, $U^{x^*}$ and $d_{x^*}U^{x^*}$ are invertible and
    \bes
    \forall x\in\vvv,\ \ W\circ U^{x^*}(x) - W(x^*) = \sum_{i=1}^nt_i^{x^*}(x_i-x_i^*)^{\nu_i^{x^*}}.
    \ees
\end{assumption}
We consider the partition $\uuu = \uuu^{odd}\sqcup\uuu^{even}$ where $x^*\in\uuu^{even} \iff \forall i,\ \nu_i^{x^*}\in 2\N$. Now we shall say that $x^*\in\uuu^{even}$ is of index $j\in\lb0,n\rb$ if $\sharp\{i\ |\ t_i^{x^*}<0\}=j$ and therefore $\uuu^{even}=\bigsqcup_{j=0}^n\uuu^{(j)}$ where $\uuu^{(j)}$ is the set of critical points of $V$ with even order in each direction and of index $j$, we also denote $n_0 = \sharp\uuu^{(0)}$ the number of minima of $W$. In the following when speaking of saddle point, we only refer to critical points of index $1$, that is elements of $\uuu^{(1)}$, and we will mostly denote them by the letter $\s$ and its variations.

In works in the lineage of \cite{BoLePMi22}, for Morse functions, the important critical points to study are the minima as they generates the eigenvalues and some saddle points because they are the only ones such that $B(\s,r)\cap\{W < W(\s)\}$ has exactly two connected components (for $r>0$ small enough, where $B(\s,r)\subset\R^n$ denotes the open ball of center $\s$ and radius $r$), meaning that if a process starts from one side of this set, it has to cross the level set $\{W = W(\s)\}$ in order to go the other side. One can show that for any other critical point, this set is connected therefore the process need not to climb the landscape given by $W$ in order to cross wells. We need to prove a similar result for functions that are not Morse anymore, but satisfy Assumption \ref{ass:morsedegen2}.

\begin{proposition}\label{prop:saddle}
    Let $\s\in\uuu$, for $r>0$ small enough, the set $B(\s,r)\cap\{W < W(\s)\}$ has exactly two connected components if $\s\in\uuu^{(1)}$ and one otherwise.
\end{proposition}

\subsection{Proof of Proposition \ref{prop:saddle}}

We first need to prove the following lemma, adapted from \cite[Lemma 3.1]{AsBoMi22}, but before that, we define a family of paths. Let $a=(a_1,\ldots,a_n)\in[1,+\infty)^n$, $x,y\in\R^n$, we denote
\bes
\begin{array}{c|ccc}
     \gamma^a_{x,y}:&[0,1]&\to&\R^n\\
     &t&\mapsto&(t^{a_i}(x_i-y_i))_i + y.
\end{array}
\ees
We see that $\gamma^a_{x,y}(0) = y$ and $\gamma^a_{x,y}(1) = x$, and when $a\equiv 1$, this is just the standard linear parametrization of the segment $[y,x]$.
\begin{lemma}\label{lem:path}
    Let $\phii$ be a local smooth diffeomorphism of $\R^n$ defined in a neighborhood of $b\in\R^n$ and $a=(a_1,\ldots,a_n)\in[1,+\infty)^n$. Then there exists $r_b>0$ such that for all $0<r<r_b$,
    \bes
    \forall x\in\phii(B(b,r)),\ \gamma^a_{x,\phii(b)}([0,1])\subset\phii(B(b,r)).
    \ees
\end{lemma}

\bp
In the following, we will denote $\underline a = \inf_ia_i$ and $\overline a = \sup_ia_i$. We have to show that for all $x\in B(0,r)$ and $t\in[0,1]$, the point $\gamma^a_{\phii(b+x),\phii(b)}(t)$ belongs to $\phii(B(b,r))$. In other words,
\bes
\forall t\in[0,1],\ g(t) = |\phii^{-1}(\gamma^a_{\phii(b+x),\phii(b)}(t))-b|^2\in[0,r^2).
\ees
First, we see that $g(0) = 0$ and $g(1) = |x|^2 < r^2$. Then, there exists $E$ a set large enough such that
\bes
\begin{aligned}
    g(t) &= |\phii^{-1}(\gamma^a_{\phii(b+x),\phii(b)}(t))-b|^2\\
    &\leq \vvert{d\phii^{-1}}_{L^\infty(E)}^2 |\gamma^a_{\phii(b+x),\phii(b)}(t)-\phii(b)|^2\\
    &= \vvert{d\phii^{-1}}_{L^\infty(E)}^2 |t^{a_i}(\phii(b+x)_i-\phii(b)_i)|^2\\
    &\leq \vvert{d\phii^{-1}}_{L^\infty(E)}^2\vvert{d\phii}_{L^\infty(B(b,r))}^2 t^{2\underline a}|x|^2.
\end{aligned}
\ees
Thus there exists $C>0$ such that for all $\eps>0$ and all $t\leq\eps$,
\bes
g(t) \leq C\eps^{2\underline a}|x|^2.
\ees
Choosing $\eps$ such that $C\eps^{2\underline a} \leq 1$, we have that $g(t) < r^2$ for $t\in[0,\eps]$.
Moreover, the Taylor formula implies
\bes
\begin{aligned}
    g'(t) &= 2\< \p_t(\phii^{-1}\circ\gamma^a_{\phii(b+x),\phii(b)})(t),\phii^{-1}\circ\gamma^a_{\phii(b+x),\phii(b)}(t) - b\>\\
    &= 2\< d_{\gamma^a_{\phii(b+x),\phii(b)}(t)}\phii^{-1}({\gamma^a_{\phii(b+x),\phii(b)}}'(t)),\phii^{-1}\circ\gamma^a_{\phii(b+x),\phii(b)}(t) - b\>\\
    &= 2\< d_{\phii(b)+O(t^{\underline a}x)}\phii^{-1}\big((a_it^{a_i-1}d_b\phii_i(x))_i + O(t^{\underline a-1}x^2)\big),\\&\phantom{************}d_{\phii(b)}\phii^{-1}\big((t^{a_i}d_b\phii_i(x))_i\big) + O(t^{2\underline a}x^2)\>.
\end{aligned}
\ees
Here we denoted $\phii_i$ the $i$-th component of $\phii$. Denoting $A = d_{\phii(b)}\phii^{-1}$ and $\Gamma_1(t) = \diag(t^{a_i-1})$, $\Gamma_2(t) = \diag(a_it^{a_i-1})$, we observe that
\bes
g'(t) = 2t\<A\Gamma_2(t)A^{-1}x,A\Gamma_1(t)A^{-1}x\> + \vvert{A\Gamma_1(t)A^{-1}x}O(t^{\underline a}x^2) + O(t^{2\underline a-1}x^3).
\ees
But using that $A$ is invertible and thus $A^*A$ is positive definite, we have that
\bes
\begin{aligned}
    \<A^*A\Gamma_2(t)A^{-1}x,\Gamma_1(t)A^{-1}x\> &\gtrsim \<\Gamma_2(t)A^{-1}x,\Gamma_1(t)A^{-1}x\>\\
    &\gtrsim t^{2\overline a-2}\<(AA^*)^{-1}x,x\>\\
    &\gtrsim t^{2\overline a-2}|x|^2.
\end{aligned}
\ees
Thus we have for $\eps$ defined above
\bes
\exists 0<c_{\eps}<1,r_{\eps}>0,\forall |x|<r_{\eps},\forall t\geq\eps,\ \ g'(t) \geq c_{\eps}|x|^2
\ees
and it leads to
\bes
\forall t\geq\eps,\ g(t) \leq (1-(1-t)c_{\eps})|x|^2 < r_{\eps}^2.
\ees

\ep

We can now prove the proposition.

\bp
Up to a translation, we can consider that $\s=0$ and $W(\s)=0$. Thanks to Assumption \ref{ass:morsedegen2}, we know that there exists $U, (t_i)_i\subset\R^*, (\nu_i)_i\subset\N\setminus\{0,1\}$ such that $U$ is a smooth diffeomorphism in a neighborhood of $0$, $U(0) = 0$ and for $x$ in that neighborhood,
\bes
W\circ U (x) = \sum_{i=1}^nt_ix_i^{\nu_i}.
\ees

We denote $X_0 = \{W<0\}$, let $\overline \nu = \sup\negthinspace_i\nu_i$, $a = (\overline\nu/\nu_i)_i\in[1,+\infty)^n$ and $x\in U^{-1}(X_0\cap B(0,r))$ with $0<r<r_0$, $r_0$ given by Lemma \ref{lem:path}. We see that
\bes
\forall s\in[0,1],\ W\circ U (\gamma^a_{x,0}(s)) = \sum_{i=1}^nt_i(s^{\frac{\overline\nu}{\nu_i}}x_i)^{\nu_i} = s^{\overline\nu}W\circ U(x),
\ees
hence
\be\label{eq:a}\tag{a}
\gamma^a_{x,0}((0,1])\subset U^{-1}(X_0).
\ee
Moreover, applying Lemma \ref{lem:path} to $\phii = U^{-1}$ and $b=0$, we have
\be\label{eq:b}\tag{b}
\gamma^a_{x,0}((0,1])\subset U^{-1}(B(0,r)),
\ee
thus, combining \eqref{eq:a} and \eqref{eq:b}, we have
\be\label{eq:c}\tag{c}
\gamma^a_{x,0}((0,1])\subset U^{-1}(X_0\cap B(0,r)).
\ee
Now, let us notice that
\be\label{eq:d}\tag{d}
\gamma^a_{x,0}(s)\xrightarrow[s\to0]{}0
\ee
and
\be\label{eq:e}\tag{e}
\forall \eta \ll 1,\ U^{-1}(X_0\cap B(0,r))\cap B(0,\eta) = U^{-1}(X_0)\cap B(0,\eta)
\ee
because $U^{-1}$ is an open map since it is a diffeomorphism. Thus, up to proving that
\be\label{eq:effe}\tag{f}
U^{-1}(X_0)\cap B(0,\eta) \mbox{ is connected}
\ee
we have
\bes
X_0\cap B(0,r) \mbox{ is connected}
\ees
combining \eqref{eq:c}, \eqref{eq:d}, \eqref{eq:e} and \eqref{eq:effe}, because connectedness is invariant by diffeomorphism. Therefore it just remains to prove \eqref{eq:effe}. Now let us show that if $0\notin\uuu^{(1)}$, then $U^{-1}(\{W<0\})\cap B(0,\eta) = \{x\in B(0,\eta),\ W\circ U(x)<0\}$ is connected for $\eta>0$ small enough.

We consider the case where $0\in\uuu^{odd}$ which means that there is at least one odd $\nu$. Without any loss of generality, we can consider that $\nu_1$ is odd and $t_1<0$ (in the case $t_1>0$, we just replace $\eta/2$ by $-\eta/2$ in the following).

Start from a point $(x_1,\ldots,x_n)$ in $\{x\in B(0,\eta),\ W\circ U(x)<0\}$, the goal is to connect it to $(\eta/2,0,\ldots,0)$ via a path that remains within the set. In the following, by saying that we link $a$ to $b$ we mean that we create a segment from one to another (so a path of the form $t\mapsto a + t(b-a)$), one can check that each time we do that, the whole path stays in $\{x\in B(0,\eta),\ W\circ U(x)<0\}$.

The first step of the path is to link all the positive $t_ix_i^{\nu_i}$ to $0$, hence either $x_1=0$ or $t_1x_1^{\nu_1}<0$ and all other non-zero contributions to $W\circ U$ are negative ones. Then we halve all the $x_i$ (so we link $x_i$ to $x_i/2$), this way, the ending point is in $B(0,\eta/2)$ and we know that its first coordinate is non-negative. Now we link $x_1$ to $\eta/2$, because we were in $B(0,\eta/2)$ and we had $x_1\geq0$, we indeed remain in $B(0,\eta)$. And at last, we link all the other coordinates to $0$, this way, we finally reached the aimed point.

For $0\in\uuu^{even}$ the setting is extremely similar to the well-known Morse case and the proof is exactly the same.

Recalling that for the sets we considered throughout the proof (open subsets of the Euclidean space), being connected and arc-connected is equivalent, we have proven the proposition.

\ep

\subsection{Labeling in the generic case}\label{sec:labeling}

Now, among the saddle points of $W$ near a given minima $\m$, not all are important for the study of the eigenvalue associated with $\m$. Heuristically, the crucial ones are those of minimum height (that means they minimize $W(\s)$) such that a process stuck around $\m$ needs to cross in order to fall into a well of lower energy (that is, a well associated with some $\m'$ such that $W(\m') < W(\m)$).

\begin{defin}\label{def:labeling}
    We say that $\s\in\uuu$ is a separating saddle point if, for every $r>0$ small enough, $\{x\in B(\s,r),\ W(x)<W(\s)\}$ is composed of two connected components that are contained in two different connected components of $\{x\in\R^n,\ W(x)<W(\s)\}$. Hence, from Proposition \ref{prop:saddle} we have that the set of these points is a subset of $\uuu^{(1)}$, we will then denote it by $\vvv^{(1)}$, we also denote by separating saddle value of $W$ a point in $V(\vvv^{(1)})$.

    We say that $E\subset\R^n$ is a critical component of $W$ if there exists $\sigma\in V(\vvv^{(1}))$ such that $E$ is a connected component of $\{W\leq\sigma\}$ and $\p\negthinspace E\cup \vvv^{(1)}\neq\emptyset$.
\end{defin}

\begin{figure}[h]
    \centering
    \includegraphics[width=\linewidth]{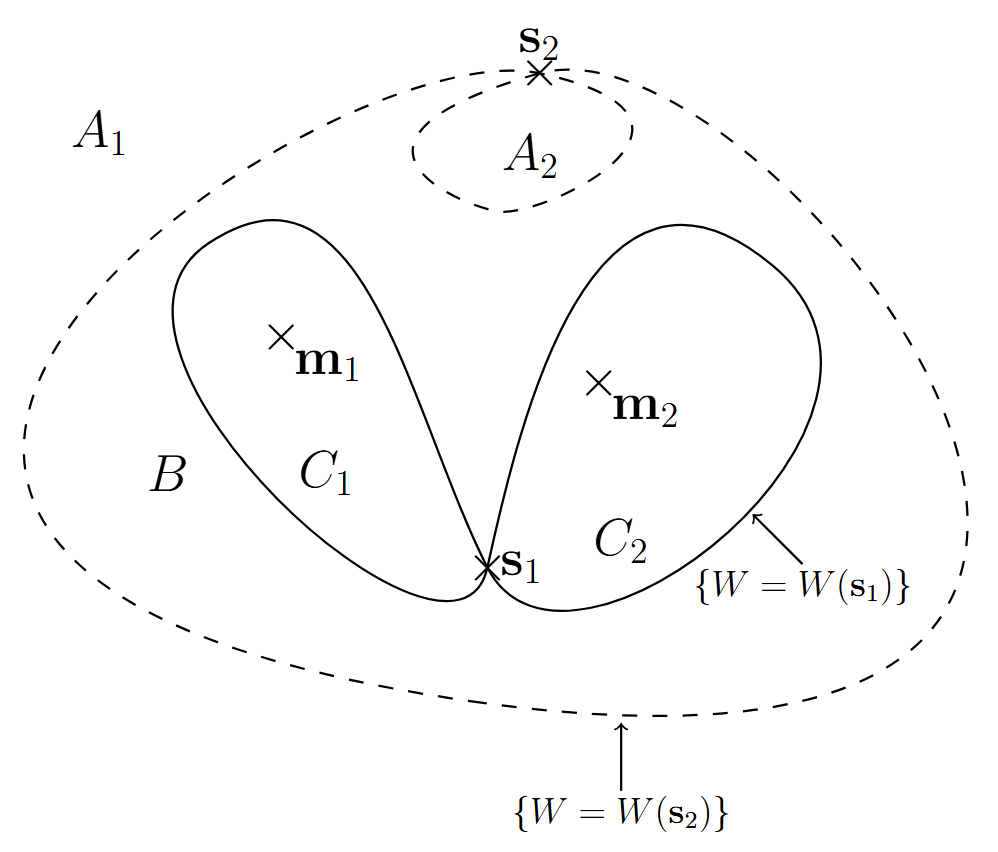}
    \caption{Example of a separating and a non-separating saddle point}
    \label{fig:separating}
\end{figure}

Consider a potential $W$ having the level sets of Figure \ref{fig:separating}, with $\m_1,\m_2$ minima, $\s_1,\s_2$ saddle points such that
\bes
W(\m_1) < W(\m_2) < W(\s_1) < W(\s_2).
\ees
And the capital letters denote connected components of some subsets the following way
\bes
A_1\sqcup A_2 = \{W > W(\s_2)\},\ B = \{W(\s_1) < W < W(\s_2)\},\ C_1\sqcup C_2 = \{W < W(\s_1)\}.
\ees
Observe that $\s_1$ is a separating saddle point but not $\s_2$. A process starting near $\m_2$ can go to a well of lower energy by crossing $\s_1$ and going to $\m_1$. However, a process starting near $\m_1$ cannot reach below $\m_1$ by crossing $\s_2$ because both connected components of $\{x\in B(\s_2,r),\ W(x)<W(\s_2)\}$ are in the same connected component of $\{x\in \R^n,\ W(x)<W(\s_2)\}$, and so the process would just go back to $\m_1$.

Let us now recall the labeling procedure. Under the assumptions \ref{ass.confin} and \ref{ass:morsedegen2}, we have that $W(\vvv^{(1}))$ is finite. We denote $N = \sharp W(\vvv^{(1})) + 1$ and by $\sigma_2 > \sigma_3 > \ldots > \sigma_N$ the elements of $W(\vvv^{(1}))$, for convenience, we also introduce a fictive infinite saddle value $\sigma_1 = +\infty$. Starting from $\sigma_1$, we recursively associate to each $\sigma_i$ a finite family of local minima $(\m_{i,j})_j$ and a finite family of critical component $(E_{i,j})_j$:
\begin{itemize}
    \item Let $X_{\sigma_1} = \{x\in\R^n,\ W(x) < \sigma_1\} = \R^n$. We let $\m_{1,1}$ be any global minimum of $V$ and $E_{1,1} = \R^n$. In the following we will write $\um = \m_{1,1}$.
    \item Next, we consider $X_{\sigma_2} = \{x\in \R^n,\ W(x) < \sigma_2\}$. This is the union of its finitely many connected components. Exactly one contains $\um$ and the other components are denoted by $E_{2,1},\ldots,E_{2,N_2}$. They are all critical and, in each component $E_{2,j}$ we pick up a point $\m_{2,j}$ which is a global minimum of $W_{|E_{2,j}}$.
    \item Suppose now that the families $(\m_{k,j})_j$ and $(E_{k,j})_j$ have been constructed until rank $k = i-1$. The set $X_{\sigma_i} = \{x\in\R^n,\ W(x) < \sigma_i\}$ has again finitely many connected components and we label $E_{i,j}$, $j\in\lb1,N_i\rb$, those of these components that do not contain any $\m_{k,j}$ for $k<i$. They are all critical and, in each $E_{i,j}$, we pick up a point $\m_{i,j}$ which is a global minimum of $W_{|E_{i,j}}$.
\end{itemize}
At the end of this procedure, all the minima have been labeled. Throughout, we denote by $\s_1$ a fictive saddle point such that $W(\s_1) = \sigma_1 = \infty$ and, for any set $A$, $\ppp(A)$ denotes the power set of $A$. From the above labeling, we define two mappings
\bes
E:\uuu^{(0)}\to\ppp(\R^n)\ \ \mbox{ and }\ \ \j:\uuu^{(0)}\to\ppp(\vvv^{(1)}\cup\{\s_1\})
\ees
as follows: for every $i\in\lb1,N\rb$ and every $j\in\lb1,N_i\rb$,
\bes
E(\m_{i,j}) = E_{i,j},
\ees
and
\bes
\j(\um) = \{\s_1\}\ \ \mbox{ and }\ \ \j(\m_{i,j}) = \p\negthinspace E_{i,j}\cap \vvv^{(1)}\ \mbox{ for }i\geq2.
\ees
In particular, we have $E(\um) = \R^n$ and, for all $i\in\lb1,N\rb$, $j\in\lb1,N_i\rb$, one has $\emptyset\neq\j(\m_{i,j})\subset\{W=\sigma_i\}$. We then define the mappings
\bes
\bsigma:\uuu^{(0)}\to W(\vvv^{(1)})\cup\{\sigma_1\}\ \ \mbox{ and }\ \ S:\uuu^{(0)}\to(0,+\infty],
\ees
by
\be\label{eq:defS}
\forall \m\in\uuu^{(0)},\ \ \bsigma(\m) = W(\j(\m))\ \ \mbox{ and }\ \ S(\m) = \bsigma(\m) - W(\m),
\ee
where, with a slight abuse of notation, we have identified $W(\j(\m))$ with its unique element. Note that $S(\m) = \infty$ if and only if $\m = \um$. 

We now consider the following generic assumption in order to lighten the result and the proof.
\be\label{ass.gener2}\tag{Gener}
\begin{array}{l}
    (\ast)\mbox{ for any }\m\in\uuu^{(0)},\m \mbox{ is the unique global minimum of } W_{|E(\m)}\\
    (\ast)\mbox{ for all }\m\neq\m'\in\uuu^{(0)},\j(\m)\cap\j(\m')=\emptyset.
\end{array}
\ee
In particular, \eqref{ass.gener2} implies that $W$ uniquely attains its global minimum at $\um$. This assumption allows us to avoid some heavy constructions regarding the set $\uuu$ and lighten the definition of the quasimodes, see \cite{Mi19} for the general setting. But it seems that its not a true obstruction and that we can pursue the computations without this assumption as described in \cite[Section 6]{BoLePMi22} and \cite{No24}, in the spirit of \cite{Mi19}.

This assumption comes from \cite[(1.7)]{BoGaKl05_01} and \cite[Assumption 3.8]{HeKlNi04_01} where they were hard to handle properly. Then it changed to become \cite[Hypothesis 5.1]{HeHiSj11_01} when finally reaching the form of \cite[Assumption 4]{LePMi20}.

One can show that \eqref{ass.gener2} is weaker than \cite{BoGaKl05_01}'s, \cite{HeKlNi04_01}'s and \cite{HeHiSj11_01}'s assumptions. More precisely, they supposed that $\j(\m)$ is a singleton while here we have no restriction on the size of $\j(\m)$.

\section{Some technical results}

For the four following lemmas, by $Q = O(r(h))$ (for $Q$ operator and $r$ some positive real function) we mean there exists $C>0$ such that $\vvert{Q} \leq Cr(h)$.

\begin{lemma}\label{lem:1}
For $s\geq0$
\bes
(h^{2-\frac2{\overline{\nu}}}g_1(h) + \d_V^*G\d_V)^{-s} = O((h^{2-\frac2{\overline{\nu}}}g_1(h))^{-s}).
\ees
\end{lemma}
\bp We use that $\d_V^*G\d_V\geq0$ and functional calculus.

\ep

\begin{lemma}\label{lem:2}
    $\alpha\cdot \d_V\Pi(h^{2-\frac2{\overline{\nu}}}g_1(h) + \d_V^*G\d_V)^{-1/2} = O(1)$.
\end{lemma}
\bp
\bes
\begin{split}
    (\ast_1) :&= \vvert{\alpha\cdot \d_V\Pi(h^{2-\frac2{\overline{\nu}}}g_1(h) + \d_V^*G\d_V)^{-1/2}u}^2\\
    &= \<(\alpha\cdot \d_V\Pi)^*\alpha\cdot \d_V\Pi(h^{2-\frac2{\overline{\nu}}}g_1(h) + \d_V^*G\d_V)^{-1/2}u,(h^{2-\frac2{\overline{\nu}}}g_1(h) + \d_V^*G\d_V)^{-1/2}u\>\\
    &= \<\d_V^*G\d_V\Pi(h^{2-\frac2{\overline{\nu}}}g_1(h) + \d_V^*G\d_V)^{-1/2}u,(h^{2-\frac2{\overline{\nu}}}g_1(h) + \d_V^*G\d_V)^{-1/2}u\>\\
    &\leq \vvert{\d_V^*G\d_V(h^{2-\frac2{\overline{\nu}}}g_1(h) + \d_V^*G\d_V)^{-1}}\vvert{u}^2,
\end{split}
\ees
having that $\d_V^*G\d_V$ commutes with $(h^{2-\frac2{\overline{\nu}}}g_1(h) + \d_V^*G\d_V)^{-1/2}$ from functional calculus, hence the result.

\ep

In this Lemma we can replace $\alpha\Pi$ by $G^{1/2}$ and obtain the same result.

\begin{lemma}\label{lem:3}
    $\d_V(h^{2-\frac2{\overline{\nu}}}g_1(h) + \d_V^*G\d_V)^{-1/2} = O(g_1(h)^{-1/2})$.
\end{lemma}
\bp
\bes
\begin{split}
    (\ast_2) :&= \vvert{\d_V(h^{2-\frac2{\overline{\nu}}}g_1(h) + \d_V^*G\d_V)^{-1/2}u}^2\\
    &= \vvert{G^{-1/2}G^{1/2}\d_V(h^{2-\frac2{\overline{\nu}}}g_1(h) + \d_V^*G\d_V)^{-1/2}u}^2\\
    &= \<G^{-1}G^{1/2}\d_V(h^{2-\frac2{\overline{\nu}}}g_1(h) + \d_V^*G\d_V)^{-1/2}u,G^{1/2}\d_V(h^{2-\frac2{\overline{\nu}}}g_1(h) + \d_V^*G\d_V)^{-1/2}u\>\\
    &\leq g_1(h)^{-1}\vvert{G^{-1/2}\d_V(h^{2-\frac2{\overline{\nu}}}g_1(h) + \d_V^*G\d_V)^{-1/2}u}^2\\
    &\leq g_1(h)^{-1}\vvert{u}^2
\end{split}
\ees
using Assumption \ref{ass.G}, and a result similar to Lemma \ref{lem:2}.

\ep

\begin{lemma}\label{lem:4}
    For $i\in\lb1,d\rb$, denoting $\d_{V,i} = h\p_{x_i} + \p_iV$ so that $\d_V = (\d_{V,i})_i$, we have for all $i,j\in\lb1,d\rb$
    \bes
    \d_{V,i}^*\d_{V,j}(h^{2-\frac2{\overline{\nu}}}g_1(h) + \d_V^*G\d_V)^{-1} = O(h^{\frac{2}{\overline{\nu}} - 1}g_1(h)^{-3/2}g_2(h)^{1/2}).
    \ees
\end{lemma}
\bp
In the following, we denote $R = (h^{2-\frac2{\overline{\nu}}}g_1(h) + \d_V^*G\d_V)^{-1}$ for clarity.

Using the commutation rules
\bes
[\d_{V,i}^*,\d_{V,j}] = -2h\p_{i,j}^2V;\quad [\d_{V,i}^*,R] = R[\d_V^*G\d_V,\d_{V,i}^*]R
\ees
we have
\bes
\begin{split}
    \d_{V,i}^*\d_{V,j}R &= [\d_{V,i}^*,\d_{V,j}]R + \d_{V,j}\d_{V,i}^*R\\
    &= -2h\p_{i,j}^2VR + \d_{V,j}R\d_{V,i}^* + \d_{V,j}[\d_{V,i}^*,R]\\
    &= -2h\p_{i,j}^2VR + \d_{V,j}R\d_{V,i}^* + \d_{V,j}R[\d_V^*G\d_V,\d_{V,i}^*]R.
\end{split}
\ees
The last commutator gives
\bes
\begin{split}
    [\d_V^*G\d_V,\d_{V,i}^*] &= [\d_V^*,\d_{V,i}^*]G\d_V + \d_V^*[G,\d_{V,i}^*]\d_V + \d_V^*G[\d_V,\d_{V,i}^*]\\
    &= \d_V^*h\p_iG\d_V + 2h\d_V^*G(\p_{i,k}^2V)_{k}
\end{split}
\ees
since $[\d_{V,i}^*,\d_{V,j}^*] = 0$ for all $i,j$. Using Assumptions \ref{ass.confin} and \ref{ass.G} $ii)$ we have
\bes
\begin{split}
    \vvert{\d_{V,i}^*\d_{V,j}R} &\leq 2Ch\vvert{R} + \vvert{\d_{V,j}R^{1/2}}\vvert{R^{1/2}\d_{V,i}^*}\\
    &\phantom{*****} + Ch\vvert{\d_{V,j}R^{1/2}}\vvert{R^{1/2}\d_V^*G^{1/2}}\vvert{G^{1/2}\d_VR^{1/2}}\vvert{R^{1/2}}\\
    &\phantom{*****} + 2Ch\vvert{\d_{V,j}R^{1/2}}\vvert{R^{1/2}\d_V^*G^{1/2}}\vvert{G^{1/2}}\vvert{R}.
\end{split}
\ees
Here we used that
\bes
\begin{split}
    \vvert{R^{1/2}\d_V^*\p_iG\d_VR^{1/2}} &= \sup_{u\neq0}\frac{\<R^{1/2}\d_V^*\p_iG\d_VR^{1/2}u,u\>}{\vvert{u^2}}\\
    &= \sup_{u\neq0}\frac{\<\p_iG\d_VR^{1/2}u,\d_VR^{1/2}u\>}{\vvert{u^2}}\\
    &\lesssim \sup_{u\neq0}\frac{\<G\d_VR^{1/2}u,\d_VR^{1/2}u\>}{\vvert{u^2}}\\
    &= \vvert{R^{1/2}\d_V^*G\d_VR^{1/2}}\\
    &\leq \vvert{R^{1/2}\d_V^*G^{1/2}}\vvert{G^{1/2}\d_VR^{1/2}}
\end{split}
\ees

Now with Lemmas \ref{lem:2} and \ref{lem:3}, we obtain
\bes
\vvert{\d_{V,i}^*\d_{V,j}R} \lesssim h^{\frac{2}{\overline{\nu}} - 1}g_1(h)^{-1} + g_1(h)^{-1} + h^{\frac{1}{\overline{\nu}}}g_1(h)^{-1} + h^{\frac{2}{\overline{\nu}} - 1}g_1(h)^{-3/2}g_2(h)^{1/2}.
\ees
hence the result.

\ep

Therefore, we also have
\bes
\Delta_V(h^{2-\frac2{\overline{\nu}}}g_1(h) + \d_V^*G\d_V)^{-1} = O(h^{\frac{2}{\overline{\nu}} - 1}g_1(h)^{-3/2}g_2(h)^{1/2}).
\ees

\subsection{Proof of Lemma \ref{lem:soussolaffine}}\label{sec:proofsoussolaffine}
Up to a translation, let us consider that $V\geq 0$. According to Assumption \ref{ass.confin}, we also have that there exists $R>0$ such that for all $|x|\geq R$, $|\nabla V(x)| \geq \frac1C$.

For $x_0\in\R^d$, we consider the maximal solution to the Cauchy problem
\be\label{1eq:edo}
\left\{
\begin{aligned}
&\dot x(t) = -\nabla V(x(t)),
\\&x(0) = x_0.
\end{aligned}
\right.
\ee
We then see that
\bes
\frac{d}{dt}(V(x(t))) = -|\nabla V(x(t))|^2.
\ees

Let us first show that for $\lambda>0$, $x_0\in\R^d$ such that $V(x_0) \leq\lambda$ and $|x_0| > R + \lambda C$, we have a solution to \eqref{1eq:edo} defined on $\R_+$ and for all $t\geq0$, $|x(t)|\geq R$. Consider $t$ the supremum such that this is true, and by contradiction assume that $t<\infty$. Then
\bes
\begin{aligned}
    |x(t) - x(0)| &= \big|\int_0^t-\nabla V(x(s))ds\big|\ \leq \int_0^t|\nabla V(x(s))|ds\\
    &\leq C\int_0^t|\nabla V(x(s))|^2ds = C(V(x_0) - V(x(t)) \leq \lambda C,
\end{aligned}
\ees
which leads to $|x(t)|\leq |x_0| + \lambda C$. But we also have that
\bes
|x(t)| \geq |x_0| - |x(t) - x_0| > R + \lambda C - \lambda C = R
\ees
this is a contradiction with the maximality of $t$.

Now assume that such an $x_0$ exists, therefore we have
\bes
0\leq V(x(t)) = V(x_0) + \int_0^td(V(x(s))) = V(x_0) - \int_0^t|\nabla V(x(s))|^2ds \leq \lambda - \frac1{C^2}t,
\ees
but because $x(t)$ is global, this is absurd for $t$ big enough. This leads to
\be\label{1eq:sublevel}
\{V\leq\lambda\}\subset \overline B(0,R+\lambda C)
\ee
where the right-hand side term is the closed ball centered at $0$ of radius $R+\lambda C$. Now let $x\in\R^d$ be such that $|x| > R$, therefore there exists $\eps>0$ such that $\lambda = \frac1C(|x| - R - \eps) > 0$. Hence we can write $|x| = R + \lambda C + \eps$, thus by \eqref{1eq:sublevel}, we know that $V(x) > \lambda = \frac1C|x| + a$, with $a = -\frac{R+\eps}{C}$. Therefore, there exists $\tilde a\in\R$ such that
\be\label{1eq:borneinfV}
\forall |x| > R,\ \ V(x) > \frac1C|x| + \tilde a.
\ee
Because $|x|\leq R$ is compact, we can find $b\in\R$ such that \eqref{1eq:borneinfV} is true on the whole space.

\subsection{Proof of Proposition \ref{prop:existence}}\label{ssec:existence}

We use \cite[Theorem 1.7]{Kh12} which states that a Lyapunov function for the associated deterministic differential equation ensures the result announced in the proposition. The system we have to study is
\bes
\left\{
\begin{aligned}
&dx_t = \alpha(x_t,v_t)dt,\\
&dv_t = \beta(x_t,v_t)dt - 4\Sigma^T\Sigma v_tdt.
\end{aligned}
\right.
\ees
Consider the function $t\mapsto f(x_t,v_t)$,
\bes
\frac{d(f(x_t,v_t))}{dt} = \alpha\cdot\p_xV + 2\beta\cdot\Sigma^T\Sigma v_t - 8|\Sigma^T\Sigma v_t|^2 \leq \frac h2 (\div_x\alpha + \div_v\beta) \lesssim f(x_t,v_t).
\ees
Here we used \eqref{eq:assskew} and Assumption \ref{ass.accretive}. We define
\bes
L:t\mapsto\sqrt{f(x_t,v_t) - \min V + 1}.
\ees
It is clear that this is a Lyapunov function satisfying $L' \lesssim L$. Using Assumption \ref{ass.confin}, it also satisfies the other required hypotheses to apply \cite[Theorem 1.7]{Kh12}, namely: $\lim_{R\to\infty}\inf_{|x,v|\geq R}L = \infty$ and $L$ is globally Lipschitz. For the first one, it is a direct consequence of Lemma \ref{lem:soussolaffine}, and the global Lipschitz property is induced by the boundedness of the Hessian of $V$.

\ep

\subsection{Proof of Proposition \ref{prop:accretive}}\label{ssec:accretive}

The idea is to mimic the proof of \cite[Theorem 15.1]{He13}.

Let $h>0$ be fixed. To show that $P$ admits a maximal accretive extension, it is first necessary to show that it admits an accretive extension, this comes from the skew-adjointness of $X$, as well as from the non-negativity of $N$. It therefore remains to show the maximal side, for that we use the criterion which tells us that $P$ is maximal accretive if $T=P+(2h\Tr(\Sigma^T\Sigma)+1)\Id$ has a dense image.

Let $u\in L^2(\R^{d+d'})$ such that
\be\label{eq:uTphii}
\forall \phii\in\ccc_c^\infty(\R^{d+d'}),\ \< u,T\phii\> =0.
\ee
We then must show that $u=0$. As $P$ is real, we can assume also is $u$. We split $X= \underbrace{\alpha\cdot h\p_x + \beta\cdot h\p_v}_{=X_1}\underbrace{+\frac h2(\div_x\alpha+\div_v\beta)}_{=X_0}$ where we can see $X_1$ is a homogeneous differential operator of order 1 and $X_0$ is a mere $\ccc^\infty$ function.

Under Assumption \ref{ass.Hormander}, $T^*$ is hypoelliptic and using that thanks to \eqref{eq:uTphii}, $T^*u=0$ in $\ddd'$, we have that $u\in\ccc^\infty(\R^{d+d'})$. We then consider $\zeta$ a $\ccc^\infty_c(\R)$ function satisfying $0\leq\zeta\leq1$, $\zeta = 1$ on $[-1,1]$ and $\zeta=0$ outside of $[-2,2]$ and we denote for $k\in\R_+^*$, $\zeta_k(x,v) = \zeta\Big(\frac{f(x,v)}{k}\Big)$ where $f$ is defined in \eqref{eq:f}. We thus have
\bes
\begin{split}
    \<u,T\zeta_k^2 u\> &= \<u,X\zeta_k^2 u\> + \<\zeta_k u,(N+(2h\Tr(\Sigma^T\Sigma)+1)\Id)\zeta_k u\> + \<u, [N,\zeta_k]\zeta_k u\>
\end{split}
\ees
\bes
\begin{split}
    \<u, [N,\zeta_k]\zeta_k u\> &= -h^2\<\p_v\zeta_k,2u\p_v(\zeta_k u)-\p_v(u\zeta_k u)\>\\
    &= -h^2\<\p_v\zeta_k,u^2\p_v\zeta_k+2\zeta_k(u\p_v u-\p_vu\ u)\>\\
    &= -h^2\vvert{u\p_v\zeta_k}^2
\end{split}
\ees
\bes
\begin{split}
    \<\zeta_k u,(N+(2h\Tr(\Sigma^T\Sigma)+1)\Id)\zeta_k u\> &= h^2\vvert{\p_v(\zeta_k u)}^2 + \vvert{2|\Sigma^T\Sigma v|\zeta_k u}^2 + \vvert{\zeta_k u}^2
\end{split}
\ees
\bes
\begin{split}
    \<u,X\zeta_k^2 u\> &= \<\zeta_k u,X\zeta_k u\> + \<u,[X,\zeta_k]\zeta_k u\>\\
    &= \<u,[X_1,\zeta_k]\zeta_k u\>\\
    &= \frac{h^2}{2k}\<u,(\div_x\alpha + \div_v\beta)\zeta'(f/k)\zeta_k u\>
\end{split}
\ees
where we used Assumption \ref{ass.skewadj} for the last estimate and the skew-adjointness of $X$ the line above. Using \eqref{eq:uTphii} with $\phii = \zeta_k^2u$, we obtain
\bes
\begin{split}
    (\ast):&= h^2\vvert{\p_v(\zeta_k u)}^2 + \vvert{2|\Sigma^T\Sigma v|\zeta_k u}^2 + \vvert{\zeta_k u}^2\\
    &= h^2\vvert{u\p_v\zeta_k}^2 - \frac{h^2}{2k}\<u,(\div_x\alpha + \div_v\beta)\zeta'(f/k)\zeta_k u\>.
\end{split}
\ees
Using that $\p_v\zeta_k = \frac2k\Sigma^T\Sigma v\zeta'(f/k)$ and
\bes
h^2\vvert{\p_v(\zeta_k u)}^2 + \vvert{2|\Sigma^T\Sigma v|\zeta_k u}^2-\frac{h^2}{k^2}\vvert{2|\Sigma^T\Sigma v|\zeta'(f/k)u}^2\geq0
\ees
for $k$ large enough, we have
\bes
\vvert{\zeta_k u}^2 \leq - \frac{h^2}{2k}\<u,(\div_x\alpha + \div_v\beta)\zeta'(f/k)\zeta_k u\>.
\ees
Therefore, using that $\supp\zeta' \cap (-1,1) = \emptyset$, we have $\frac1k(\div_x\alpha + \div_v\beta)\zeta'(f/k)\to0$. Hence, with the dominated convergence theorem, $u=0$ because $\zeta_ku\to u$.

\ep

\subsection{Proof of Lemma \ref{lem:situ2}}\label{ssec:lemma}

To this purpose, we define a valuation $\gamma$ over the monomials in the variables $h,x$ and $v$ by: $\gamma(ch^jx^av^b) = 2j+a+b$ for all $c\in\R^*$. We now prove by induction the following property
\bes\label{HR}\tag{$\Pr_n$}
\text{The terms of valuation }n\text{ in }\ell\text{ are even in }v.
\ees
We start at $n=1$, it has been proven when looking at the principal order of the eikonal equation that assuming $\xi_x \neq 0$ implies $\xi_v = 0$. Let $n\geq2$, assume $\Pr_m$ is true for all $m<n$, then look at \eqref{eq:decompjab} for $2j+a+b = n$ and $b=2b'+1$ is odd. We know that $\p_v\ell = O(v+h)$ and $\ell = O(x + v^2 + h)$. Thus, since the valuation $\gamma$ is multiplicative ($\gamma(uv) = \gamma(u)\gamma(v)$), all the factors of terms appearing in $p_{j,a,2b'+1}$ come from monomials of valuation less than $n$, hence they are even in $v$ by the hypothesis. Then, $\p_v\ell$ will give odd terms and $\ell$ even ones, therefore $p_{j,a,2b'+1}$ is even in $v$, but it is of the form $h^jx^av^{2b'+1}$ thus $p_{j,a,2b'+1} = 0$. Moreover, for all $k\in\lb1,a\rb$,
\bes
\gamma(\ell_{j,a-k,2b'+1}) \leq \gamma(\ell_{j,a+1,2b'-1}) = \gamma(\ell_{j-1,a+1,2b'+1}) = n-1.
\ees
By the induction hypothesis, the monomials $\ell_{j,a-k,2b'+1}, \ell_{j,a+1,2b'-1}, \ell_{j-1,a+1,2b'+1}$ are $0$ and we obtain
\be\label{eq:vdvmoinsdv2}
v\p_v\ell_{j,a,2b'+1} - \p_v^2\ell_{j-1,a,2b'+3} = 0
\ee
for all $j,a,b'\in\N$ such that $2j+a+2b'+1 = n$. Therefore taking $j=0$ we have
\bes
\forall a,b'\in\N,a+2b'+1=n,\quad v\p_v\ell_{0,a,2b'+1} = 0.
\ees
Using now a quick induction on $j$ along with \eqref{eq:vdvmoinsdv2} we complete the induction on $n$.

\ep

\end{document}